%% file: ms-2024nov13.tex
\documentclass[11pt]{article} 

\usepackage[utf8]{inputenc} 

\usepackage{geometry} 
\usepackage{graphicx} 

\usepackage{booktabs} 
\usepackage{array} 
\usepackage{paralist} 
\usepackage{verbatim} 
\usepackage{subfig} 

\usepackage{fancyhdr} 
\usepackage{sectsty}
\allsectionsfont{\sffamily\mdseries\upshape} 

\usepackage[nottoc,notlof,notlot]{tocbibind} 
\usepackage[titles,subfigure]{tocloft} 

\usepackage{amsmath, amssymb, bm}
\usepackage[T1]{fontenc}
\usepackage{textcomp}

\newtheorem{theorem}{Theorem}[section]
\newtheorem{lemma}[theorem]{Lemma}
\newtheorem{proposition}[theorem]{Proposition}

\newtheorem{definition}[theorem]{Definition}

\newtheorem{corollary}[theorem]{Corollary}

\numberwithin{equation}{section} 

\title{Dynamical Spin Limit Shape of Young Diagram and \\ 
Spin Jucys--Murphy Elements for Symmetric Groups
\thanks{This work was supported by JSPS KAKENHI Grant Number 
JP22K03346.}}
\author{{Akihito HORA}
\thanks{Department of Mathematics, Faculty of Science, Hokkaido University, Sapporo 
060-0810, Japan; hora@math.sci.hokudai.ac.jp}}
\date{Dedicated to Professor Takeshi Hirai for his 88th birthday}

\begin{document}
\maketitle

\begin{abstract}
The branching rule for spin irreducible representations of symmetric groups gives 
rise to a Markov chain on the spin dual $(\widetilde{\mathfrak{S}}_n)^\wedge_{\mathrm{spin}}$ 
of symmetric group $\mathfrak{S}_n$ through restriction and induction of 
spin irreducible representations. 
This further produces a continuous time random walk $(X_s^{(n)})_{s\geqq 0}$ on 
$(\widetilde{\mathfrak{S}}_n)^\wedge_{\mathrm{spin}}$ by introducing an appropriate 
pausing time. 
Taking diffusive scaling limit for these random walks under $s=tn$ and $1/\sqrt{n}$ 
reduction as $n\to\infty$, we consider 
a concentration phenomenon at each macroscopic time $t$. 
Since an element of $(\widetilde{\mathfrak{S}}_n)^\wedge_{\mathrm{spin}}$ is 
regarded as a strict partition of $n$ with $\pm 1$ indices, the limit shapes of 
profiles of strict partitions appear. 
In this paper, we give a framework in which initial concentration at $t=0$ is 
propagated to concentration at any $t>0$. 
We thus obtain the limit shape $\omega_t$ depending on macroscopic time $t$, 
and describe the time evolution by using devices in free probability theory. 
Included is the case where Kerov's transition measure has non-compact support 
but determinate moment problem. 
A spin version of Biane's formula for spin Jucys--Murphy elements is shown, 
which plays an important role in our analysis. 
\end{abstract}

\section{Introduction}

Young diagrams are ubiquitous objects appearing in a wide variety of scenes in mathematics. 
In this paper, we consider statistical properties of an ensemble consisting of a huge number 
of Young diagrams for both static and dynamic models. 
Let $\mathcal{P}_n$ denote the set of partitions of a natural number $n$. 
For example, 
\[
\mathcal{P}_4 = \bigl\{ (4), (3,1), (2,2), (2,1,1), (1,1,1,1) \bigr\}. 
\] 
The set of partitions of natural numbers is 
\[ 
\mathcal{P} = \bigcup_{n=0}^\infty \mathcal{P}_n, \qquad \mathcal{P}_0 = \{\varnothing\}. 
\] 
\begin{figure}[hbt]
\centering
\include{Ylambda-fig}
\vspace{-5mm}
\caption{$Y(\lambda)$ for $\lambda = (3,1)\in \mathcal{P}_4$}
\label{fig:1-1}
\end{figure}
The size of $\lambda\in\mathcal{P}$ is denoted by $|\lambda|$, that is, 
$\lambda\in\mathcal{P}_n \Leftrightarrow |\lambda| =n$. 
The number of parts of $\lambda$ is denoted by $l(\lambda)$. 
The set of equivalence classes of irreducible representations of $\mathfrak{S}_n$, 
the symmetric group of degree $n$, is parametrized by $\mathcal{P}_n$. 
Let $m_k(\lambda)$ denote the number of parts of length $k$ in $\lambda\in \mathcal{P}_n$. 
Then $\lambda$ is determined by $(m_k(\lambda))_{k=1}^\infty$: 
$\lambda = (1^{m_1(\lambda)} 2^{m_2(\lambda)} \cdots)$. 
Young diagram $Y(\lambda)$ is assiged to $\lambda\in \mathcal{P}_n$ as in Figure~\ref{fig:1-1}. 
Let $\mathcal{SP}_n$ denote the subset of $\mathcal{P}_n$ consisting of the strict partitions. 
For example, \ $\mathcal{SP}_4 = \bigl\{ (4), (3,1)\bigr\}$. 
Set  
\[ 
\mathcal{SP} = \bigcup_{n=0}^\infty \mathcal{SP}_n, \qquad \mathcal{SP}_0 = \{\varnothing\}. 
\] 
Namely, 
$\lambda \in \mathcal{SP} \iff m_k(\lambda) \leqq 1 \ (\forall k\in\mathbb{N})$. 
Assigned to $\lambda\in \mathcal{SP}_n$ is a shifted Young diagram $S(\lambda)$ where 
rows are shifted one by one as in Figure~\ref{fig:1-2}. 
In this paper, since we treat asymptotic properties of spin representations of 
the symmetric group $\mathfrak{S}_n$ as $n$ tends to $\infty$, we are interested in 
the growth of the shape of $S(\lambda)$. 
\begin{figure}[hbt]
\centering
\include{Slambda-fig}
\vspace{-5mm}
\caption{$S(\lambda)$ for $\lambda = (3,1)\in \mathcal{SP}_4$}
\label{fig:1-2}
\end{figure}

Limit shapes of growing Young diagrams have been studied by many authors. 
Let $M_{\mathrm{U}}^{(n)}$ be the uniform probability on $\mathcal{P}_n$: 
\[ 
M_{\mathrm{U}}^{(n)}(\{\lambda\}) = \frac{1}{p(n)} \quad (\lambda\in\mathcal{P}_n), 
\qquad p(n) = |\mathcal{P}_n|.
\] 
Similarly, $M_{\mathrm{US}}^{(n)}$ denotes the uniform probability on $\mathcal{SP}_n$. 
For $0<\delta <1$, let $M_{\mathrm{U}}^\delta$ and $M_{\mathrm{US}}^\delta$ 
be the probabilities on $\mathcal{P}$ and $\mathcal{SP}$ 
defined by 
\begin{align*}
&M_{\mathrm{U}}^\delta (\{\lambda\}) = \frac{1}{C} \,\delta^{|\lambda|} \quad 
(\lambda\in\mathcal{P}), \qquad 
C= \sum_{\lambda\in\mathcal{P}} \delta^{|\lambda|} = \sum_{n=0}^\infty p(n) \delta^n, \\ 
&M_{\mathrm{US}}^\delta (\{\lambda\}) = 
\frac{1}{C^\prime} \,\delta^{|\lambda|} \quad 
(\lambda\in\mathcal{SP}), \qquad 
C^\prime = \sum_{\lambda\in\mathcal{SP}} \delta^{|\lambda|}. 
\end{align*}
In accordance with terminology in statistical mechanics, probability spaces like 
$(\mathcal{P}_n, M_{\mathrm{U}}^{(n)})$ and $(\mathcal{P}, M_{\mathrm{U}}^\delta)$ 
are often referred to as micro-canonical and grand-canonical respectively. 
\begin{center}
\begin{tabular}{ll} 
\toprule 
micro-canonical  \quad & $|\lambda|$ is fixed \\ 
canonical & $l(\lambda)$ is limited \\ 
grand-canonical  & no limitation  \\ \bottomrule 
\end{tabular}
\end{center}
Put $Y(\lambda)$ for $\lambda\in \mathcal{P}_n$ on the first quadrant of the $xy$ plane as 
indicated in the middle of Figure~\ref{fig:1-1}. 
The profile of $Y(\lambda)$ is regarded as the graph of $y=Y(\lambda)(x)$. 
For $\lambda\in\mathcal{P}$ and $r>0$, 
let $Y(\lambda)^r$ denote $Y(\lambda)$ reduced by $1/r$: 
\begin{equation}\label{eq:1-1}
Y(\lambda)^r (x) = 
\frac{1}{r} \, Y(\lambda) \bigl( r x\bigr), \qquad 
\lambda\in \mathcal{P}.
\end{equation}
Let $N>0$, $\delta(N)$ and $\delta^\prime (N)$ be related to satisfy 
\[ 
\sum_{\lambda\in\mathcal{P}} |\lambda| M_{\mathrm{U}}^{\delta(N)}(\{\lambda\}) = 
\sum_{\lambda\in\mathcal{SP}} |\lambda| M_{\mathrm{US}}^{\delta^\prime (N)}(\{\lambda\}) 
=N. 
\] 
We have $\delta(N) \to 1$ and $\delta^\prime (N) \to 1$ as $N\to\infty$. 
The following laws of large numbers \eqref{eq:1-2}, \eqref{eq:1-3} 
for $Y(\lambda)^{\sqrt{n}}$ and \eqref{eq:1-2.5}, \eqref{eq:1-3.5} for $Y(\lambda)^{\sqrt{N}}$ 
are due to Vershik (\cite{Ver96}): \ 
for $\forall\epsilon >0$, 
\begin{align}
&M_{\mathrm{U}}^{(n)} \Bigl( \Bigl\{ \lambda\in \mathcal{P}_n \,\Big|\, 
\sup_{x>0} \bigl| Y(\lambda)^{\sqrt{n}} (x)- \psi_{\mathrm{U}} (x)\bigr| 
> \epsilon \Bigr\} \Bigr) \ 
\xrightarrow[\,n\to\infty\,] \ 0, 
\label{eq:1-2} \\ 
&M_{\mathrm{U}}^{\delta(N)} \Bigl( \Bigl\{ \lambda\in \mathcal{P} \,\Big|\, 
\sup_{x>0} \bigl| Y(\lambda)^{\sqrt{N}} (x)- \psi_{\mathrm{U}} (x)\bigr| 
> \epsilon \Bigr\} \Bigr) \ 
\xrightarrow[\,n\to\infty\,] \ 0, 
\label{eq:1-2.5} \\ 
&M_{\mathrm{US}}^{(n)} \Bigl( \Bigl\{ \lambda\in \mathcal{SP}_n \,\Big|\, 
\sup_{x>0} \bigl| Y(\lambda)^{\sqrt{n}} (x)- \psi_{\mathrm{US}} (x)\bigr| 
> \epsilon \Bigr\} \Bigr) \ 
\xrightarrow[\,n\to\infty\,] \ 0, 
\label{eq:1-3} \\ 
&M_{\mathrm{US}}^{\delta^\prime (N)} \Bigl( \Bigl\{ \lambda\in \mathcal{SP} \,\Big|\, 
\sup_{x>0} \bigl| Y(\lambda)^{\sqrt{N}} (x)- \psi_{\mathrm{US}} (x)\bigr| 
> \epsilon \Bigr\} \Bigr) \ 
\xrightarrow[\,n\to\infty\,] \ 0
\label{eq:1-3.5} 
\end{align}
where 
\begin{equation}\label{eq:1-5}
\psi_{\mathrm{U}} (x) = - \frac{\sqrt{6}}{\pi} \log ( 1- e^{-\pi x/ \sqrt{6}}), 
\quad 
\psi_{\mathrm{US}} (x) = \frac{2\sqrt{3}}{\pi} \log ( 1+ e^{-\pi x/ 2\sqrt{3}}), 
\quad x>0. 
\end{equation}
The functions of \eqref{eq:1-5} are often called Vershik curves. 
In a similar way to shift from $Y(\lambda)$ to $S(\lambda)$, let us consider 
the transformation on the $xy$ plane 
\begin{equation}\label{eq:1-6}
(x, y) \ \longmapsto \ \bigl( \frac{x+y}{\sqrt{2}}, \frac{y}{\sqrt{2}}
\bigr). 
\end{equation}
Under \eqref{eq:1-6}, the graph $y =\psi_{\mathrm{US}}(x)$ of \eqref{eq:1-5} coincides 
with the right half (the part in $y\leqq x$) of $y=\psi_{\mathrm{U}}(x)$ of \eqref{eq:1-5} 
(see in \S\S\ref{subsec:A1}). 

Next we take the Plancherel measure $M_{\mathrm{Pl}}^{(n)}$ defined by 
\begin{equation}\label{eq:1-7}
M_{\mathrm{Pl}}^{(n)} \bigl( \{ \lambda\}\bigr) = \frac{(\dim\lambda)^2}{n!}, 
\qquad \lambda\in \mathcal{P}_n
\end{equation}
where $\dim\lambda$ is the dimension of an irreducible representation of 
$\mathfrak{S}_n$ corresponding to $\lambda\in \mathcal{P}_n$. 
Now we put $Y(\lambda)$ for $\lambda\in \mathcal{P}_n$ in the region of $y\geqq |x|$ of 
the $xy$ plane as in the right of Figure~\ref{fig:1-1}. 
Here the length of an edge of each box is $\sqrt{2}$ so that every vertex 
lies in $\mathbb{Z}^2$. 
The profile in this display style is again denoted by $Y(\lambda)$. 
The profile rescaled by $r >0$ as in \eqref{eq:1-1} is also denoted by $Y(\lambda)^r$. 
In this case, we have 
\begin{equation}\label{eq:1-7.5}
\int_{\mathbb{R}} \bigl( Y(\lambda)^{\sqrt{n}}(x) -|x|\bigr) \, dx = 2, \qquad 
\lambda\in\mathcal{P}_n.
\end{equation}
The law of large numbers due to Vershik--Kerov \cite{VeKe77} and 
Logan--Shepp \cite{LoSh77} tells us that: \ 
for $\forall\epsilon >0$
\begin{align}
&M_{\mathrm{Pl}}^{(n)} \Bigl( \Bigl\{ \lambda\in \mathcal{P}_n \,\Big|\, 
\sup_{x\in\mathbb{R}} \bigl| Y(\lambda)^{\sqrt{n}} (x)- 
\varOmega_{\mathrm{VKLS}} (x)\bigr| > \epsilon \Bigr\} \Bigr) \ 
\xrightarrow[\,n\to\infty\,] \ 0, 
\label{eq:1-9} \\ 
&\varOmega_{\mathrm{VKLS}} (x) = \begin{cases} \frac{2}{\pi} \bigl( 
x\arcsin\frac{x}{2} + \sqrt{4-x^2}\bigr), & |x|\leqq 2\\ 
\quad |x|, & |x|>2. \end{cases}
\label{eq:1-10}
\end{align}
The function of \eqref{eq:1-10} is often called the VKLS curve. 
A shifted version of \eqref{eq:1-9} was given by Ivanov (\cite{Iva04}, \cite{Iva06}). 
See also \cite{DeS17} and \cite{MaSn20}. 
This shifted version will be referred to later after touching upon spin representations. 
Biane extended limit shapes of Young diagrams like \eqref{eq:1-9} arising from 
group representations to various statistical ensembles of Young diagrams by 
formulating the notion of \textit{approximate factorizaton property} (\cite{Bia01}). 
This notion palys an important role in the present paper and are described later in detail 
(Definition~\ref{def:AFP}). 
Related works are in \cite{Sni06}, \cite{DoFe16}, \cite{MaSn20}, in which central limit theorems 
are also obtained. 
All the results mentioned up to here fall within static models. 

Dynamic models describe macroscopic time evolution of limit shapes and provide us 
with interesting mathematical problems. 
A typical scheme can be formulated as follows. 
Take a parameter $N>0$ indicating magnitude of the system considered, and run a continuous 
time Markov chain $(X_s^{(N)})$ which describes transitions between Young diagrams. 
Assume that the rescaled profile $Y(\lambda)^{\sqrt{N}}$ concentrates at a limit shape 
as $N\to\infty$ under the distribution of $X_0^{(N)}$. 
Then we ask whether, for a macroscopic time $t>0$, the concentration of 
$Y(\lambda)^{\sqrt{N}}$ as $N\to\infty$ is propagated under the distribution of 
$X_{tN}^{(N)}$. 
If it is the case, we intend to describe time evolution of the limit shape $\psi(t,x)$ 
along macroscopic time $t$ in a sufficiently concrete manner. 
When $(X_s^{(N)})$ keeps probability $M^{(N)}$ invariant and rescaled profiles 
converge to limit shape $\psi(x)$ as $N\to\infty$ under $M^{(N)}$, it is expected 
that $\psi(t,x)$ converges to $\psi(x)$ as $t\to\infty$. 

In \cite{FuSa10}, Funaki--Sasada studied hydrodynamic limit in such a diffusive 
time/space scaling and obtained nonlinear partial differential equations for the 
limit shape $\psi(t,x)$. 
Their continuous time Markov chains are constructed in grand canonical setting 
and have $M_{\mathrm{U}}^{\delta(N)}$ and $M_{\mathrm{US}}^{\delta^\prime (N)}$ 
as invariant probabilities. 
According as the partitions or the strict partitions, the obtained PDE's for 
$\psi = \psi (t, x)$ are 
\begin{align}
&\partial_t \psi = \partial_x \Bigl( \frac{\partial_x\psi}{1-\partial_x\psi}
\Bigr) + \frac{\pi}{\sqrt{6}} \frac{\partial_x\psi}{1-\partial_x\psi}, 
\label{eq:1-11} \\ 
&\partial_t\psi = \partial_{xx}\psi + \frac{\pi}{2\sqrt{3}} 
\partial_x\psi (1+\partial_x\psi)
\label{eq:1-12}
\end{align}
respectively. 
See also \cite{Fun16} for further models. 
To see the situation of $t$ tends to $\infty$, set $\partial_t\psi =0$ in 
\eqref{eq:1-11} and \eqref{eq:1-12}. 
Then they become ODE's which have the Vershik curves 
\eqref{eq:1-5} as solutions. 
In a similar way the two of \eqref{eq:1-5} are related, the transformation 
\eqref{eq:1-6} connects the PDE's \eqref{eq:1-11} and \eqref{eq:1-12} 
(see \S\S\ref{subsec:A1}). 

A Markov chain arising from restiction and induction of irreducible representations of 
symmetric groups in canonical setting keeps the Plancherel measure invariant 
(see \eqref{eq:1-41}). 
Fulman (\cite{Ful04}, \cite{Ful05}) introduced this Markov chain to obtain CLT for 
Plancherel and Jack measures via Stein's method. 
Borodin--Olshanski (\cite{BoOl09}) used this chain to construct a diffusion process 
on the Thoma simplex. 
In earlier works (\cite{Hor15}, \cite{Hor16}), we studied time evolution of the 
limit shapes obtained from this continuous time Markov chain under diffusive scaling limit. 
We gave an explicit form of time evolution for the Kerov transition measure and PDE 
for its Stieltjes transform $G = G(t, z)$ as 
\begin{equation}\label{eq:1-12.5}
\frac{\partial G}{\partial t} = -G\, \frac{\partial G}{\partial z} +G + 
\frac{1}{G}\, \frac{\partial G}{\partial z} 
\end{equation}
by using notions of free probability theory. 
Moreover, in \cite{Hor20}, we extended such a diffusive limit to non-Markovian models 
by admitting non-exponential pausing time at each state. 

As dynamical models related to limit shapes but different from ours, a wide variety 
of growth models have been studied. 
In addition to the famous Plancherel growth process due to Kerov (\cite{Ker99}, \cite{Ker03}), 
we refer to \cite{Str08} and \cite{BoBuOl15} also. 
The Plancherel growth process will be mentioned explicitly after Corollary~\ref{cor:PDE} in 
comparison to our result. 

The aim of this paper is to investigate such a dynamical model of diffusive limit for 
continuous time random walks as in \cite{Hor15}, \cite{Hor16} and \cite{Hor20} in 
a new framework based on \textit{spin irreducible representations} of symmetric groups. 
Another important ingredient is to reveal an interesting property of a spin version of 
the Jucys--Murphy elements for symmetric groups which forms an essential step in 
the proof of the main theorem. 

Let us briefly review some necessary notions in spin representations of a symmetric group. 
We mention here \cite{HoHu92}, \cite{Kle05} and \cite{Hir18} for references. 
A projective representation of group $G$ is by definition a homomorphism from 
$G$ to a projective linear group $PGL(V)$. 
Generally speaking, importance of projective representations is seen in that they 
appear quite naturally, for example, in: 

\noindent 
-- formulation of a state vector as an element of a Hilbert space \textit{modulo scalar 
multiple\,} in quantum mechanics

\noindent 
-- theory of Mackey's machine constructing an irreducible representation of a group 
by induction from an irreducible representaion of its subgroup. 

Fundamentals of projective representations for symmetric groups were established 
by Schur. 
As described below, we have two classes: ordinary representations and spin representations. 
The relation of the two is not direct. 
While they have some common features, their representation theories are different 
in many aspects. 
Especially, we focus on the branching rule and its asymptotic behavior 
(see Figure~\ref{fig:branch}). 

For a finite group $G$, there exists a finite group $\widetilde{G}$ called a 
representation group of $G$. 
It is a central extension of $G$: 
\[
\{e\} \longrightarrow Z \longrightarrow \widetilde{G} 
\xrightarrow{\ \Phi\ } \  G \longrightarrow \{e\}, 
\] 
where $Z$ is isomorphic to the Schur multiplier. 
Considerations of projective representations of $G$ result in those of linear representations 
of $\widetilde{G}$ in which the elements of $Z$ act as scalars. 
Let $G$ be the symmetric group $\mathfrak{S}_n$ of degree $n$. 
If $n\geqq 4$, its representation groups are constructed in two ways as double 
covers of $\mathfrak{S}_n$, and in fact the two are not isomorphic unless $n=6$. 
In this paper, $\widetilde{\mathfrak{S}}_n$ denotes the group generated by 
$z, r_1, r_2, \cdots, r_{n-1}$ under the fundamental relations 
\begin{equation}\label{eq:1-13}
z^2 =e, \quad r_iz = zr_i, \quad r_i^2 =e, \quad (r_i r_{i+1})^3 =e, 
\quad r_ir_j = zr_jr_i \ (|i-j| \geqq 2).
\end{equation}
Here $e$ is the identity element. 
The group $\widetilde{\mathfrak{S}}_n$ is a double cover of $\mathfrak{S}_n$, 
and is a representation group of $\mathfrak{S}_n$ if $n\geqq 4$: 
\[
\mathfrak{S}_n \cong \widetilde{\mathfrak{S}}_n /Z, \qquad 
Z = \langle z\rangle \cong \mathbb{Z}_2. 
\] 
See \eqref{eq:2-13} for the other representation group $\widetilde{S}_n$. 
Since there exists a correspondence between the actions of irreducible 
representation matrices as indicated in \eqref{eq:2-14}, it suffices to treat 
linear representations of one of the representation groups. 
The canonical homomorphism 
$\Phi : \widetilde{\mathfrak{S}}_n \longrightarrow \mathfrak{S}_n$ 
is given by 
\begin{equation}\label{eq:1-14}
\Phi(z) = e, \quad \Phi(r_i) = s_i 
\end{equation}
where $s_i = (i\ i\!+\!1)$ is a simple transposition in $\mathfrak{S}_n$. 
Define $\mathrm{sgn}\,x$ to be $\mathrm{sgn}\,\Phi(x)$ for 
$x\in \widetilde{\mathfrak{S}}_n$ also. 
If $(\tau, V)$ is a linear representation of $\widetilde{\mathfrak{S}}_n$ where 
$Z$ acts as scalars, $\tau(z)$ is equal to $id_V$ or $-id_V$. 
We call $\tau$ ordinary or spin according as $\tau(z)= id_V$ or 
$\tau(z)= -id_V$. 
An ordinary representation of $\widetilde{\mathfrak{S}}_n$ is identified with 
a linear representation of $\mathfrak{S}_n$. 
The set $(\widetilde{\mathfrak{S}}_n)^{\wedge}$ of linear equivalence classes 
of irreducible linear representations of $\widetilde{\mathfrak{S}}_n$ is 
divided into the classes of ordinary representations and the classes of 
spin representations: 
\begin{equation}\label{eq:1-15}
(\widetilde{\mathfrak{S}}_n)^{\wedge} = 
(\widetilde{\mathfrak{S}}_n)^{\wedge}_{\mathrm{ord}} \sqcup 
(\widetilde{\mathfrak{S}}_n)^{\wedge}_{\mathrm{spin}}.
\end{equation}
The set of partitions $\mathcal{P}_n$ parametrizes 
$(\widetilde{\mathfrak{S}}_n)^{\wedge}_{\mathrm{ord}}$. 
For $\lambda\in \mathcal{SP}_n$, we write as $\lambda\in \mathcal{SP}_n^+$ or 
$\lambda\in \mathcal{SP}_n^-$ 
according as $n-l(\lambda)$ is even or odd. 
In this paper, Nazarov's parameters (\cite{Naz90}) 
\begin{equation}\label{eq:1-16}
(\lambda, \gamma) \qquad \text{where} \quad \lambda\in \mathcal{SP}_n, 
\quad \gamma \in \{1,-1\}
\end{equation}
are often used to parametrize $(\widetilde{\mathfrak{S}}_n)^{\wedge}_{\mathrm{spin}}$, 
where $(\lambda, 1) = (\lambda, -1)$ for $\lambda\in \mathcal{SP}_n^+$ and 
$(\lambda, 1) \neq (\lambda, -1)$ for $\lambda\in \mathcal{SP}_n^-$. 
Let $\tau_{\lambda, \gamma}$ be a spin irreducible representation of 
$\widetilde{\mathfrak{S}}_n$ realizing parameter $(\lambda, \gamma)$. 
Then $\tau_{\lambda, \gamma}$ is self-associate (that is, equivalent to the 
sign multiple one) for $\lambda\in \mathcal{SP}_n^+$, while it is not self-associate and 
$\tau_{\lambda, -1} \cong \mathrm{sgn}\cdot \tau_{\lambda, 1}$ holds for 
$\lambda\in \mathcal{SP}_n^-$. 

For $\lambda\in \mathcal{SP}$, we display $S(\lambda)$ as the right of Figure~\ref{fig:1-2} 
and construct a shift-symmetric diagram $D(\lambda)$ by pasting its inversion 
to the left side. 
See the left of Figure~\ref{fig:1-3}. 
Note that we here adopt slightly different pasting from 
\cite[Definition, p.184]{HoHu92}, \cite{Mat18} and \cite{MaSn20}. 
We have $|D(\lambda)| =2n$ for $\lambda\in \mathcal{SP}_n$ and hence 
\[ 
\int_{\mathbb{R}} \bigl( D(\lambda)^{\sqrt{2n}}(x) -|x|\bigr) \, dx =2, 
\qquad \lambda\in \mathcal{SP}_n
\] 
as \eqref{eq:1-7.5}. 
When $S(\mu)$ is formed by adding one box to $S(\lambda)$ for 
$\lambda\in \mathcal{SP}_n$ and $\mu\in \mathcal{SP}_{n+1}$, we write as 
$\lambda\nearrow\mu$ or $\mu\searrow\lambda$. 
If the box $\mu/\lambda$ sits in the $i$th row and $j$th column, its content 
$j-i$ is denoted by $c(\mu/\lambda)$. 
For example, when $\mu =(3,2) \searrow \lambda =(3,1)$, $c(\mu/\lambda) =1$ 
with $i=2$, $j=3$ as in Figure~\ref{fig:2-1} (and the right of Figure~\ref{fig:1-3}). 
Note that $c(\mu /\lambda)$ and $-c(\mu /\lambda)-1$ are valleys of 
the profile of $D(\lambda)$. 

\begin{figure}[hbt]
\centering 
\include{Dlambda-fig1}
\vspace{-15mm}
\caption{$D(\lambda)$ for $\lambda = (3,1)\in \mathcal{SP}_4$, \ 
$\lambda\nearrow\mu= (3,2) \in \mathcal{SP}_5$}
\label{fig:1-3}
\end{figure}
\begin{figure}[hbt]
\centering 
\include{content-fig}
\vspace{-10mm}
\caption{contents for $\mu = (3,2)\in \mathcal{SP}_5$}
\label{fig:2-1}
\end{figure}

The Jucys--Murphy elements of $\mathfrak{S}_n$ are defined by 
\begin{equation}\label{eq:1-17.5}
J_k = (1\ k) +(2\ k)+ \cdots + (k-1\ k) \quad 
(2\leqq k \leqq n), \qquad J_1 =0, 
\end{equation}
which play important roles in representation theory of symmetric groups and have a 
wide variety of applications including topology and physics. 
In particular, Biane (\cite{Bia95}, \cite{Bia98}) used Jucys--Murphy elements to reveal 
beautiful connection between  permutations and free random variables. 
For example, the distribution of $J_n /\sqrt{n-1}$ with respect to state $\delta_e$ on 
$\mathbb{C}[\mathfrak{S}_n]$ converges to the standard semicircle distribution, 
which opens the door to free central limit theorem. 
The role of Jucys--Murphy elements in this paper lies in the stream due to Biane. 

Let us see how the spin version of $J_k$ is introduced. 
Using the generators of \eqref{eq:1-13}, define transposition $[i\ j]$ in 
$\widetilde{\mathfrak{S}}_n$ by 
\begin{equation}\label{eq:1-18}
[i\ j] = z^{j-i-1} r_{j-1}\cdots r_{i+1}r_ir_{i+1}\cdots r_{j-1}, \quad 
[j \ i] = z [i\ j],  \qquad 1\leqq i<j\leqq n.
\end{equation}
Picking up distinct $i_1, \cdots, i_r$ from $\{1,\cdots, n\}$, define an 
$r$-cycle by 
\begin{equation}\label{eq:1-19}
[i_1\ i_2\ \cdots \ i_r] = [i_{r-1}\ i_r] \cdots [i_2\ i_r] [i_1\ i_r].
\end{equation}
Both \eqref{eq:1-18} and \eqref{eq:1-19} are consistent if they are 
projected onto $\mathfrak{S}_n$ by $\Phi$ of \eqref{eq:1-14}: 
\[ 
\Phi( [i_1\ i_2\ \cdots \ i_r] ) = (i_1\ i_2\ \cdots \ i_r). 
\] 
Similarly to \eqref{eq:1-17.5}, the Jucys--Murphy elements of $\widetilde{\mathfrak{S}}_n$ 
are defined by 
\begin{equation}\label{eq:1-21}
\widetilde{J}_k = [1\ k] +[2\ k]+ \cdots + [k-1\ k] \quad 
(2\leqq k \leqq n), \qquad \widetilde{J}_1 =0. 
\end{equation}
See \cite{Kle05}. 
As seen from Lemma~\ref{lem:2-2}, $\widetilde{J}_k^{\;2}$ commutes with 
$r_i$ if $k\neq i, i+1$. 

Take the restriction homomorphism 
$\widetilde{E}_n : \mathbb{C}[\widetilde{\mathfrak{S}}_{n+1}] 
\longrightarrow \mathbb{C}[\widetilde{\mathfrak{S}}_n]$, which can be interpreted 
as a conditional expectation, 
\begin{equation}\label{eq:1-22} 
\widetilde{E}_n x = \begin{cases} x, & x\in \widetilde{\mathfrak{S}}_n \\ 
0, & x\in \widetilde{\mathfrak{S}}_{n+1} \setminus \widetilde{\mathfrak{S}}_n.
\end{cases}
\end{equation}
Here we have natural inclusion 
$\widetilde{\mathfrak{S}}_n \subset \widetilde{\mathfrak{S}}_{n+1}$ 
by identifying the generators $r_i, z$ of $\widetilde{\mathfrak{S}}_n$ 
with those of $\widetilde{\mathfrak{S}}_{n+1}$. 
The Kerov transition measure of $D(\lambda)$ is denoted by 
$\mathfrak{m}_{D(\lambda)}$. 
See \S\S\ref{subsec:A2} for the definition of Kerov transition measures. 
The irreducible character corresponding to 
$\xi\in (\widetilde{\mathfrak{S}}_n)^{\wedge}$ is denoted by $\chi^\xi$. 
In \cite{Bia03}, Biane gave a certain trace formula for Jucys--Murphy 
elements of a symmetric group, which reflects their spectral properties. 
The following is regarded as a spin version of Biane's formula. 

\begin{theorem}\label{th:JMtr}
The Jucys--Murphy element 
$\widetilde{J}_{n+1} = [1 \ n+\!1] + \cdots + [n \ n+\!1] \in 
\mathbb{C}[\widetilde{\mathfrak{S}}_{n+1}]$ 
satisfies 
\begin{equation}\label{eq:1-23}
\frac{\chi^{(\lambda, \gamma)}
\bigl( \widetilde{E}_n \widetilde{J}_{n+1}^{\ 2k}\bigr)}
{\dim(\lambda,\gamma)} = 
\sum_{\mu\in \mathcal{SP}_{n+1}: \mu\searrow\lambda} 
\Bigl( \frac{c(\mu /\lambda) \bigl( c(\mu /\lambda)+1\bigr)}{2}\Bigr)^k 
\mathfrak{m}_{D(\lambda)}\bigl( \{ c(\mu /\lambda), -c(\mu /\lambda)-1\}\bigr)
\end{equation}
for $(\lambda, \gamma) \in (\widetilde{\mathfrak{S}}_n)^\wedge_{\mathrm{spin}}$ 
and $k\in\mathbb{N}$. 
\hfill $\square$
\end{theorem}

As stated in Proposition~\ref{prop:2-1},  
\[
\frac{1}{2} c(\mu /\lambda) \bigl( c(\mu /\lambda)+1\bigr)
\] 
in \eqref{eq:1-23} 
is an eigenvalue of an action of $\widetilde{J}_{n+1}^{\ 2}$. 

There is a bijective correspondence $f \longleftrightarrow M$ between the set of 
normalized central positive-definite functions
\footnote{Equivalently, unital tracial states on group algebra 
$\mathbb{C}[\widetilde{\mathfrak{S}}_n]$ by linearity.}  
on $\widetilde{\mathfrak{S}}_n$ 
and the set of probabilities on $(\widetilde{\mathfrak{S}}_n)^{\wedge}$ 
through the Fourier expansion: 
\begin{equation}\label{eq:1-25}
f(x) = \sum_{\xi\in (\widetilde{\mathfrak{S}}_n)^{\wedge}} 
M(\{\xi\}) \frac{\chi^\xi(x)}{\dim\xi}, \qquad 
x\in \widetilde{\mathfrak{S}}_n.
\end{equation}
Note that 
\[
\frac{\chi^\xi (z)}{\dim\xi} = \begin{cases} 1, & 
\xi\in (\widetilde{\mathfrak{S}}_n)^{\wedge}_{\mathrm{ord}} \\ 
-1, & \xi\in (\widetilde{\mathfrak{S}}_n)^{\wedge}_{\mathrm{spin}}. 
\end{cases}
\] 
The function $f$ in \eqref{eq:1-25} is called ordinary or spin according as 
$f(zx) = f(x)$ or $-f(x)$. 
By the decomposition of \eqref{eq:1-15}, ordinary and spin $f$ 
correspond to probability $M$ on 
$(\widetilde{\mathfrak{S}}_n)^{\wedge}_{\mathrm{ord}}$ and 
$(\widetilde{\mathfrak{S}}_n)^{\wedge}_{\mathrm{spin}}$ respectively. 
This paper discusses a concentration phenomenon for a family of 
probabilities on $(\widetilde{\mathfrak{S}}_n)^{\wedge}_{\mathrm{spin}}$ 
as $n\to\infty$. 

For $x\in\widetilde{\mathfrak{S}}_n$, consider the cycle type 
$\rho\in \mathcal{P}_n$ of $\Phi(x)\in\mathfrak{S}_n$ through \eqref{eq:1-14}. 
Let $\mathcal{OP}_n$ denote the subset of $\mathcal{P}_n$ consisting of partitions whose 
parts have all odd length and set 
\[ 
\mathcal{OP} = \bigcup_{n=0}^\infty \mathcal{OP}_n, \qquad 
\mathcal{OP}_0 = \{\varnothing\}. 
\] 
It is known that, if $\rho$ is the cycle type of $\Phi(x) = \Phi(zx)$, 
\begin{equation}\label{eq:1-26}
x\text{ and }zx\text{ are not conjugate in } \widetilde{\mathfrak{S}}_n  
\iff \rho\in \mathcal{OP}_n \sqcup \mathcal{SP}_n^-. 
\end{equation}
If $f$ is a spin central function on $\widetilde{\mathfrak{S}}_n$ and 
the type of $x\in \widetilde{\mathfrak{S}}_n$ does not belong to 
$\mathcal{OP}_n\sqcup \mathcal{SP}_n^-$, then $f(x)=0$ holds since 
$f(x) = f(zx) =-f(x)$. 
For $x\in \widetilde{\mathfrak{S}}_n$, let $\mathrm{type}\,x$ denote 
the type of $\Phi(x)\in\mathfrak{S}_n$, and $\mathrm{supp}\, x$ 
denote the support of $\Phi(x)$ (that is, the letters not fixed by $\Phi(x)$). 
If $x, y\in \widetilde{\mathfrak{S}}_n$ satisfy 
$\mathrm{supp}\,x \cap \mathrm{supp}\, y =\varnothing$, then 
\begin{equation}\label{eq:1-27}
\mathrm{type}\,x = (\rho, 1^{n-|\rho|}), \quad 
\mathrm{type}\,y = (\sigma, 1^{n-|\sigma|})
\end{equation}
imply $\mathrm{type}\,xy = (\rho\sqcup\sigma, 1^{n-|\rho|-|\sigma|})$. 
Here $\rho$ of $(\rho, 1^{n-|\rho|})$ denotes the parts whose lengths are 
not $1$ (that is $m_1(\rho) =0$). 
Since in asymptotic representation theory one often takes the so-called dual 
approach in which parameters of irreducible representations play the role 
of variables while an element of a group or the type of a conjugacy class 
remains fixed, such a notation of a partition as $(\rho, 1^{n-|\rho|})$ is 
useful. 
Set $\widetilde{\mathfrak{S}}_\infty = 
\bigcup_{n=4}^\infty \widetilde{\mathfrak{S}}_n$. 
For $x\in \widetilde{\mathfrak{S}}_\infty$, since 
$x\in \widetilde{\mathfrak{S}}_n$ holds if $n\in\mathbb{N}$ is large 
enough, we have $\mathrm{type}\,x = (\rho, 1^{n-|\rho|})$ where 
$m_1(\rho) =0$. 

Take probability $M^{(n)}$ on $(\widetilde{\mathfrak{S}}_n)^\wedge_{\mathrm{spin}}$ 
and consider a concentration phenomenon for the family $\{ M^{(n)}\}$ as 
$n\to\infty$. 
We look at $\{ f^{(n)}\}$, the corresponding characteristic functions 
through \eqref{eq:1-25}. 
When $M^{(n)}$ and $f^{(n)}$ are connected by \eqref{eq:1-25}, 
$f^{(n)}$ is spin and central. 
As a type of elements of $\widetilde{\mathfrak{S}}_n$, take 
$(\rho, 1^{n-|\rho|})$ such that $\rho\in \mathcal{P}$ and $m_1(\rho)=0$. 
Since this type does not belong to $\mathcal{SP}$ for large $n$ enough, we have only to 
consider values of $f^{(n)}$ on the conjugacy classes of type 
\begin{equation}\label{eq:1-28}
(\rho, 1^{n-|\rho|}), \quad \rho\in \mathcal{OP}, \quad m_1(\rho)=0.
\end{equation}

\begin{definition}\label{def:AFP}
Let $f^{(n)}$ be a spin normalized central 
positive-definite function on $\widetilde{\mathfrak{S}}_n$ for 
$n\in\mathbb{N}$. 
The following condition for $\{f^{(n)}\}$ is called 
approximate factorization property after Biane \cite{Bia01}
\footnote{Letting $M^{(n)}$ be the probability on 
$(\widetilde{\mathfrak{S}}_n)^\wedge_{\mathrm{spin}}$ corresponding to $f^{(n)}$ 
by \eqref{eq:1-25}, we say also $\{M^{(n)}\}$ satisfies 
approximate factorization property.}: 
for $x, y\in \widetilde{\mathfrak{S}}_\infty\setminus \{e\}$ such that 
$\mathrm{supp}\,x\cap\mathrm{supp}\,y= \varnothing$ and 
\eqref{eq:1-27} with $\rho, \sigma$ satisfying \eqref{eq:1-28}, 
\begin{equation}\label{eq:1-29}
f^{(n)}(xy)- f^{(n)}(x) f^{(n)}(y) = o \Bigl( 
n^{-\frac{1}{2}(|\rho|-l(\rho)+|\sigma|-l(\sigma))} \Bigr) 
\qquad (n\to\infty).
\end{equation}
Note that the absolute value of the left hand side of \eqref{eq:1-29} does 
not change under $x \leftrightarrow zx$, $y \leftrightarrow zy$. 
\hfill $\square$
\end{definition}

The following theorem is a result for static models. 
In the statement we need free cumulant-moment formula, notion of continuous diagrams 
and Markov transforms, which are briefly explained in \S\S\ref{subsec:A2} 
and \S\S\ref{subsec:A3}. 
Theorem~\ref{th:static} is proved in \S\ref{sect:3} by way of consideration of 
the spin Jucys--Murphy operators. 

\begin{theorem}\label{th:static}
Let $f^{(n)}$ be a spin normalized central positive-definite function on 
$\widetilde{\mathfrak{S}}_n$ and assume the sequence $\{f^{(n)}\}$ 
(or $\{M^{(n)}\}$ connected by \eqref{eq:1-25}) satisfies 
approximate factorization property of Definition~\ref{def:AFP}. 
Assume also that there exists 
\begin{equation}\label{eq:1-30}
\lim_{n\to\infty} n^{\frac{k-1}{2}} f^{(n)} \bigl( [1\ 2\ \cdots \ k]\bigr) 
= r_{k+1}, \qquad k\in\mathbb{N}, \text{ odd, } \geqq 3.
\end{equation}
Setting 
\begin{equation}\label{eq:1-33}
r_1 =r_3 = r_5 =\cdots =0, \qquad r_2 =1,  
\end{equation}
we determine sequence $(m_j)$ from sequence $(r_j)$ by the free cumulant-moment 
formula \eqref{eq:a3-2}. 
Then, $(m_j)$ is a moment sequence and there exists probability $\mu$ on $\mathbb{R}$ 
such that 
\begin{equation}\label{eq:1-35}
m_j = \int_{\mathbb{R}} x^j \mu (dx)  \qquad (\forall j\in\mathbb{N}).
\end{equation}
(Hence $r_k$ is the $k$th free cumulant of $\mu$.) 
Furthermore, the following (1) and (2) hold. 

(1) If $r_{k+1}$ in \eqref{eq:1-30} satisfies 
\begin{equation}\label{eq:1-31}
\bigl\{|r_{k+1}|^{\frac{1}{k+1}}\bigr\} \text{ is bounded}, 
\end{equation}
$\mu$ is uniquely determined and has compact support. 
Take $\omega\in\mathbb{D}$ corresponding to $\mu$ via the 
Markov transform (see \eqref{eq:a2-3}), that is 
\begin{equation}\label{eq:1-310}
\mathfrak{m}_\omega = \mu. 
\end{equation}
Then we have concentration at the limit shape $\omega$ under $\{M^{(n)}\}$: \ 
for $\forall\epsilon >0$, 
\begin{equation}\label{eq:1-32}
\lim_{n\to\infty} M^{(n)} \Bigl( \Bigl\{ (\lambda, \gamma)\in 
(\widetilde{\mathfrak{S}}_n)^\wedge_{\mathrm{spin}} \,\Big|\, 
\sup_{x\in\mathbb{R}} \bigl| D(\lambda)^{\sqrt{2n}}(x) - \omega (x) \bigr| 
> \epsilon \Bigr\}\Bigr) =0.
\end{equation}
Here the right half ($x\geqq 0$) of the graph $y=\omega (x)$ becomes the limit shape 
of $S(\lambda)^{\sqrt{2n}}$.

(2) If \eqref{eq:1-31} is not assumed but the moment problem for \eqref{eq:1-35} 
is determinate, then the uniquely determined $\mu$ yields the limit shape $\omega$ 
by \eqref{eq:1-310}, and we have \eqref{eq:1-32}. 
Here $\omega$ is an element of $\mathcal{D}$ (see \eqref{eq:a2-2}) whose transition 
measure does not necessarily have compact support. 
\hfill $\square$ 
\end{theorem}

The law of large numbers for spin Plancherel measures was given by Ivanov 
(\cite{Iva04}, \cite{Iva06}). 
See also \cite{DeS17}. 
The spin Plancherel measure $M^{(n)}_{\mathrm{Pl,spin}}$ is described in 
\eqref{eq:1-44}, \eqref{eq:1-45} and \S\S\ref{subsec:5-2}. 
The corresponding positive-definite function is 
\[
\sum_{\xi\in (\widetilde{\mathfrak{S}}_n)^{\wedge}_{\mathrm{spin}}} 
M^{(n)}_{\mathrm{Pl,spin}}(\{\xi\}) \frac{\chi^\xi(g)}{\dim\xi} = 
\begin{cases} 1, & g=e\\ -1, & g=z\\ 0, & \text{otherwise}. \end{cases}
\]
The paper \cite{MaSn20} of  Matsumoto--\'Sniady is a comprehensive work in which 
they developed not only limit shapes but also fluctuations of shifted Young diagrams 
from the viewpoint of approximate factorization property. 
Their method is based on \'Sniady's theory for the case of linear representations of 
symmetric groups (see \cite{Sni06}). 
Since our aim is analysis of limit shapes, we have only to assume the weakest version 
\eqref{eq:1-29} for approximate factorization property. 

We are now in a position to proceed to our dynamic model. 
Given a finite group $G$ and its subgroup $H$, one naturally gets a Markov chain on $G$ 
reflecting the branching rule (\cite{Ful04}, \cite{Ful05}, \cite{BoOl09}). 
For $\xi\in \widehat{G}$ and $\eta\in \widehat{H}$, we can set 
\[ 
c_{\eta, \xi} = [ \mathrm{Res}^G_H \xi : \eta] = [\mathrm{Ind}_H^G \eta : \xi] 
\] 
by Frobenius' reciprocity to have branching rules 
\begin{equation}\label{eq:1-37} 
\mathrm{Res}^G_H \xi \cong \bigoplus_{\eta\in\widehat{H}} [c_{\eta, \xi}]\eta, 
\qquad 
\mathrm{Ind}^G_H \eta \cong \bigoplus_{\xi\in\widehat{G}} [c_{\eta, \xi}]\xi. 
\end{equation}
Form three stochastic matrices $P^\downarrow$, $P^\uparrow$ and $P$ 
\begin{align}
&P^\downarrow_{\xi, \eta} = \frac{c_{\eta, \xi} \dim\eta}{\dim\xi},
 \qquad 
P^\uparrow_{\eta, \xi} = \frac{c_{\eta, \xi} \dim\xi}{[G:H] \dim\eta}, 
\qquad \xi\in\widehat{G}, \ \eta\in\widehat{H}, 
\label{eq:1-38} \\ 
&P = P^\downarrow P^\uparrow \label{eq:1-39}, 
\end{align}
with the sizes of $|\widehat{G}|\times |\widehat{H}|$, 
$|\widehat{H}|\times |\widehat{G}|$ and 
$|\widehat{G}|\times |\widehat{G}|$ respectively. 
The Markov chain (Res-Ind chain) on $\widehat{G}$ with transition probability 
matrix $P$ is symmetric with respect to the Plancherel measure 
\begin{equation}\label{eq:1-40}
M_{\mathrm{Pl}} (\{\xi\}) = \frac{(\dim\xi)^2}{|G|} 
\end{equation}
on $\widehat{G}$ (see \eqref{eq:1-7} for the case of $G=\mathfrak{S}_n$). 
In other words, setting a diagonal matrix $D$ with diagonal entries $M_{\mathrm{Pl}}(\{\xi\})$, 
we have 
\begin{equation}\label{eq:1-41} 
^t(DP) = DP \qquad \text{where} \quad 
D = \mathrm{diag} (M_{\mathrm{Pl}}). \qquad 
\end{equation}

The branching rule for $G= \widetilde{\mathfrak{S}}_n$ and 
$H= \widetilde{\mathfrak{S}}_{n-1}$ is well known. 
Decomposed as in \eqref{eq:1-15}, $(\widetilde{\mathfrak{S}}_n)^\wedge_{\mathrm{ord}}$ 
and $\widehat{\mathfrak{S}_n}$ obey the same rule. 
For $(\widetilde{\mathfrak{S}}_n)^\wedge_{\mathrm{spin}}$, by using Nazarov's parameter 
\eqref{eq:1-16}, we have for $\gamma\in \{1, -1\}$ 
\begin{align}
\mathrm{Res}^{\widetilde{\mathfrak{S}}_n}_{\widetilde{\mathfrak{S}}_{n-1}}
&(\lambda, \gamma) \cong \notag \\ 
&\bigoplus_{\mu\in \mathcal{SP}_{n-1}^- : \lambda\searrow\mu} 
\bigl( (\mu, 1)\oplus (\mu, -1)\bigr) \ \oplus 
\delta_{1, \lambda_{l(\lambda)}} (\lambda^-, \gamma) & 
&\text{for} \quad \lambda\in \mathcal{SP}_n^+, 
\label{eq:1-42} \\ 
&\bigoplus_{\mu\in \mathcal{SP}_{n-1}^+ : \lambda\searrow\mu} 
(\mu, \gamma) \ \oplus 
\delta_{1, \lambda_{l(\lambda)}} (\lambda^-, \gamma) & 
&\text{for} \quad \lambda\in \mathcal{SP}_n^-. 
\label{eq:1-43}
\end{align}
Here $\lambda^-\in \mathcal{SP}_{n-1}$ denotes the partition formed by removing the last part 
of $\lambda\in \mathcal{SP}_n$ such that $\lambda_{l(\lambda)} =1$. 
We refer to \cite{Naz90} and \cite{HoHu92} for details on these branching rules. 
In \eqref{eq:1-43}, $\lambda^-\in \mathcal{SP}_{n-1}^-$ holds and $\gamma$ of 
$(\lambda^-,\gamma)$ agrees with $\gamma$ in the left hand side. 
The branching rule of \eqref{eq:1-42} and \eqref{eq:1-43} naturally produces 
a graph (like the Young graph) with the set of Nazarov's parameters as its 
vertex set. 
Two vertices $(\lambda, \gamma)$ and $(\mu, \delta)$ are joined by an edge 
if  $[(\lambda, \gamma) : (\mu, \delta)] >0$. 
Since the branching rule is multiplicity-free, the graph is simple as its beginning is 
presented in Figure~\ref{fig:branch}. 
If $\gamma =\pm 1$ are merged, one gets a graph with $\mathcal{SP}$ as its vertex set, 
which admits edges of multiplicity $1$ or $2$. 
This graph is called the Schur graph (see \cite{Bor99} and \cite{Pet10}). 

The Plancehrel measure \eqref{eq:1-40} on $(\widetilde{\mathfrak{S}}_n)^\wedge$ 
is decomposed according to \eqref{eq:1-15}: 
\begin{align}
&M^{(n)}_{\mathrm{Pl}} = \frac{1}{2} M^{(n)}_{\mathrm{Pl,ord}} + 
\frac{1}{2} M^{(n)}_{\mathrm{Pl,spin}}, \label{eq:1-44} \\ 
&M^{(n)}_{\mathrm{Pl,ord}} (\{\lambda\}) = \frac{(\dim\lambda)^2}{n!}, \qquad 
M^{(n)}_{\mathrm{Pl,spin}} \bigl( \{ (\lambda,\gamma)\}\bigr) = 
\frac{\dim(\lambda, \gamma)^2}{n!}. \label{eq:1-45} 
\end{align}
Probabilities $M^{(n)}_{\mathrm{Pl,ord}}$ and $M^{(n)}_{\mathrm{Pl,spin}}$ are on 
$(\widetilde{\mathfrak{S}}_n)^\wedge_{\mathrm{ord}}$ and 
$(\widetilde{\mathfrak{S}}_n)^\wedge_{\mathrm{spin}}$ respectively. 
The transition matrix of the Res-Ind chain on $(\widetilde{\mathfrak{S}}_n)^\wedge$ 
has a form of 
\begin{equation}\label{eq:1-46}
P = \begin{pmatrix} P_{\mathrm{ord}} & O \\ O & P_{\mathrm{spin}} \end{pmatrix}. 
\end{equation}
Let $P^{(n)}$ denote the part of $P_{\mathrm{spin}}$. 
Then $P^{(n)}$ is symmetric with respect to the Plancherel measure 
$M^{(n)}_{\mathrm{Pl,spin}}$ on $(\widetilde{\mathfrak{S}}_n)^\wedge_{\mathrm{spin}}$: 
\begin{equation}\label{eq:1-47}
M^{(n)}_{\mathrm{Pl,spin}}(\xi) P^{(n)}_{\xi\eta} = 
M^{(n)}_{\mathrm{Pl,spin}}(\eta) P^{(n)}_{\eta\xi}, \qquad 
\xi, \eta\in (\widetilde{\mathfrak{S}}_n)^\wedge_{\mathrm{spin}}. 
\end{equation}
The Markov chain on $(\widetilde{\mathfrak{S}}_n)^\wedge_{\mathrm{spin}}$ induced by 
$P^{(n)}$ is the Res-Ind chain we treat in the sequel. 
Let $(Z_k^{(n)})_{k\in\{0,1,2,\cdots\}}$ be the Markov chain on 
$(\widetilde{\mathfrak{S}}_n)^\wedge_{\mathrm{spin}}$ with transition matrix $P^{(n)}$ 
and initial distribution $M_0^{(n)}$. 
We introduce pausing time sequence $(\tau_j)_{j\in\mathbb{N}}$ which is IID and 
independent of $(Z_k^{(n)})$ also, each $\tau_j$ obeying distribution $\psi(dx)$ 
on $[0, \infty)$. 
This yields counting process $(N_s)_{s\geqq 0}$ as 
\[ 
N_s = \begin{cases} j, & 
\tau_1+\cdots +\tau_j \leqq s< \tau_1+\cdots +\tau_{j+1} \\ 
0, & s< \tau_1, \end{cases} \qquad N_0 =0 \ \text{a.s.} 
\] 
We assume $\psi\bigl( (0, \infty)\bigr) >0$ to get that $\tau_1+\cdots +\tau_j$ 
diverges to $\infty$ a.s. as $j\to\infty$. 
Consider a continuous time random walk on 
$(\widetilde{\mathfrak{S}}_n)^\wedge_{\mathrm{spin}}$ by setting 
\begin{equation}\label{eq:1-49}
X_s^{(n)} = Z_{N_s}^{(n)}, \qquad s\geqq 0.  
\end{equation}
We do not stick to Markovianity of $(X_s^{(n)})$ since the pausing time distribution 
need not be exponential. 
The transition probability of $(X_s^{(n)})$ satisfies the following: \ 
for $\xi, \eta\in (\widetilde{\mathfrak{S}}_n)^\wedge_{\mathrm{spin}}$, 
\begin{align}
&\mathrm{Prob}( X_s^{(n)} =\eta \,|\, X_0^{(n)} =\xi) \notag \\ 
&= \sum_{j=0}^\infty \mathrm{Prob}( Z_{N_s}^{(n)}=\eta \,|\, 
N_s=j, \ Z_0^{(n)}=\xi) \, \mathrm{Prob}( N_s=j\,|\, Z_0^{(n)}=\xi) 
\notag \\ 
&= \sum_{j=0}^\infty \mathrm{Prob}( Z_j^{(n)}=\eta \,|\, 
Z_0^{(n)}=\xi) \, \mathrm{Prob}( 
\tau_1+\cdots +\tau_j< s\leqq \tau_1+\cdots +\tau_{j+1}) 
\notag \\ 
&= \sum_{j=0}^\infty (P^{(n)\,j})_{\xi\eta} \int_{[0, s]} 
\psi\bigl( (s-u, \infty)\bigr) \psi^{\ast j} (du). \label{eq:1-50}
\end{align}
Here $\psi^{\ast j}$ denotes the $j$-fold convolution of $\psi$. 
By \eqref{eq:1-47},  $(X^{(n)}_s)_{s\geqq 0}$ also keeps the spin Plancherel 
measure $M^{(n)}_{\mathrm{Pl,spin}}$ invariant. 
Let $M^{(n)}_s$ be the distribution of $X^{(n)}_s$ and $f^{(n)}_s$ denote 
the spin normalized central positive-definite function on $\widetilde{\mathfrak{S}}_n$ 
determined by \eqref{eq:1-25}. 
We see from \eqref{eq:1-50} that, for 
$\eta\in (\widetilde{\mathfrak{S}}_n)^\wedge_{\mathrm{spin}}$, 
\begin{equation}\label{eq:1-51}
M_s^{(n)}(\{\eta\}) = \mathrm{Prob}(X_s^{(n)}=\eta) = 
\sum_{j=0}^\infty (M_0^{(n)} P^{(n)\,j})_\eta 
\int_{[0, s]} \psi\bigl( (s-u, \infty)\bigr) \psi^{\ast j} (du). 
\end{equation}

Now we state the main theorem of this paper. 
We need notions of free convolution $\mu\boxplus\nu$ and free compression 
$\mu_c$ for probabilities on $\mathbb{R}$, which are briefly explained in 
Appendix~\ref{subsec:A3}. 

\begin{theorem}\label{th:dynamic}
Assume $\{f^{(n)}_0\}_{n\in\mathbb{N}}$ corresponding to a sequence of 
initial distributions on $(\widetilde{\mathfrak{S}}_n)^\wedge_{\mathrm{spin}}$'s 
satisfies approximate factorization property of Definition~\ref{def:AFP}. 
Assume also that there exists 
\begin{equation}\label{eq:1-52}
\lim_{n\to\infty} n^{\frac{k-1}{2}} f^{(n)}_0 \bigl( [1\ 2\ \cdots\ k]\bigr) 
= r_{k+1}, \qquad 
\forall k\in\mathbb{N}, \text{ odd, } \geqq 3, 
\end{equation}
Making up sequence $(r_j)_{j=1}^\infty$ by adding \eqref{eq:1-33} and determining 
sequence $(m_j)_{j=0}^\infty$ by the 
free cumulant-moment formula \eqref{eq:a3-2}
\footnote{Theorem~\ref{th:static} assures $(m_j)$ is a moment sequence.}, 
let $\mu_0$ be a probability on $\mathbb{R}$ having $(m_j)$ as its moment 
sequence. 
For continuous time random walk $(X_s^{(n)})$ of \eqref{eq:1-49}, assume 
that pausing time distribution $\psi$ has mean $m$ and its characteristic 
function $\varphi$ satisfies an integrability condition
\footnote{$\varphi (u) = \int_{[0, \infty)} e^{iux}\psi(dx)$}  
\begin{equation}\label{eq:1-53}
\int_{\{|u|\geqq a\}} \Bigl| \frac{\varphi(u)}{u}\Bigr| du < \infty 
\qquad \text{for} \quad \exists a>0.  
\end{equation}
Then the following (1) and (2) hold. 

(1) If $r_{k+1}$ in \eqref{eq:1-52} satisfies 
\begin{equation}\label{eq:1-54}
\bigl\{|r_{k+1}|^{\frac{1}{k+1}}\bigr\} \text{ is bounded}, 
\end{equation}
$\mu_0$ is uniquely determined and has compact support. 
For any $t\geqq 0$, concentration occurs also under 
$\{M^{(n)}_{tn}\}_{n\in\mathbb{N}}$ and the limit shape $\omega_t\in\mathbb{D}$ 
is obtained: \ for $\forall \epsilon >0$, 
\begin{equation}\label{eq:1-55}
\lim_{n\to\infty} M^{(n)}_{tn} \Bigl( \Bigl\{ (\lambda,\gamma)\in 
(\widetilde{\mathfrak{S}}_n)^\wedge_{\mathrm{spin}} \,\Big|\, \sup_{x\in\mathbb{R}}
\bigl| D(\lambda)^{\sqrt{2n}}(x) - \omega_t(x)\bigr| >\epsilon \Bigr\}\Bigr) 
=0.
\end{equation}
We have $\mu_0 = \mathfrak{m}_{\omega_0}$. 
The limit shape $\omega_t$ is characterized in terms of its transition measure 
$\mathfrak{m}_{\omega_t}$ and the centered semicircle distribution $\gamma$ 
with variance $1$ by the formula. 
\begin{equation}\label{eq:1-56}
\mathfrak{m}_{\omega_t} = (\mathfrak{m}_{\omega_0})_{e^{-t/m}} 
\boxplus \gamma_{1-e^{-t/m}} \qquad (t>0), 
\end{equation}
or equivalently 
\begin{equation}\label{eq:1-57}
R_1(\mathfrak{m}_{\omega_t}) =0, \quad 
R_2(\mathfrak{m}_{\omega_t}) =1, 
\quad R_{k+1}(\mathfrak{m}_{\omega_t}) = 
R_{k+1}(\mathfrak{m}_{\omega_0}) e^{-kt/m} \ (k\geqq 2)
\end{equation}
by using free cumulants $R_k$'s
\footnote{We have $r_k = R_k(\mathfrak{m}_{\omega_0})$.}. 

(2) Assume $r_{k+1}$ in \eqref{eq:1-52} satisfies: \ $\exists C>0$ such that 
\begin{equation}\label{eq:1-571}
\sup_{k \text{ : odd, } \geqq 3} |r_{k+1}|^{\frac{1}{k+1}} \leqq C (k+1)
\end{equation}
instead of \eqref{eq:1-54}. 
Then, $\mu_0$ is uniquely determined and 
the concentration phenomenon is propagated at time $t>0$. 
Namely, we obtain the limit shape $\omega_t \in\mathcal{D}$ which satisfies 
\eqref{eq:1-55}. 
Moreover, \eqref{eq:1-56} and \eqref{eq:1-57} hold
\footnote{If $\mathrm{supp}\,\mathfrak{m}_{\omega_0}$ is not compact, 
the condition for free cumulants \eqref{eq:1-57} defines \eqref{eq:1-56}. 
The free cumulants are well-defined since the moment problem is determinate.}.
\hfill $\square$
\end{theorem}

Our assumption imposed on the pausing time excludes, for example, deterministic 
or heavy-tailed ones. 
In this paper, since mean $m$ is a completely fixed parameter, we might set $m=1$ 
without any loss of generality. 

The probability $\gamma$ has the form of 
\begin{equation}\label{eq:1-58}
\gamma (A) = \int_{A\cap [-2, \, 2]} \frac{1}{2\pi} 
\sqrt{4 -x^2}\, dx.
\end{equation}
After the Markov transform, $\gamma$ is the transition measure of 
$\varOmega_{\mathrm{VKLS}}$ in \eqref{eq:1-10}. 
The right half of $\varOmega_{\mathrm{VKLS}}$ agrees with the limit shape of 
$S(\lambda)^{\sqrt{2n}}$ under spin Plancherel 
measures $\{M^{(n)}_{\mathrm{Pl,spin}}\}$. 
Properties of probabilities convoluted with the semicircle distribution (like \eqref{eq:1-56}) 
were investigated in \cite{Bia97}. 

Let us consider the Stieltjes transform of $\mathfrak{m}_{\omega_t}$ of \eqref{eq:1-56}: 
\begin{equation}\label{eq:1-60}
G(t, z) = \int_{\mathbb{R}} \frac{1}{z-x}\, \mathfrak{m}_{\omega_t}(dx). 
\end{equation}
We restate time evolution of the limit shapes established in Theorem~\ref{th:dynamic} 
in terms of $G(t, z)$ as follows. 
The procedure of deriving partial differential equation \eqref{eq:1-61} from the 
expression of \eqref{eq:1-57} of free cumulants is the same as we did in \cite{Hor15} 
(and \cite{Hor16}) and hence is omitted. 

\begin{corollary}\label{cor:PDE}
In the framework of Theorem~\ref{th:dynamic}, PDE 
\begin{equation}\label{eq:1-61}
m \frac{\partial G}{\partial t} =  -G\,\frac{\partial G}{\partial z} + G + 
\frac{1}{G}\, \frac{\partial G}{\partial z}
\end{equation}
is satisfied by $G(t, z)$ of \eqref{eq:1-60}, the Stieltjes transform of \eqref{eq:1-56}. 
\hfill $\square$
\end{corollary}

As a long time situation, by setting $\partial G /\partial t =0$ in \eqref{eq:1-61}, 
we have ODE 
\begin{equation}\label{eq:1-62}
-G\,\frac{\partial G}{\partial z} + G + \frac{1}{G}\, \frac{\partial G}{\partial z} =0.  
\end{equation}
The Stieltjes transform of $\mathfrak{m}_{\varOmega_{\mathrm{VKLS}}} = \gamma$ 
of \eqref{eq:1-58} 
\[ 
\int_{\mathbb{R}} \frac{1}{z-x}\, \gamma (dx) = 
\frac{1}{2} \bigl( z-\sqrt{z^2-4}\bigr) 
\] 
is a solution of \eqref{eq:1-62}. 
As Corollary\ref{cor:PDE} shows, the Stieltjes transform of 
$\mathfrak{m}_{\omega_t}$ is governed by a nonlinear PDE. 
We do not have a functional (differential and/or integral) equation for $\omega_t$. 
On the other hand, the description of $\mathfrak{m}_{\omega_t}$ by using free 
convolution and free cumulants lets us see a piece of additive structure in some sense. 

Before closing Introduction, in order to compare \eqref{eq:1-12.5} (or \eqref{eq:1-61}) 
with the well-known corresponding equation 
for the Plancherel growth process (see \eqref{eq:1-65} below), 
let us recall the latter one in a similar context to ours. 
Consider the up matrix of \eqref{eq:1-38} for symmetric groups: 
\[
P^\uparrow _{\lambda\mu} = \frac{\dim\mu}{(|\lambda|+1) \dim\lambda} 
= \mathfrak{m}_{Y(\lambda)}\bigl( \{ c(\mu/\lambda)\}\bigr), \qquad 
\lambda\nearrow\mu, 
\] 
and Markov chain $(Z_k)_{k=0}^\infty$ on $\mathcal{P}$ with 
transition matrix $(P^\uparrow _{\lambda\mu})$ and initial distribution $\delta_\varnothing$. 
Taking Poisson process $(N_s)_{s\geqq 0}$ as a counting process of jumps 
(with exponentially distributed pausing time of mean $1$), we have 
continuous time Markov process $X_s = Z_{N_s}$, whose distribution at time $s$ 
is a Poissonized Plancherel measure on $\mathcal{P}$: 
\[
M_s = \mathrm{Prob}(X_s =~\cdot~) = 
\sum_{k=0}^\infty \frac{e^{-s} s^k}{k!} M_{\mathrm{Pl}}^{(k)}.
\] 
Let $N>0$ be a parameter of magnitude of the system, which means we consider 
rescaled diagram $Y(\lambda)^{\sqrt{N}}$, and set $s=tN$ for macroscopic time $t$ 
(diffusive scaling). 
The macroscopic time evolution of the Plancherel growth process due to Kerov 
(\cite{Ker99}, \cite{Ker03}) is described as follows: for $\forall\epsilon >0$, 
\[ 
M_{tN} \Bigl( \Bigl\{ \lambda\in \mathcal{P} \,\Big|\, 
\sup_{x\in\mathbb{R}} \bigl| Y(\lambda)^{\sqrt{N}} (x)- 
\sqrt{t}\varOmega_{\mathrm{VKLS}} \bigl( \frac{x}{\sqrt{t}}\bigr)\bigr| > 
\epsilon \Bigr\} \Bigr) \ \xrightarrow[\,N\to\infty\,] \ 0
\] 
where the transition measure of the limit shape is given by 
\begin{equation}\label{eq:1-64}
\frac{1}{2\pi t} \sqrt{4t-x^2}\, 1_{[-2\sqrt{t}, 2\sqrt{t}]}(x)\, dx.
\end{equation}
The Stieltjes transform $G(t, z)$ of \eqref{eq:1-64} satisfies complex Burgers equation 
\begin{equation}\label{eq:1-65}
\frac{\partial G}{\partial t} = - G\, \frac{\partial G}{\partial z}.
\end{equation}

The rest of this paper is organized as follows. 
In \S\ref{sect:2}, after reviewing some items on the spin Jucys--Murphy elements 
of symmetric groups, we show Theorem~\ref{th:JMtr}. 
Proof for the staic model, Theorem~\ref{th:static}, is given in \S\ref{sect:3}. 
The dynamic model, Theorem~\ref{th:dynamic}, is proved in \S\ref{sect:4}. 
In \S\ref{sect:5}, some concrete computations are mentioned for Vershik curves 
and ensembles arising from spin characters of $\mathfrak{S}_\infty$ 
as typical examples. 
Appendix includes some reviews on terminology of free probability and 
asymptotic representation theory.

\section{Spin Jucys-Murphy element for symmetric group}\label{sect:2}

Spin Jucys--Murphy elements are defined by \eqref{eq:1-21}. 
In particular, one of the Jucys--Murphy elements of 
$\widetilde{\mathfrak{S}}_{n+1}$ is 
\[ 
\widetilde{J}_{n+1} = [1\ n\!+\!1] + [2\ n\!+\!1] + \cdots + 
[n\ n\!+\!1]. 
\] 
The aim of this section is to show Theorem~\ref{th:JMtr} for the moments of 
$\widetilde{J}_{n+1}$. 
Take $(\mu, \delta)\in (\widetilde{\mathfrak{S}}_{n+1})^\wedge_{\mathrm{spin}}$ 
and a corresponding irreducible representation 
$(\tau_{\mu,\delta}, V_{\mu,\delta})$ of $\widetilde{\mathfrak{S}}_{n+1}$. 
In view of 
$\widetilde{\mathfrak{S}}_n\subset\widetilde{\mathfrak{S}}_{n+1}$, 
we have decomposition into irreducible $\widetilde{\mathfrak{S}}_n$-modules as 
\begin{equation}\label{eq:2-2}
V_{\mu,\delta} = \bigoplus_{(\lambda, \gamma) : 
[(\mu,\delta): (\lambda,\gamma)] >0} W_{\lambda,\gamma}, 
\qquad 
(\mu, \delta)\in (\widetilde{\mathfrak{S}}_{n+1})^\wedge_{\mathrm{spin}}
\end{equation}
by the branching rule \eqref{eq:1-42} and \eqref{eq:1-43}. 
Since $\widetilde{J}_{n+1}^{\ 2}$ commutes with $\widetilde{\mathfrak{S}}_n$ 
as seen from \eqref{eq:2-4}, Jucys--Murphy operator 
$\tau_{\mu,\delta}(\widetilde{J}_{n+1}^{\ 2})$ acts as a scalar on each 
summand of \eqref{eq:2-2}. 
We have $\mu\searrow\lambda$ if $[(\mu,\delta): (\lambda,\gamma)] >0$. 

\begin{proposition}\label{prop:2-1}
For $[(\mu,\delta): (\lambda,\gamma)] >0$ in \eqref{eq:2-2}, 
$\tau_{\mu,\delta}(\widetilde{J}_{n+1}^{\ 2})$ acts as a scalar 
\begin{equation}\label{eq:2-3}
\frac{1}{2} c(\mu/\lambda) \bigl( c(\mu/\lambda)+1\bigr) 
\end{equation}
on $W_{\lambda,\gamma}$. 
\hfill $\square$
\end{proposition}

Although the fact of Proposition~\ref{prop:2-1} is in \cite[22.3]{Kle05}, we here 
proceed in a framework of classical group theory. 

\begin{lemma}\label{lem:2-2}
As an element of $\mathbb{C}[\widetilde{\mathfrak{S}}_{n+1}]$ for $n\geqq 3$, 
we have 
\begin{equation}\label{eq:2-4}
\widetilde{J}_{n+1}^{\ 2} = \widetilde{A}_{(3,1^{n-2})} - 
\widetilde{A}_{(3,1^{n-3})} +n
\end{equation}
where $\widetilde{A}_{(3,1^{n-2})}$ {\rm [}resp. $\widetilde{A}_{(3,1^{n-3})}${\rm ]} 
is a conjugacy class indicator summed over the conjugacy class in 
$\widetilde{\mathfrak{S}}_{n+1}$ {\rm [}resp. $\widetilde{\mathfrak{S}}_n${\rm ]} 
containing $3$-cycle $[1\ 2\ 3]$. 
\hfill $\square$
\end{lemma}

\textit{Proof} \ 
Since $z^2 =1$, $[1\ 2\ 3]$ is invariant under a cyclic permutation of the letters. 
If $i, j, k\ (\leqq n)$ are distinct, $[1\ 2\ 3]$ and $[i\ j\ k]$ are conjugate in 
$\widetilde{\mathfrak{S}}_n$. 
On the other hand, $[1\ 2\ 3]$ and $z[1\ 2\ 3]$ are not conjugate in 
$\widetilde{\mathfrak{S}}_n$. 
Hence we have 
\begin{align}
&\widetilde{A}_{(3,1^{n-2})} = \frac{1}{3} 
\sum_{i,j,k\in \{1, \cdots, n+1\} : \text{distinct}} [i\ j\ k], \qquad 
\widetilde{A}_{(3,1^{n-3})} = \frac{1}{3} 
\sum_{i,j,k\in \{1, \cdots, n\} : \text{distinct}} [i\ j\ k], 
\notag \\ 
&\widetilde{A}_{(3,1^{n-2})} - \widetilde{A}_{(3,1^{n-3})} = 
\sum_{i,j\in \{1,\cdots, n\} : \text{distinct}} [i\ j\ n\!+\!1].
\label{eq:2-7}
\end{align} 
Comparing \eqref{eq:2-7} with 
\[ 
\widetilde{J}_{n+1}^{\ 2} = n+ 
\sum_{i,j\in \{1,\cdots, n\} : \text{distinct}} [j\ n\!+\!1] 
[i\ n\!+\!1], 
\] 
we get \eqref{eq:2-4}. 
\hfill $\blacksquare$

\begin{lemma}\label{lem:2-3}
In \eqref{eq:2-2}, the operators \ 
$\tau_{\mu,\delta}(\widetilde{A}_{(3,1^{n-2})})$, \  
$\tau_{\mu,\delta}(\widetilde{A}_{(3,1^{n-3})})$ \ 
act as scalars 
\begin{equation}\label{eq:2-10}
\frac{1}{3} (n+1)n (n-1) \frac{\chi^{(\mu,\delta)}([1\ 2\ 3])}
{\dim(\mu,\delta)}, \qquad 
\frac{1}{3} n (n-1)(n-2) \frac{\chi^{(\lambda,\gamma)}([1\ 2\ 3])}
{\dim(\lambda,\gamma)}
\end{equation}
respectively on $W_{\lambda,\gamma}$. 
\hfill $\square$
\end{lemma}

\textit{Proof} \ 
They act as scalars $\alpha, \beta$ respectively by Schur's lemma so that 
\begin{align*}
&\alpha \dim(\mu,\delta) = \mathrm{tr}\, \tau_{\mu,\delta} 
( \widetilde{A}_{(3,1^{n-2})}) = 
\chi^{(\mu,\delta)} ([1\ 2\ 3]) \frac{(n+1)n(n-1)}{3}, \\ 
&\beta \dim(\lambda,\gamma) = \mathrm{tr} \,\tau_{\mu,\delta} 
( \widetilde{A}_{(3,1^{n-3})})\big|_{W_{\lambda,\gamma}} = 
\mathrm{tr} \,\tau_{\lambda,\gamma} 
( \widetilde{A}_{(3,1^{n-3})}) = 
\chi^{(\lambda,\gamma)} ([1\ 2\ 3]) \frac{n(n-1)(n-2)}{3}. 
\end{align*}
These imply \eqref{eq:2-10}. 
\hfill $\blacksquare$

\medskip

\textit{Proof of Proposition~\ref{prop:2-1}}
Combining Lemma~\ref{lem:2-2} and Lemma~ \ref{lem:2-3}, we verify the value of 
\[ 
\frac{(n+1)n(n-1)}{3} \frac{\chi^{(\mu,\delta)}([1\ 2\ 3])}
{\dim\mu} - 
\frac{n(n-1)(n-2)}{3} \frac{\chi^{(\lambda,\gamma)}([1\ 2\ 3])}
{\dim\lambda} +n 
\] 
is equal to \eqref{eq:2-3}, 
where $\mu\in \mathcal{SP}_{n+1}$, $\lambda\in \mathcal{SP}_n$ with $\mu\searrow\lambda$. 
\hfill $\blacksquare$

\medskip

The above verification of the character values is performed e.g. by using Ivanov's formula 
in \cite{Iva06} or \cite[Theorem~6.4.11]{Hir18}. 
Ivanov's formula for spin irreducible characters of a symmetric group is given 
for the other representation group $\widetilde{S}_n$. 
Instead of \eqref{eq:1-13}, the fundamental relations for generators 
$z, r^\prime_1, \cdots, r^\prime_{n-1}$ of $\widetilde{S}_n$ are 
\begin{equation}\label{eq:2-13}
z^2 =e, \quad r^\prime_iz =zr^\prime_i, \quad {r^\prime_i}^2 =z, 
\quad (r^\prime_i r^\prime_{i+1})^3 = z, \quad 
r^\prime_i r^\prime_j = zr^\prime_j r^\prime_i \ 
(|i-j|\geqq 2).
\end{equation}
The correspondence between the spin irreducible representations are given by 
\begin{equation}\label{eq:2-14}
\text{representation matrix of } r^\prime_k\in \widetilde{S}_n 
\text{ is the $i$ multiple one of } r_k\in \widetilde{\mathfrak{S}}_n.
\end{equation}
See \cite[p.199]{HoHu92}. 
In \cite{Hir18}, $\widetilde{\mathfrak{S}}_n$ of \eqref{eq:1-13} is adopted as a 
representation group of $\mathfrak{S}_n$.

We get spectral decomposition of Jucys--Murphy operator 
$\tau_{\mu,\delta}(\widetilde{J}_{n+1}^{\ 2})$ by \eqref{eq:2-2} and 
Proposition~\ref{prop:2-1}, and hence will see the following. 
Let $g_\lambda$ denote the number of standard tableaux of shifted Young 
diagram $S(\lambda)$. 

\begin{proposition}\label{prop:2-5}
We have 
\begin{align}
\frac{\chi^{(\lambda,\gamma)}(\widetilde{E}_n \widetilde{J}_{n+1}^{\ 2k})}
{\dim(\lambda,\gamma)} &= 
\sum_{(\mu,\delta)\in (\widetilde{\mathfrak{S}}_{n+1})^\wedge_{\mathrm{spin}} : 
(\mu,\delta)\searrow (\lambda,\gamma)} 
\Bigl( \frac{c(\mu/\lambda) \bigl( c(\mu/\lambda)+1\bigr)}{2}\Bigr)^k 
\frac{\dim (\mu,\delta)}{(n+1)\dim (\lambda,\gamma)} 
\label{eq:2-23} \\ 
&= \sum_{\mu\in \mathcal{SP}_{n+1} : \mu\searrow\lambda} 
\Bigl( \frac{c(\mu/\lambda) \bigl( c(\mu/\lambda)+1\bigr)}{2}\Bigr)^k 
\frac{2 g_\mu}{(n+1) g_\lambda}. \label{eq:2-24}
\end{align}
for $(\lambda, \gamma)\in (\widetilde{\mathfrak{S}}_n)^\wedge_{\mathrm{spin}}$ 
and $k\in\mathbb{N}$. 
\hfill $\square$ 
\end{proposition}

Though \eqref{eq:2-23} is derived from a common argument (see \cite{Bia03}, \cite{Hor05}, 
\cite[Theorem~9.23]{HoOb07}), let us mention its proof for reader's convenience. 

\begin{lemma}\label{lem:2-6}
For finite group $G$ and its subgroup $H$, take restriction map 
$E : \mathbb{C}[G] \longrightarrow \mathbb{C}[H]$ as \eqref{eq:1-22}. 
Let $\xi\in\widehat{G}$, $\eta\in\widehat{H}$ and $b\in\mathbb{C}[G]$. 
Under the branching rule \eqref{eq:1-37}, if operator $\tau_\xi(b)$ acts as 
scalar $\beta_{\xi,\eta}$ on the $\eta$-isotypical component 
$(V_\eta)^{\oplus c_{\eta,\xi}}$, then 
\begin{equation}\label{eq:2-25}
\frac{\chi^\eta (Eb)}{\dim\eta} = \sum_{\xi\in\widehat{G}} 
\beta_{\xi,\eta} \, \frac{c_{\eta,\xi} \dim\xi}{[G:H] \dim\eta}.
\end{equation}
holds. 
\hfill $\square$
\end{lemma}

\textit{Proof} \ 
Let $E^\prime : \mathbb{C}[H] \longrightarrow \mathbb{C}[G]$ be the 
canonical embedding (extension). 
The minimal central projection of $\mathbb{C}[H]$ associated with $\eta\in\widehat{H}$ 
is given by 
\[ 
p_\eta (x) = \frac{\dim\eta}{|H|} \,\overline{\chi^\eta(x)}, \qquad x\in H. 
\] 
Then, we have 
\begin{align}
\chi^\eta (Eb) &= \sum_{x\in H} (Eb)(x) \chi^\eta (x) = 
\frac{|H|}{\dim\eta} \sum_{x\in H} (Eb)(x) \overline{p_\eta (x)} = 
\frac{|H|}{\dim\eta} \sum_{x\in G} b(x) \overline{(E^\prime p_\eta) (x)} 
\notag \\ 
&= \frac{|H| |G|}{\dim\eta} \sum_{\xi\in \widehat{G}} 
\frac{\dim\xi}{|G|^2} \mathrm{tr} \bigl( \tau_\xi(b) \mathrm{Res}^G_H 
\tau_\xi (p_\eta)\bigr). \label{eq:2-27}
\end{align}
Here the last equality comes from self-adjointness of $p_\eta$ and 
\[ 
\tau_\xi (E^\prime p_\eta) = \sum_{x\in G} (E^\prime p_\eta)(x) \tau_\xi(x) 
=\sum_{x\in H} p_\eta(x) (\mathrm{Res}^G_H \tau_\xi)(x) = 
(\mathrm{Res}^G_H \tau_\xi)(p_\eta) \] 
with the Plancherel theorem on $G$. 
Since the projection onto $\eta$-isotypical component in \eqref{eq:1-37} is 
given by $(\mathrm{Res}^G_H \tau_\xi)(p_\eta)$, \eqref{eq:2-27} continues as 
\[ 
= \frac{1}{[G:H] \dim\eta} \sum_{\xi\in\widehat{G}} \dim\xi\, 
\beta_{\xi,\eta} c_{\eta,\xi} \dim\eta.
\] 
Dividing by $\dim\eta$, we get \eqref{eq:2-25}. 
\hfill $\blacksquare$

\medskip

\textit{Proof of Proposition~\ref{prop:2-5}} \ 
Applying Proposition~\ref{prop:2-1} and Lemma~\ref{lem:2-6} to 
$G= \widetilde{\mathfrak{S}}_{n+1}$, $H= \widetilde{\mathfrak{S}}_n$ and 
$b = \widetilde{J}_{n+1}^{\ 2k}$, we readily get the first equality \eqref{eq:2-23}. 

To transform it into the second expression \eqref{eq:2-24}, recall a well-known formula 
giving the dimension of an irreducible representation: \ for 
$\lambda\in \mathcal{SP}_n$ and $\gamma\in \{ 1, -1\}$, 
\[ 
\dim (\lambda, \gamma) = 2^{\frac{n-l(\lambda)-\epsilon(\lambda)}{2}} g_\lambda, 
\qquad 
\epsilon(\lambda) = \begin{cases} 1, & n-l(\lambda) \text{ is odd} \\ 
0, & n-l(\lambda) \text{ is even}. \end{cases}
\] 
Let us divide the cases according as values of $\epsilon(\lambda)$ and 
$\lambda_{l(\lambda)}$. 

(i) For $\epsilon(\lambda)=1$ and $\lambda_{l(\lambda)} =1$, since 
$(\mu,\delta)\searrow (\lambda,\gamma)$ implies the unique $(\mu,\delta)$ 
such that $\epsilon(\mu)=0$ and $l(\mu)=l(\lambda)$, we have 
\begin{align*}
&\dim (\lambda,\gamma) = 2^{\frac{n-l(\lambda)-1}{2}} g_\lambda, 
\\ 
&\dim(\mu,\delta) = 2^{\frac{n+1-l(\mu)}{2}} g_\mu = 
2^{\frac{n+1-l(\lambda)}{2}} g_\mu 
\end{align*}
and hence 
\begin{equation}\label{eq:2-33}
\frac{\dim (\mu,\delta)}{\dim (\lambda,\gamma)} = 
\frac{2 g_\mu}{g_\lambda}. 
\end{equation}

(ii) For $\epsilon(\lambda) =1$ and  $\lambda_{l(\lambda)} \geqq 2$, 
$(\mu,\delta)\searrow (\lambda,\gamma)$ implies two cases: 
 
 \noindent 
(a) the unique $(\mu,\delta)$ such that 
$\epsilon(\mu)=1$, $\delta =\gamma$, $l(\mu)=l(\lambda)+1$, 
 
 \noindent 
(b) the unique $(\mu,\delta)$ such that $\epsilon(\mu)=0$, $l(\mu)=l(\lambda)$.
 
\noindent Case (a) Since 
\[ 
\dim(\mu,\delta) = 2^{\frac{n+1-l(\mu)-1}{2}} g_\mu = 
2^{\frac{n-l(\lambda)-1}{2}} g_\mu 
\] 
holds, we have $\dim(\mu,\delta) / \dim (\lambda,\gamma) = g_\mu / g_\lambda$. 
However, since $c(\mu/\lambda) =0$ holds, this case does not contribute to the 
sum of \eqref{eq:2-24}. 

\noindent Case (b) Since 
\[ 
\dim(\mu,\delta) = 2^{\frac{n+1-l(\mu)}{2}} g_\mu = 
2^{\frac{n+1-l(\lambda)}{2}} g_\mu 
\] 
holds, we have $\dim(\mu,\delta) / \dim (\lambda,\gamma) = 2 g_\mu / g_\lambda$, 
that is, \eqref{eq:2-33}. 

(iii) For $\epsilon(\lambda) =0$ and $\lambda_{l(\lambda)} =1$, we have 
\begin{equation}\label{eq:2-36}
\dim (\lambda, \gamma) = 2^{\frac{n-l(\lambda)}{2}} g_\lambda, 
\end{equation}
and $(\mu,\delta)\searrow (\lambda,\gamma)$ implies two 
$(\mu, \delta)$ ($\delta = 1, -1$) such that $\epsilon(\mu) =1$ and $l(\mu)=l(\lambda)$. 
Since 
\begin{equation}\label{eq:2-37}
\dim (\mu,\delta) = 2^{\frac{n+1-l(\mu)-1}{2}} g_\mu = 
2^{\frac{n-l(\lambda)}{2}} g_\mu, 
\end{equation}
we have $\dim(\mu,\delta) / \dim (\lambda,\gamma) = g_\mu / g_\lambda$. 
Combine two cases of $\delta =1$ and $ -1$ to have a term in \eqref{eq:2-24}. 

(iv) For $\epsilon(\lambda) =0$ and $\lambda_{l(\lambda)} \geqq 2$, 
$(\mu,\delta)\searrow (\lambda,\gamma)$ implies two cases: 

\noindent 
(a) two $(\mu,\delta)$'s such that $\epsilon(\mu)=1$ and $l(\mu)=l(\lambda)$, 

\noindent 
(b) the unique $(\mu,\delta)$ such that $\epsilon(\mu)=0$ and $l(\mu)=l(\lambda)+1$. 

\noindent Case (a) is similar to (iii) yielding \eqref{eq:2-37} and \eqref{eq:2-36}. 
In Case (b), we have 
\[ 
\dim (\mu,\delta) = 2^{\frac{n+1-l(\mu)}{2}} g_\mu = 
2^{\frac{n-l(\lambda)}{2}} g_\mu
\] 
and hence 
$\dim(\mu,\delta) / \dim (\lambda,\gamma) = g_\mu / g_\lambda$. 
However, since $c(\mu/\lambda) =0$ holds, this case does not contribute to the 
sum of \eqref{eq:2-24}. 

We have thus verified \eqref{eq:2-24}. 
\hfill $\blacksquare$

\medskip

To rewrite the ratio $g_\mu /g_\lambda$ in \eqref{eq:2-24} of Proposition~\ref{prop:2-5}, 
we use a hook formula expressing $g_\lambda$ for $\lambda\in \mathcal{SP}$. 
The hook length ($=$ arm length$+$leg length$+1$) at box $b$ of Young diagram 
$Y$ is denoted by $h_Y(b)$ (see Figure~\ref{fig:2-2}). 
\begin{figure}[hbt]
\centering
\include{hook-length-fig}
\vspace{-3mm}
\caption{hook length of box $\ast$}
\label{fig:2-2}
\end{figure}

\textbf{Hook length formula} \ (see \cite[Proposition10.6]{HoHu92}) 
\footnote{Since our $D(\lambda)$ is formed by slightly different pasting from a shift 
symmetric diagram in \cite{HoHu92}, the product in \eqref{eq:2-39} 
is taken over $D(\lambda)/S(\lambda)$.}: \ 
\begin{equation}\label{eq:2-39}
g_\lambda = \frac{n!}{\prod_{b\in D(\lambda)/S(\lambda)} h_{D(\lambda)}(b)}, 
\qquad \lambda\in \mathcal{SP}_n.
\end{equation}

\begin{lemma}\label{lem:2-7} 
We have 
\begin{equation}\label{eq:2-40}
\frac{g_\mu}{(n+1) g_\lambda} = \begin{cases} 
\frac{1}{2} \mathfrak{m}_{D(\lambda)} (\{ c(\mu/\lambda), -c(\mu/\lambda)-1\}), 
& c(\mu/\lambda) >0 \\ 
\mathfrak{m}_{D(\lambda)} (\{ c(\mu/\lambda)\}), & c(\mu/\lambda) =0. 
\end{cases}
\end{equation}
for $\lambda\in \mathcal{SP}_n$, $\mu\in \mathcal{SP}_{n+1}$ such that $\lambda\nearrow\mu$. 
\hfill $\square$
\end{lemma}

\textit{Proof} \ 
Since the hook length at a corner box is $1$, \eqref{eq:2-39} yields 
\begin{equation}\label{eq:2-41}
\frac{g_\mu}{(n+1) g_\lambda} = 
\frac{\prod_{b\in D(\lambda)/ S(\lambda)} h_{D(\lambda)}(b)}
{\prod_{b\in D(\mu)/ S(\mu)} h_{D(\mu)}(b)} = 
\prod_{b\in D(\lambda)/ S(\lambda)} \frac{h_{D(\lambda)}(b)}
{h_{D(\mu)}(b)}.
\end{equation}
We consider three cases according as how box $\mu/\lambda$ is added. 
One of the valleys of the profile of $D(\lambda)$ coincides with $c(\mu/\lambda)$. 

\noindent 
(i) Case of $\lambda_{l(\lambda)} =1$. 
In $D(\lambda)$, $0$ is not a valley, $-1$ is a peak, the number of valleys is even, 
and $c(\mu/\lambda) >0$. 

\noindent 
(ii) Case of $\lambda_{l(\lambda)} \geqq 2$. 
In $D(\lambda)$, $0$ is a valley, the number of valleys is odd, and there are 
(iia) case of $c(\mu/\lambda) =0$ \quad (iib) case of $c(\mu/\lambda) >0$. 

\begin{figure}[hbt]
\centering 
\include{zoneI-III-fig}
\caption{Case (i) for $\lambda = (7,6,4,3,1)\nearrow\mu = (7,6,5,3,1)$, 
$c(\mu/\lambda) = x_2 =4$}
\label{fig:2-3}
\end{figure}

Case (i) Let $D(\lambda)$ have interlacing coordinates 
\[ 
-x_r-1 < -y_{r-1}-1 < \cdots < -y_1-1 < -x_1-1 <-1 < x_1< y_1< 
\cdots < y_{r-1} < x_r 
\] 
%
in which $c(\mu/\lambda) = x_k$. 
When box $b\in D(\lambda)/S(\lambda)$ runs outside Zones I, II, III in 
Figure~\ref{fig:2-3}, we have $h_{D(\lambda)}(b) = h_{D(\mu)}(b)$. 
Letting box $b$ run along Zone I from the bottom, we notice the moments 
at which hook length decreases irregularly, not by $1$. 
Then we have 
\begin{equation}\label{eq:2-43}
\prod_{b\in \text{Zone I}} \frac{h_{D(\lambda)}(b)}{h_{D(\mu)}(b)} = 
\frac{(y_k-x_k)\cdots (y_{r-1}-x_k)}{(x_{k+1}-x_k)\cdots (x_r-x_k)}.
\end{equation}
Similarly when $b$ runs along Zone III, we have 
\begin{equation}\label{eq:2-44}
\prod_{b\in \text{Zone III}} \frac{h_{D(\lambda)}(b)}{h_{D(\mu)}(b)} = 
\frac{(x_k+y_{r-1}+1)\cdots (x_k+y_k+1)}{(x_k+x_r+1)\cdots (x_k+x_{k+1}+1)}.
\end{equation}
When $b$ runs along Zone II from the left, noticing the moments at which 
hook length varies irregularly, we see 
\begin{equation}\label{eq:2-45}
\prod_{b\in \text{Zone II}} \frac{h_{D(\lambda)}(b)}{h_{D(\mu)}(b)} = 
\frac{(x_k-y_{k-1})\cdots (x_k-y_1)x_k (x_k+y_1+1)\cdots (x_k+y_{k-1}+1)}
{(x_k-x_{k-1})\cdots (x_k-x_1)(x_k+x_1+1)\cdots (x_k+x_{k-1}+1) 2x_k}.
\end{equation}
Taking product of \eqref{eq:2-43}, \eqref{eq:2-44} and \eqref{eq:2-45}, 
we rewrite \eqref{eq:2-41} as 
\begin{equation}\label{eq:2-46}
\frac{(x_k-y_1) \cdots (x_k-y_{k-1})(x_k+y_1+1) \cdots (x_k+y_{r-1}+1)}
{2 (x_k-x_1)\cdots \widehat{(x_k-x_k)}\cdots (x_k-x_r) (x_k+x_1+1)\cdots 
\widehat{(x_k+x_k+1)}\cdots (x_k+x_r+1)}.
\end{equation}
Here $\widehat{(\ast\ast\ast)}$ implies to remove the indicated factor. 

On the other hand, partial fraction expansion of 
\[ 
G_{D(\lambda)}(z) = \frac{(z+y_{r-1}+1)\cdots (z+y_1+1) (z+1) 
(z-y_1)\cdots (z-y_{r-1})}
{(z+x_r+1)\cdots (z+x_1+1) (z-x_1)\cdots (z-x_r)} 
\] 
yields, through \eqref{eq:a2-1}, 
\begin{align}
\mathfrak{m}_{D(\lambda)}(\{x_k\}) &= 
\frac{(x_k+y_{r-1}+1)\cdots (x_k+y_1+1) (x_k+1) (x_k-y_1) \cdots (x_k-y_{r-1})}
{(x_k+x_r+1)\cdots (x_k+x_1+1) (x_k-x_1) \cdots \widehat{(x_k-x_k)}\cdots 
(x_k-x_r)}, \label{eq:2-48} \\ 
\mathfrak{m}_{D(\lambda)}(\{-x_k-1\}) &= 
\frac{(y_{r-1}-x_k)\cdots (y_1-x_k) x_k (y_1+x_k+1)\cdots (y_{r-1}+x_k+1)}
{(x_y-x_k)\cdots \widehat{(x_k-x_k)}\cdots (x_1-x_k) (x_1+x_k+1)\cdots 
(x_r+x_k+1)}. \label{eq:2-49}
\end{align}
Comparing \eqref{eq:2-48}+\eqref{eq:2-49} with \eqref{eq:2-46}, we get the 
first of the right hand side of \eqref{eq:2-40}. 

\begin{figure}[hbt]
\include{zoneVI-VIII-fig}
\caption{Case (ii) for $\lambda = (7,6,4,3,2)\nearrow\mu = (7,6,5,3,2)$; 
$c(\mu/\lambda)=4$ (left), $=0$ (right)}
\label{fig:2-4}
\end{figure}

Case (ii) Let $D(\lambda)$ have interlacing coordinates 
\[ 
-x_r-1 < -y_r-1 < \cdots < -x_1-1 < -y_1-1 < 0< y_1< x_1< \cdots y_r< x_r. 
\] 
(iia) When box $b$ runs outside Zones IV and V in Figure~\ref{fig:2-4} (right), 
we have $h_{D(\lambda)}(b) = h_{D(\mu)}(b)$. 
Letting $b$ run along Zone IV, we have 
\[ 
\prod_{b\in \text{Zone IV}} \frac{h_{D(\lambda)}(b)}{h_{D(\mu)}(b)} = 
\frac{y_r\cdots y_1}{x_r\cdots x_1}.
\] 
Also letting $b$ run along Zone V, we have 
\[ 
\prod_{b\in \text{Zone V}} \frac{h_{D(\lambda)}(b)}{h_{D(\mu)}(b)} = 
\frac{(y_r+1)\cdots (y_1+1)}{(x_r+1)\cdots (x_1+1)}.
\] 
Hence \eqref{eq:2-41} is equal to 
\[ 
\frac{y_1\cdots y_r (y_1+1)\cdots (y_r+1)}
{x_1\cdots x_r (x_1+1)\cdots (x_r+1)}.
\] 
On the other hand, 
\begin{equation}\label{eq:2-54}
G_{D(\lambda)}(z) = 
\frac{(z+y_r+1)\cdots (z+y_1+1)(z-y_1)\cdots (z-y_r)}
{(z+x_r+1)\cdots (z+x_1+1)z (z-x_1)\cdots (z-x_r)}
\end{equation}
yields 
\[ 
\mathfrak{m}_{D(\lambda)}(\{0\}) = 
\frac{(y_r+1)\cdots (y_1+1)(-y_1)\cdots (-y_r)}
{(x_r+1)\cdots (x_1+1)(-x_1)\cdots (-x_r)}.
\] 
Hence we get the second of the right hand side of \eqref{eq:2-40}. 

(iib) Set $c(\mu/\lambda) =x_k$. 
When box $b$ runs outside Zones VI, VII, VIII in Figure~\ref{fig:2-4} (left), 
we have $h_{D(\lambda)}(b) = h_{D(\mu)}(b)$. 
Similarly to Case (i), we have 
\begin{align*}
\prod_{b\in \text{Zone VI}} \frac{h_{D(\lambda)}(b)}{h_{D(\mu)}(b)} 
&= 
\frac{(y_r-x_k)\cdots (y_{k+1}-x_k)}{(x_r-x_k)\cdots (x_{k+1}-x_k)}, \\ 
\prod_{b\in \text{Zone VIII}} \frac{h_{D(\lambda)}(b)}{h_{D(\mu)}(b)} 
&= 
\frac{(x_k+y_r+1)\cdots (x_k+y_{k+1}+1)}
{(x_k+x_r+1)\cdots (x_k+x_{k+1}+1)}, \\ 
\prod_{b\in \text{Zone VII}} \frac{h_{D(\lambda)}(b)}{h_{D(\mu)}(b)} 
&= 
\frac{(x_k-y_k)\cdots (x_k-y_1)(x_k+y_1+1)\cdots (x_k+y_k+1)}
{(x_k-x_{k-1})\cdots (x_k-x_1)(x_k+1) (x_k+x_1+1)\cdots 
(x_k+x_{k-1}+1) 2x_k}. 
\end{align*}
Hence \eqref{eq:2-41} is equal to 
\begin{equation}\label{eq:2-59} 
\frac{(x_k-y_1)\cdots (x_k-y_r)(x_k+y_1+1)\cdots (x_k+y_r+1)}
{2 (x_k-x_1)\cdots \widehat{(x_k-x_k)}\cdots (x_k-x_r) 
(x_k+x_1+1)\cdots \widehat{(x_k+x_k+1))}\cdots (x_k+x_r+1)
x_k (x_k+1)}.
\end{equation}
On the other hand, since $G_{D(\lambda)}(z)$ is the same with \eqref{eq:2-54}, 
we have 
\begin{align}
\mathfrak{m}_{D(\lambda)}(\{x_k\}) &= 
\frac{(x_k+y_r+1)\cdots (x_k+y_1+1)(x_k-y_1)\cdots (x_k-y_r)}
{(x_k+x_r+1)\cdots (x_k+x_1+1)x_k (x_k-x_1) \cdots 
\widehat{(x_k-x_k)}\cdots (x_k-x_r)}, \label{eq:2-60} \\ 
\mathfrak{m}_{D(\lambda)}(\{-x_k-1\}) &= 
\frac{(x_k-y_r)\cdots (x_k-y_1)(x_k+y_1+1)\cdots (x_k+y_r+1)}
{(x_k-x_r) \cdots \widehat{(x_k-x_k)}\cdots (x_k-x_1)
(x_k+1)(x_k+x_1+1)\cdots (x_k+x_r+1)}. \label{eq:2-61} 
\end{align}
Reduce $2x_k+1$ in \eqref{eq:2-60}+\eqref{eq:2-61}. 
Comparing this with \eqref{eq:2-59} yields the first of the right hand side 
of \eqref{eq:2-40}. 
\hfill $\blacksquare$

\medskip

\textit{Proof of Theorem~\ref{th:JMtr}} \ 
Combine Proposition~\ref{prop:2-5} and Lemma~\ref{lem:2-7} to complete the 
proof of the theorem. 
\hfill $\blacksquare$

\medskip

Since the right hand side of \eqref{eq:1-23} is expressed by the transition measure 
of $D(\lambda)$, it contains visual information on the interlacing coordinates of 
$D(\lambda)$, while the left hand side of \eqref{eq:1-23} helps us carry out 
asymptotic arguments in combinatorics like certain random walk counting 
as is developed in the sequel.

\section{Proof for static model}\label{sect:3}

Under the assumption of Theorem~\ref{th:static}, we show some asymptotic 
behavior of $f^{(n)}$ on the center of $\mathbb{C}[\widetilde{\mathfrak{S}}_n]$
\footnote{We mean asymptotic behavior in the sense of dual approach. 
Fix type $(\rho, 1^{n-|\rho|})$ of a conjugacy class or let it vary in a bounded range, 
and let $n$ tend to $\infty$.}. 
Let us fix some notations. 

Set 
\begin{equation}\label{eq:3-1}
\overline{\mathcal{P}} = \{\rho\in \mathcal{P} \,|\, m_1(\rho) =0\}, \qquad 
\overline{\mathcal{OP}} = \{\rho\in \mathcal{OP} \,|\, m_1(\rho) =0\}.
\end{equation}
For $\rho\in \overline{\mathcal{P}}$ such that $n\geqq |\rho|$, set 
\begin{equation}\label{eq:3-2}
[ \rho, 1^{n-|\rho|}] = [ 1 \cdots \rho_1] [\rho_1\!+\! 1\cdots \rho_1\!+\!\rho_2] 
\cdots [\rho_1\!+\!\cdots\!+\!\rho_{l(\rho)-1}\!+\! 1 \ \cdots \ 
\rho_1\!+\!\cdots\! +\! \rho_{l(\rho)}] \in \widetilde{\mathfrak{S}}_n, 
\end{equation}
where the cycle notation in $\widetilde{\mathfrak{S}}_n$ was defined in \eqref{eq:1-19}. 
If $|\rho|=n$, \eqref{eq:3-2} is simply written as $[\rho]$. 
Further set 
\begin{align}
&\widetilde{C}_{(\rho, 1^{n-|\rho|})} = \text{conjugacy class of } 
\widetilde{\mathfrak{S}}_n \text{ to which } [\rho, 1^{n-|\rho|}] \text{ belongs}, 
\label{eq:3-3} \\ 
&\widetilde{A}_{(\rho, 1^{n-|\rho|})} = A_{\widetilde{C}_{(\rho, 1^{n-|\rho|})}} 
= \sum_{x\in \widetilde{C}_{(\rho, 1^{n-|\rho|})}} x \ \in 
\mathbb{C}[\widetilde{\mathfrak{S}}_n]. \label{eq:3-4} 
\end{align}
The case of $3$-cycles for \eqref{eq:3-4} already appeared in \eqref{eq:2-4}. 

For $\rho\in \overline{\mathcal{OP}}$, \eqref{eq:3-2} does not depend on the order of the 
product of cycles. 
Indeed, if $i_1, \cdots, i_p$, $j_1, \cdots, j_q$ ($p, q\geqq 2$) are distinct, we see 
\begin{equation}\label{eq:3-5}
[i_1 \cdots i_p] [j_1 \cdots j_q] = z^{(p-1)(q-1)} 
[j_1 \cdots j_q] [i_1 \cdots i_p].
\end{equation}
Since 
$[l\ k_1] [k_1\ k_2 \cdots k_p] [l\ k_1] = z^{p-1} [l\ k_2 \cdots k_p]$ 
holds if $k_1, \cdots, k_p, l$ are distinct, odd length cycles with common length 
are mutually conjugate: $[k_1 \cdots k_p] \sim [l_1 \cdots l_p]$. 
We see from \eqref{eq:1-26} 
\begin{align*}
&\rho\in\overline{\mathcal{OP}} \ \Longrightarrow \ z\widetilde{C}_{(\rho,1^{n-|\rho|})}
\text{ and } \widetilde{C}_{(\rho,1^{n-|\rho|})} 
\text{ are different conjugacy classes}, \\ 
&\rho\in \overline{\mathcal{P}}\setminus\overline{\mathcal{OP}} \text{ and } n\geqq |\rho|+2 \ 
\Longrightarrow \ 
z\widetilde{C}_{(\rho,1^{n-|\rho|})} = \widetilde{C}_{(\rho,1^{n-|\rho|})}. 
\end{align*}

In the following Lemma~\ref{lem:3-2}, Lemma~\ref{lem:3-3} and 
Lemma~\ref{lem:3-4}, let $f^{(n)}$ be a spin normalized central positive-definite 
function on $\widetilde{\mathfrak{S}}_n$ and assume that their sequence 
$\{ f^{(n)}\}$ satisfies 
approximate factorization property \eqref{eq:1-29} and condition \eqref{eq:1-30}. 
Also let $\rho, \rho^\prime \in \overline{\mathcal{OP}}$ in \eqref{eq:3-1}. 

\begin{lemma}\label{lem:3-2}
We have 
\begin{equation}\label{eq:3-8}
f^{(n)}([ \rho, 1^{n-|\rho|}] ) = O( n^{-\frac{|\rho|-l(\rho)}{2}}) \qquad 
(n\to\infty)
\end{equation}
for each $\rho$. 
\hfill $\square$
\end{lemma}

\textit{Proof} \ 
If $l(\rho) =2$ and $\rho = (\rho_1, \rho_2)$, we have 
\begin{align}
f^{(n)} ([\rho_1, \rho_2, 1^{n-\rho_1-\rho_2}]) &= 
f^{(n)} ([\rho_1, 1^{n-\rho_1}]) f^{(n)} ([\rho_2, 1^{n-\rho_2}]) + 
o( n^{-\frac{1}{2}(\rho_1-1+\rho_2-1)}) \notag \\ 
&= O( n^{-\frac{\rho_1-1}{2}}) O (n^{-\frac{\rho_2-1}{2}}) + 
o(n^{-\frac{1}{2}(\rho_1+\rho_2-2)}) = 
O(n^{-\frac{1}{2}(\rho_1+\rho_2-2)}).  
\label{eq:3-9}
\end{align}
The first and second equalities comes from \eqref{eq:1-29} and \eqref{eq:1-30} 
respectively. 
This implies \eqref{eq:3-8}. 
Inductively, letting $\rho = (\rho_1)\sqcup \rho^0$, we have 
\begin{align*}
f^{(n)}([\rho, 1^{n-|\rho|}]) &= 
f^{(n)}([\rho_1, 1^{n-\rho_1}]) f^{(n)}([\rho^0, 1^{n-|\rho^0|}]) + 
o( n^{-\frac{1}{2}(\rho_1-1+|\rho^0|-l(\rho^0))})  \\ 
&= O(n^{-\frac{\rho_1-1}{2}}) O( n^{-\frac{|\rho^0|-l(\rho^0)}{2}}) + 
o(n^{-\frac{|\rho|-l(\rho)}{2}}) = O(n^{-\frac{|\rho|-l(\rho)}{2}}). 
\end{align*}
This completes the proof. 
\hfill $\blacksquare$

\begin{lemma}\label{lem:3-3}
We have 
\begin{equation}\label{eq:3-11}
f^{(n)} \Bigl( 
\frac{\widetilde{A}_{(\rho, 1^{n-|\rho|})}}{|\widetilde{C}_{(\rho, 1^{n-|\rho|})}|} 
\frac{\widetilde{A}_{(\rho^\prime, 1^{n-|\rho^\prime|})}}
{|\widetilde{C}_{(\rho^\prime, 1^{n-|\rho^\prime|})}|} \Bigr) - 
f^{(n)} \bigl( [\rho\sqcup\rho^\prime, 1^{n-|\rho|-|\rho^\prime|}] \bigr) = 
O \bigl( n^{-\frac{|\rho|-l(\rho)+|\rho^\prime|-l(\rho^\prime)}{2}-1} \bigr)
\end{equation}
for each $\rho, \rho^\prime$. 
\hfill $\square$
\end{lemma}

\textit{Proof} \ 
The size of a conjugacy class of type $\mathcal{OP}_n$ is the same in $\widetilde{\mathfrak{S}}_n$ 
and $\mathfrak{S}_n$: 
\begin{equation}\label{eq:3-12}
|\widetilde{C}_{(\rho, 1^{n-|\rho|})}| = \frac{n^{\downarrow |\rho|}}{z_\rho}, 
\qquad z_\rho = \prod_j j^{m_j(\rho)} m_j(\rho)!.
\end{equation}
Set 
\begin{equation}\label{eq:3-14}
S^{(n)}_{\rho,\rho^\prime}(r) = \bigl\{ (x, y)\in 
\widetilde{C}_{(\rho, 1^{n-|\rho|})} \times 
\widetilde{C}_{(\rho^\prime, 1^{n-|\rho^\prime|})} \,\big|\, 
|(\mathrm{supp}\, x) \cap (\mathrm{supp}\,y) =r \bigr\} 
\end{equation}
for $r\in\{ 0,1,2,\cdots\}$, and take decomposition as 
\begin{equation}\label{eq:3-15}
\widetilde{C}_{(\rho, 1^{n-|\rho|})} \times 
\widetilde{C}_{(\rho^\prime, 1^{n-|\rho^\prime|})} = 
\bigsqcup_{r=0}^{|\rho|\wedge |\rho^\prime|} 
S^{(n)}_{\rho,\rho^\prime}(r).  
\end{equation}
We get 
\begin{equation}\label{eq:3-16}
\widetilde{A}_{(\rho, 1^{n-|\rho|})} 
\widetilde{A}_{(\rho^\prime, 1^{n-|\rho^\prime|})} = 
\sum_{r=0}^{|\rho|\wedge |\rho^\prime|} 
\sum_{(x, y)\in S^{(n)}_{\rho,\rho^\prime}(r)} xy.
\end{equation}
The terms of $r=0$ in the right hand side of \eqref{eq:3-16} is 
\begin{equation}\label{eq:3-17}
\sum_{(x, y)\in S^{(n)}_{\rho,\rho^\prime}(0)} xy = 
\frac{z_{\rho\cup\rho^\prime}}{z_\rho z_{\rho^\prime}} 
\widetilde{A}_{(\rho\cup\rho^\prime, 1^{n-|\rho|-|\rho^\prime|})}. 
\end{equation}
Taking $f^{(n)}$-values of \eqref{eq:3-16}, we get from 
\eqref{eq:3-12} and \eqref{eq:3-17} 
\begin{align}
f^{(n)} \Bigl( \frac{\widetilde{A}_{(\rho, 1^{n-|\rho|})} 
\widetilde{A}_{(\rho^\prime, 1^{n-|\rho^\prime|})}}
{|\widetilde{C}_{(\rho, 1^{n-|\rho|})}| 
|\widetilde{C}_{(\rho^\prime, 1^{n-|\rho^\prime|})}|} \Bigr) 
= &\bigl( 1+O(\frac{1}{n})\bigr) f^{(n)} 
([\rho\cup\rho^\prime, 1^{n-|\rho|-|\rho^\prime|}]) 
\label{eq:3-18} \\ 
&+ \frac{O(1)}{n^{|\rho|+|\rho^\prime|}} 
\sum_{r=1}^{|\rho|\wedge |\rho^\prime|} 
\sum_{x,y\in S^{(n)}_{\rho,\rho^\prime}(r)} f^{(n)}(xy). 
\label{eq:3-19}
\end{align}
The support and type of an element of $\widetilde{\mathfrak{S}}_n$ are 
defined by projecting onto $\mathfrak{S}_n$ by $\Phi$. 
If $x, y\in S^{(n)}_{\rho,\rho^\prime}(r)$ ($r\geqq 1$) and 
$\mathrm{type}(xy)$ is $(\sigma, 1^{n-|\sigma|})$, we get 
\begin{equation}\label{eq:3-20}
|\sigma| -l(\sigma) \geqq |\rho|-l(\rho) +|\rho^\prime|
-l(\rho^\prime) -2r+2 
\end{equation}
by using \cite[Lemma~3.5]{Hor05} holding in $\mathfrak{S}_n$. 
Lemma~\ref{lem:3-2} and \eqref{eq:3-20} yield 
\begin{equation}\label{eq:3-21}
x, y\in S^{(n)}_{\rho, \rho^\prime}(r), \ r\geqq 1 
\ \Longrightarrow\ 
f^{(n)}(xy) = O\bigl( n^{-\frac{1}{2}(|\rho|-l(\rho)+|\rho^\prime|
-l(\rho^\prime)-2r+2)} \bigr).
\end{equation}
It is obvious that 
\begin{equation}\label{eq:3-22}
|S^{(n)}_{\rho, \rho^\prime}(r)| = O(n^{|\rho|+|\rho^\prime|-r}).
\end{equation}
By \eqref{eq:3-21} and \eqref{eq:3-22}, we continue the terms \eqref{eq:3-19} as 
\begin{equation}\label{eq:3-23}
\frac{O(1)}{n^{|\rho|+|\rho^\prime|}} O\bigl( 
n^{|\rho|+|\rho^\prime|-r- \frac{1}{2}(|\rho|-l(\rho)
+|\rho^\prime|-l(\rho^\prime)) +r-1} \bigr) = 
O \bigl( n^{-\frac{1}{2}(|\rho|-l(\rho)+|\rho^\prime|-l(\rho^\prime))-1}
\bigr).
\end{equation}
Applying \eqref{eq:3-8} also to \eqref{eq:3-18} and combining \eqref{eq:3-23}, 
we obtain \eqref{eq:3-11}. 
\hfill $\blacksquare$

\begin{lemma}\label{lem:3-4} \ 
We have 
\[ 
f^{(n)} \Bigl( \frac{\widetilde{A}_{(\rho, 1^{n-|\rho|})}}
{|\widetilde{C}_{(\rho, 1^{n-|\rho|})}|} 
\frac{\widetilde{A}_{(\rho^\prime, 1^{n-|\rho^\prime|})}}
{|\widetilde{C}_{(\rho^\prime, 1^{n-|\rho^\prime|})}|} \Bigr) - 
f^{(n)} ([\rho, 1^{n-|\rho|}]) 
f^{(n)} ([\rho^\prime, 1^{n-|\rho^\prime|}]) = 
o\bigl( n^{-\frac{1}{2}(|\rho|-l(\rho)+|\rho^\prime|-l(\rho^\prime))}
\bigr)
\] 
for $\rho$ and $\rho^\prime$. 
\hfill $\square$
\end{lemma}

\textit{Proof} \ 
The assertion follows from approximate factorization property \eqref{eq:1-29} and 
Lemma~\ref{lem:3-3}. 
\hfill $\blacksquare$

\medskip

Let us prepare some notations. 
The canonical homomorphism of \eqref{eq:1-14} is written as 
\ $\Phi_n: \mathbb{C}[\widetilde{\mathfrak{S}}_n] \longrightarrow 
\mathbb{C}[\mathfrak{S}_n]$ \ 
with explicit subscript $n$. 
Consider the restriction map 
\ $E_n : \mathbb{C}[\mathfrak{S}_{n+1}] \longrightarrow 
\mathbb{C}[\mathfrak{S}_n]$ \ 
as in \eqref{eq:1-22}. 
The $\Phi$ and $E$ satisfies 
\begin{equation}\label{eq:3-28}
\Phi_n \widetilde{E}_n = E_n \Phi_{n+1}. 
\end{equation}
For $\sigma\in \overline{\mathcal{P}}$, let $\sigma^\circ$ denote the element of 
$\overline{\mathcal{P}}$ by removing the rows of length $2$ and the first column: 
\begin{equation}\label{eq:3-29}
\sigma = (2^{m_2}3^{m_3}4^{m_4}\cdots ) \ \longmapsto \ 
\sigma^\circ = (2^{m_3}3^{m_4}\cdots )
\end{equation}
as in Figure~\ref{fig:3-1}.

\begin{figure}[hbt]
\centering 
\include{sigma-maru-fig}
\vspace{-5mm}
\caption{$\sigma = (4,3,2,2)$ and $\sigma^\circ =(3,2)$}
\label{fig:3-1}
\end{figure}

\noindent 
For $\sigma\in \mathcal{P}_n$, let $NC(\sigma)$ denote the set of non-crossing partitions 
of $\{1,2, \cdots, n\}$ with type $\sigma$. 
Figure~\ref{fig:3-2} shows an example. 

\begin{figure}[hbt]
\centering 
\include{NC32-fig}
\vspace{-5mm}
\caption{The elements of $NC(\sigma)$, $\sigma =(3,2)$}
\label{fig:3-2}
\end{figure}

\begin{proposition}\label{prop:3-5} \ 
Assume the sequence $\{f^{(n)}\}$ in Theorem~\ref{th:static} satisfies 
approximate factorization property \eqref{eq:1-29} and \eqref{eq:1-30}. 
Then we have for $\forall k\in\mathbb{N}$ 
\begin{align}
&\lim_{n\to\infty}f^{(n)} (n^{-k} \widetilde{E}_n \widetilde{J}_{n+1}^{\;2k})
= m_{2k}, \label{eq:3-30} \\ 
&f^{(n)} (\widetilde{E}_n \widetilde{J}_{n+1}^{\;2k}\cdot 
\widetilde{E}_n \widetilde{J}_{n+1}^{\;2k}) - 
f^{(n)} (\widetilde{E}_n \widetilde{J}_{n+1}^{\;2k})^2 
= o(n^{2k}) \label{eq:3-31}
\end{align}
and hence 
\begin{equation}\label{eq:3-32}
\lim_{n\to\infty} f^{(n)} 
\bigl( (n^{-k} \widetilde{E}_n \widetilde{J}_{n+1}^{\;2k}- m_{2k})^2 
\bigr) =0
\end{equation}
also. 
\hfill $\square$
\end{proposition}

\textit{Proof} \ 
First we show \eqref{eq:3-31}. 
Since $\widetilde{E}_n \widetilde{J}_{n+1}^{\;2k}$ belongs to the center of 
$\mathbb{C}[\widetilde{\mathfrak{S}}_n]$ as seen from \eqref{eq:2-4}, 
it is expressed as a linear combination of 
$A_{\widetilde{C}}$ ($\widetilde{C}$ being a conjugacy class in 
$\widetilde{\mathfrak{S}}_n$). 
Let $\alpha_{\widetilde{C}}$ be its coefficient: 
\begin{equation}\label{eq:3-33}
\widetilde{E}_n \widetilde{J}_{n+1}^{\;2k} = 
\sum_{\widetilde{C}} \alpha_{\widetilde{C}} A_{\widetilde{C}}.
\end{equation}
Projected by $\widetilde{E}_n$, each term in the expansion of 
$\widetilde{J}_{n+1}^{\;2k}$ falls into one of the conjugacy classes unless trivial. 
This implies $\alpha_{\widetilde{C}} \geqq 0$ (like a histogram on the 
conjugacy classes). 
In a situation of fixing $k$ and letting $n$ tend to $\infty$, 
the type of an element of conjugacy class $\widetilde{C}$ appearing in 
\eqref{eq:3-33} is of bounded width and hence does not belong to $\mathcal{SP}_n$. 
Therefore, only those terms of types $\mathcal{OP}_n$ in \eqref{eq:3-33} can survive 
after taking $f^{(n)}$-values. 
We have from \eqref{eq:3-33} 
\begin{equation}\label{eq:3-34}
\widetilde{E}_n \widetilde{J}_{n+1}^{\;2k}\cdot 
\widetilde{E}_n \widetilde{J}_{n+1}^{\;2k} = 
\sum_{\widetilde{C}, \widetilde{C}^\prime} 
\alpha_{\widetilde{C}} \alpha_{\widetilde{C}^\prime} 
A_{\widetilde{C}} A_{\widetilde{C}^\prime}.
\end{equation}
In \eqref{eq:3-34}, 
if $n$ is large enough, the elements of $\widetilde{C}$ and $\widetilde{C}^\prime$ 
are never of type $\mathcal{SP}_n$. 
If $\widetilde{C}$ is not of type $\mathcal{OP}_n$, then 
since $z A_{\widetilde{C}} = A_{\widetilde{C}}$ holds, we have 
\begin{equation}\label{eq:3-35}
f^{(n)} ( A_{\widetilde{C}} A_{\widetilde{C}^\prime}) = 
f^{(n)} (z A_{\widetilde{C}} A_{\widetilde{C}^\prime}) = 
- f^{(n)} (A_{\widetilde{C}} A_{\widetilde{C}^\prime}), \quad 
\text{hence} \ 
f^{(n)}(A_{\widetilde{C}} A_{\widetilde{C}^\prime})=0.
\end{equation}
If $\widetilde{C}^\prime$ is not of type $\mathcal{OP}_n$, \eqref{eq:3-35} holds similarly. 
Taking $f^{(n)}$-values of \eqref{eq:3-33} and \eqref{eq:3-34}, we get 
\begin{multline}\label{eq:3-36}
f^{(n)}( \widetilde{E}_n \widetilde{J}_{n+1}^{\;2k}\cdot 
\widetilde{E}_n \widetilde{J}_{n+1}^{\;2k}) - 
f^{(n)}( \widetilde{E}_n \widetilde{J}_{n+1}^{\;2k})^2 \\ 
= \sum_{\widetilde{C}, \widetilde{C}^\prime : \text{type }\mathcal{OP}_n} 
\alpha_{\widetilde{C}} \alpha_{\widetilde{C}^\prime} 
\bigl\{ f^{(n)} ( A_{\widetilde{C}} A_{\widetilde{C}^\prime}) - 
f^{(n)} ( A_{\widetilde{C}}) f^{(n)}( A_{\widetilde{C}^\prime}) \bigr\}.
\end{multline}
Conjugacy classes of type $\mathcal{OP}_n$ are $\widetilde{C}_{(\rho, 1^{n-|\rho|})}$ 
and $z \widetilde{C}_{(\rho, 1^{n-|\rho|})}$. 
Since $A_{z\widetilde{C}} = z A_{\widetilde{C}}$ holds, we continue 
\eqref{eq:3-36} as 
\begin{align}
&= \sum_{\rho, \rho^\prime\in \overline{\mathcal{OP}}} 
\Bigl[ \alpha_{\widetilde{C}_{(\rho,1^{n-|\rho|})}} 
\alpha_{\widetilde{C}_{(\rho^\prime,1^{n-|\rho^\prime|})}} 
\bigl\{ f^{(n)} ( \widetilde{A}_{(\rho,1^{n-|\rho|})} 
\widetilde{A}_{(\rho^\prime,1^{n-|\rho^\prime|})}) - 
f^{(n)}(\widetilde{A}_{(\rho,1^{n-|\rho|})}) 
f^{(n)}(\widetilde{A}_{(\rho^\prime,1^{n-|\rho^\prime|})}\bigr\} 
\notag \\ 
&+ \alpha_{z \widetilde{C}_{(\rho,1^{n-|\rho|})}} 
\alpha_{\widetilde{C}_{(\rho^\prime,1^{n-|\rho^\prime|})}} 
\bigl\{ f^{(n)} (z \widetilde{A}_{(\rho,1^{n-|\rho|})} 
\widetilde{A}_{(\rho^\prime,1^{n-|\rho^\prime|})}) - 
f^{(n)}(z \widetilde{A}_{(\rho,1^{n-|\rho|})}) 
f^{(n)}(\widetilde{A}_{(\rho^\prime,1^{n-|\rho^\prime|})}\bigr\} 
\notag \\ 
&+ \alpha_{\widetilde{C}_{(\rho,1^{n-|\rho|})}} 
\alpha_{z \widetilde{C}_{(\rho^\prime,1^{n-|\rho^\prime|})}} 
\bigl\{ f^{(n)} ( \widetilde{A}_{(\rho,1^{n-|\rho|})} 
z \widetilde{A}_{(\rho^\prime,1^{n-|\rho^\prime|})}) - 
f^{(n)}(\widetilde{A}_{(\rho,1^{n-|\rho|})}) 
f^{(n)}(z \widetilde{A}_{(\rho^\prime,1^{n-|\rho^\prime|})}\bigr\} 
\notag \\ 
&+ \alpha_{z \widetilde{C}_{(\rho,1^{n-|\rho|})}} 
\alpha_{z \widetilde{C}_{(\rho^\prime,1^{n-|\rho^\prime|})}} 
\bigl\{ f^{(n)} ( z\widetilde{A}_{(\rho,1^{n-|\rho|})} 
z\widetilde{A}_{(\rho^\prime,1^{n-|\rho^\prime|})}) - 
f^{(n)}(z\widetilde{A}_{(\rho,1^{n-|\rho|})}) 
f^{(n)}(z\widetilde{A}_{(\rho^\prime,1^{n-|\rho^\prime|})}\bigr\} 
\Bigr] \notag \\ 
&= \sum_{\rho, \rho^\prime\in \overline{\mathcal{OP}}} 
(\alpha_{\widetilde{C}_{(\rho,1^{n-|\rho|})}} - 
\alpha_{z \widetilde{C}_{(\rho,1^{n-|\rho|})}}) 
(\alpha_{\widetilde{C}_{(\rho^\prime,1^{n-|\rho^\prime|})}} - 
\alpha_{z \widetilde{C}_{(\rho^\prime,1^{n-|\rho^\prime|})}}) 
\notag \\ 
&\qquad\qquad \cdot \bigl\{ f^{(n)} ( \widetilde{A}_{(\rho,1^{n-|\rho|})} 
\widetilde{A}_{(\rho^\prime,1^{n-|\rho^\prime|})}) - 
f^{(n)}(\widetilde{A}_{(\rho,1^{n-|\rho|})}) 
f^{(n)}(\widetilde{A}_{(\rho^\prime,1^{n-|\rho^\prime|})}\bigr\}. 
\label{eq:3-37}
\end{align}
Here the coefficients $\alpha_\cdot$ in \eqref{eq:3-33} enjoy the estimate of 
Lemma~\ref{lem:3-6} below. 
If $k$ is fixed, the range $\rho$ and $\rho^\prime$ run over in \eqref{eq:3-37} 
does not depend on $n$. 
Combining these with Lemma~\ref{lem:3-4}, we get 
\[ 
|\text{\eqref{eq:3-37}}| = 
O(1) n^{k+\frac{|\rho|-l(\rho)}{2}} 
n^{k+\frac{|\rho^\prime|-l(\rho^\prime)}{2}} 
o\bigl( n^{-\frac{1}{2}(|\rho|-l(\rho)+|\rho^\prime|-l(\rho^\prime))}
\bigr) = o(n^{2k}).
\] 
This implies \eqref{eq:3-31}. 
\hfill $\blacksquare$

\begin{lemma}\label{lem:3-6}
Nontrivial coefficients in \eqref{eq:3-33} have the following estimate. 
To the type $\rho\in\overline{\mathcal{OP}}$ of each $\widetilde{C}$ is assigned a unique 
$\sigma\in\overline{\mathcal{P}}_{2k}$ such that $\rho=\sigma^\circ$ 
($\,^\circ$ defined in \eqref{eq:3-29}) and 
\begin{equation}\label{eq:3-39}
\alpha_{\widetilde{C}_{(\rho,1^{n-|\rho|})}} 
|\widetilde{C}_{(\rho,1^{n-|\rho|})}| \leqq 
|NC(\sigma)| \bigl( 1+ O(\frac{1}{n})\bigr) 
n^{k+\frac{|\rho|-l(\rho)}{2}}.
\end{equation}
The conjugacy class in the left hand side can be replaced by 
$z\widetilde{C}_{(\rho,1^{n-|\rho|})}$. 
\hfill $\square$
\end{lemma}

\textit{Proof} \ 
Let $\Phi_n$ operate on \eqref{eq:3-33}: 
\begin{equation}\label{eq:3-40}
\Phi_n \widetilde{E}_n \widetilde{J}_{n+1}^{\;2k} = 
\sum_{\widetilde{C}} \alpha_{\widetilde{C}} \Phi_n(A_{\widetilde{C}}).
\end{equation}
By \eqref{eq:3-28}, the left hand side of \eqref{eq:3-40} is equal to 
\begin{equation}\label{eq:3-41}
E_n \Phi_{n+1} (\widetilde{J}_{n+1}^{\;2k}) = 
E_n ( \Phi_{n+1}(\widetilde{J}_{n+1})^{2k}) = 
E_n J_{n+1}^{\;2k}. 
\end{equation}
Since a conjugacy class of $\widetilde{\mathfrak{S}}_n$ is split or non-split, 
the right hand side of \eqref{eq:3-40} is 
\begin{align}
&= \sum_{\widetilde{C}: \text{ split}} 
\alpha_{\widetilde{C}} \Phi_n(A_{\widetilde{C}}) + 
\sum_{\widetilde{C}: \text{ non-split}} 
\alpha_{\widetilde{C}} \Phi_n(A_{\widetilde{C}}) \notag \\ 
&= \sum_{\rho\in \overline{\mathcal{OP}}} (
\alpha_{\widetilde{C}_{(\rho,1^{n-|\rho|})}} + 
\alpha_{z \widetilde{C}_{(\rho,1^{n-|\rho|})}} ) 
A_{(\rho,1^{n-|\rho|})} + 
\sum_{\nu\in\overline{\mathcal{P}}\setminus\overline{\mathcal{OP}}} 
\alpha_{\widetilde{C}_{(\nu,1^{n-|\nu|})}} \cdot 
2 A_{(\nu,1^{n-|\nu|})}
\label{eq:3-42}
\end{align}
where $A_{(\rho,1^{n-|\rho|})}, \ A_{(\nu,1^{n-|\nu|})}\in
\mathbb{C}[\mathfrak{S}_n]$. 
By using \cite[Theorem~1]{Hor06}
\footnote{Caution: The notation \lq\lq $l(\sigma)$\rq\rq\ used in \cite{Hor06} means 
$|\sigma|-l(\sigma)$.}, 
\eqref{eq:3-41} continues as
\begin{equation}\label{eq:3-43}
E_n J_{n+1}^{\;2k} = \sum_{\sigma\in\overline{\mathcal{P}}_{2k}} |NC(\sigma)| 
n^{|\sigma|-l(\sigma)} \bigl( 1+O(\frac{1}{n})\bigr) 
\frac{A_{(\sigma^\circ,1^{n-|\sigma^\circ|})}}
{|C_{(\sigma^\circ,1^{n-|\sigma^\circ|})}|}. 
\end{equation}
Dividing the sum of \eqref{eq:3-43} into two ones according as 
$\sigma^\circ\in \overline{\mathcal{OP}}$ 
or $\sigma^\circ\in \overline{\mathcal{P}}\setminus\overline{\mathcal{OP}}$, we have 
\[ 
\sigma^\circ\in \overline{\mathcal{OP}} \iff 
\text{All parts of $\sigma $ have even length}.
\] 
Compare the both sides of \eqref{eq:3-42}$=$\eqref{eq:3-43}. 
We see that a unique $\sigma\in \overline{\mathcal{P}}_{2k}$ such that $\sigma^\circ =\rho$ 
is determined for each $\rho\in\overline{\mathcal{OP}}$ and have 
\begin{equation}\label{eq:3-45}
(\alpha_{\widetilde{C}_{(\rho,1^{n-|\rho|})}} + 
\alpha_{z \widetilde{C}_{(\rho,1^{n-|\rho|})}} ) 
|C_{(\rho,1^{n-|\rho|})}| = 
|NC(\sigma)| n^{|\sigma|-l(\sigma)} \bigl( 1+O(\frac{1}{n})\bigr).
\end{equation}
Here the relation $\rho=\sigma^\circ$ such that $\sigma\in \overline{\mathcal{P}}_{2k}$ 
yields 
\[ 
\sigma = (\rho_1\!+1, \ \cdots, \ \rho_{l(\rho)}\!+1, \ 2^m), 
\qquad m= \frac{2k-(|\rho|+l(\rho))}{2}.
\] 
Hence we have 
\begin{equation}\label{eq:3-47}
|\sigma|-l(\sigma) = 2k- \{ l(\rho)+\frac{1}{2}(2k-|\rho|-l(\rho))\} 
= k+ \frac{|\rho|-l(\rho)}{2}.
\end{equation}
We get desired \eqref{eq:3-39} from \eqref{eq:3-45} and \eqref{eq:3-47}. 
\hfill $\blacksquare$

\medskip

\textit{Proof of Proposition~\ref{prop:3-5} continued} \ 
Next we show \eqref{eq:3-30}. 
Let $n$ be large enough while $k$ is fixed. 
Taking $f^{(n)}$-values of \eqref{eq:3-33} and noting only those terms 
corresponding to conjugacy classes of type $\mathcal{OP}_n$ survive, we have 
\begin{equation}\label{eq:3-48}
f^{(n)} (\widetilde{E}_n \widetilde{J}_{n+1}^{\;2k}) = 
\sum_{\widetilde{C} : \text{ type } \mathcal{OP}_n} 
\alpha_{\widetilde{C}} f^{(n)}(A_{\widetilde{C}}). 
\end{equation}
Comparing the both sides of \eqref{eq:3-42}$=$\eqref{eq:3-43}, we see that 
conjugacy classes in the sum of \eqref{eq:3-48} run over the range determined by 
\begin{multline*}
\widetilde{C} = \widetilde{C}_{(\rho, 1^{n-|\rho|})} \quad \text{or} \quad 
z\widetilde{C}_{(\rho, 1^{n-|\rho|})} \quad \text{such that} \\ 
\rho=\sigma^\circ, \quad \sigma\in\overline{\mathcal{P}}_{2k}, \quad 
\text{and all rows of $\sigma$ have even length}.
\end{multline*}
Hence \eqref{eq:3-48} continues as 
\begin{align}
&= \sum_{\sigma\in\overline{\mathcal{P}}_{2k} : \text{ even length rows}} 
\bigl\{ \alpha_{\widetilde{C}_{(\sigma^\circ, 1^{n-|\sigma^\circ|})}} 
f^{(n)}( \widetilde{A}_{(\sigma^\circ, 1^{n-|\sigma^\circ|})}) + 
\alpha_{z \widetilde{C}_{(\sigma^\circ, 1^{n-|\sigma^\circ|})}} 
f^{(n)}( z \widetilde{A}_{(\sigma^\circ, 1^{n-|\sigma^\circ|})} \bigr\} 
\notag \\ 
&= \sum_{\sigma\in\overline{\mathcal{P}}_{2k} : \text{ even length rows}} 
(\alpha_{\widetilde{C}_{(\sigma^\circ, 1^{n-|\sigma^\circ|})}} - 
\alpha_{z \widetilde{C}_{(\sigma^\circ, 1^{n-|\sigma^\circ|})}} ) 
f^{(n)}( \widetilde{A}_{(\sigma^\circ, 1^{n-|\sigma^\circ|})}). 
\label{eq:3-50}
\end{align}
Since we show $\alpha_{\widetilde{C}_{(\sigma^\circ, 1^{n-|\sigma^\circ|})}}$ 
is overwhelmingly larger than 
$\alpha_{z \widetilde{C}_{(\sigma^\circ, 1^{n-|\sigma^\circ|})}}$ 
in Lemma~\ref{lem:3-7} below, \eqref{eq:3-50} continues as 
\begin{equation}\label{eq:3-51}
= \bigl( 1+O(\frac{1}{n})\bigr) 
\sum_{\sigma\in\overline{\mathcal{P}}_{2k} : \text{ even length rows}} 
\alpha_{\widetilde{C}_{(\sigma^\circ, 1^{n-|\sigma^\circ|})}} 
f^{(n)}( \widetilde{A}_{(\sigma^\circ, 1^{n-|\sigma^\circ|})}). 
\end{equation}
By virtue of \eqref{eq:3-45}, \eqref{eq:3-51} continues as 
\begin{equation}\label{eq:3-52}
= \bigl( 1+O(\frac{1}{n})\bigr) 
\sum_{\sigma\in\overline{\mathcal{P}}_{2k} : \text{ even length rows}} 
|NC(\sigma)| n^{|\sigma|-l(\sigma)} 
f^{(n)}([\sigma^\circ, 1^{n-|\sigma^\circ|}]).
\end{equation}
Since we have 
\[ 
f^{(n)}([\rho, 1^{n-|\rho|}]) = f^{(n)}([\rho_1, 1^{n-\rho_1}]) \cdots 
f^{(n)}([\rho_{l(\rho)}, 1^{n-\rho_{l(\rho)}}]) + 
o(n^{-\frac{|\rho|-l(\rho)}{2}}), \qquad \rho\in\overline{\mathcal{OP}} 
\] 
similarly to \eqref{eq:3-9}, we get from the assumption of \eqref{eq:1-30} 
\begin{align}
f^{(n)}([ \sigma^\circ, 1^{n-|\sigma^\circ|}]) 
&= ( 1+o(1) ) r_{\sigma_1} n^{-\frac{\sigma_1-2}{2}} 
r_{\sigma_2} n^{-\frac{\sigma_2-2}{2}} \cdots 
r_{\sigma_{l(\sigma^\circ)}} n^{-\frac{\sigma_{l(\sigma^\circ)}-2}{2}} 
+ o\bigl( n^{-\frac{|\sigma^\circ|-l(\sigma^\circ)}{2}} \bigr) 
\notag \\ 
&= (1+o(1)) r_{\sigma_1}\cdots r_{\sigma_{l(\sigma^\circ)}} 
n^{-\frac{|\sigma^\circ|-l(\sigma^\circ)}{2}} + 
o \bigl( n^{-\frac{|\sigma^\circ|-l(\sigma^\circ)}{2}} \bigr). 
\label{eq:3-54}
\end{align}
Putting \eqref{eq:3-54} into \eqref{eq:3-52} and noting 
\[ 
|\sigma|-l(\sigma) - \frac{|\sigma^\circ|-l(\sigma^\circ)}{2} = 
2k- l(\sigma) - \frac{2k-2 l(\sigma)}{2} = k, 
\] 
we get 
\begin{equation}\label{eq:3-56}
f^{(n)}( \widetilde{E}_n \widetilde{J}_{n+1}^{\;2k}) = (1+o(1)) 
\sum_{\sigma\in\overline{\mathcal{P}}_{2k} : \text{ even length rows}} 
|NC(\sigma)| r_{\sigma_1} \cdots r_{\sigma_{l(\sigma^\circ)}} n^k 
+ o(n^k).
\end{equation}
Since \ 
$\sigma_{l(\sigma^\circ)+1} = \cdots = \sigma_{l(\sigma)} = 2$ \ 
holds, \eqref{eq:3-56} yields 
\begin{equation}\label{eq:3-58}
\lim_{n\to\infty} n^{-k} f^{(n)}( \widetilde{E}_n \widetilde{J}_{n+1}^{\;2k}) 
= \sum_{\sigma\in \mathcal{P}_{2k}} |NC(\sigma)| r_{\sigma_1}\cdots 
r_{\sigma_{l(\sigma)}}.
\end{equation}
The free cumulant-moment formula tells us that the right hand side of 
\eqref{eq:3-58} agrees with $m_{2k}$. 
This shows \eqref{eq:3-30} (modulo Lemma~\ref{lem:3-7}) and thus 
completes the proof of Proposition~\ref{prop:3-5}. 
\hfill $\blacksquare$

\begin{lemma}\label{lem:3-7} 
For $\sigma^\circ\in\overline{\mathcal{OP}}$ such that 
$\sigma\in\overline{\mathcal{P}}_{2k}$ in Lemma~\ref{lem:3-6}, 
\begin{equation}\label{eq:3-59}
\alpha_{z\widetilde{C}_{(\sigma^\circ, 1^{n-|\sigma^\circ|})}} / 
\alpha_{\widetilde{C}_{(\sigma^\circ, 1^{n-|\sigma^\circ|})}} 
= O(\frac{1}{n})
\end{equation}
holds. 
\hfill $\square$
\end{lemma}

\textit{Proof} \ 
The $\alpha_\cdot$'s are expanding coefficients in \eqref{eq:3-33}. 
Setting $\ast =n+1$, let us consider 
\[ 
\widetilde{J}_{n+1}^{\;2k} = 
\sum_{i_1, \cdots, i_{2k}\in\{ 1,\cdots, n\}} [i_{2k}\ \ast] 
[i_{2k-1}\ \ast] \cdots [i_2\ \ast][i_1\ \ast]. 
\] 
Counting the number of those $2k$-walks 
\begin{equation}\label{eq:3-61}
e \to [i_1\ \ast] \to [i_2\ \ast] [i_1\ \ast] \to \cdots\cdots \to 
[i_{2k}\ \ast] \cdots [i_1\ \ast] 
\end{equation}
whose end points fall into the conjugacy classes 
$\widetilde{C}_{(\sigma^\circ, 1^{n-|\sigma^\circ|})}$ and 
$z \widetilde{C}_{(\sigma^\circ, 1^{n-|\sigma^\circ|})}$ of 
$\widetilde{\mathfrak{S}}_n$, and dividing them by 
$|\widetilde{C}_{(\sigma^\circ, 1^{n-|\sigma^\circ|})}|$, 
we get 
$\alpha_{\widetilde{C}_{(\sigma^\circ, 1^{n-|\sigma^\circ|})}}$ and 
$\alpha_{z \widetilde{C}_{(\sigma^\circ, 1^{n-|\sigma^\circ|})}}$ 
respectively. 
Project a $2k$-walk of \eqref{eq:3-61} by $\Phi$, and further project it onto 
its type. 
Then we have a $2k$-walk 
\begin{equation}\label{eq:3-62}
\varnothing \to (2) \to \cdots \to \sigma^\circ 
\end{equation}
starting at $\varnothing$ and terminating at $\sigma^\circ$ on a graph with 
$\overline{\mathcal{P}}$ as its vertex set stratified according to the values of 
$|\rho|-l(\rho)$ ($\rho\in\overline{\mathcal{P}}$)
\footnote{In \cite{Hor06} we called this graph a modified Young graph.}. 
In \eqref{eq:3-62}, let us refer to an edge $\to$ along which $|\rho|-l(\rho)$ 
increases [resp. decreases] as an up [resp. down] step. 
Taking into account \eqref{eq:3-61} before falling into \eqref{eq:3-62}, 
we see the up steps are divided into two classes: 
(a) up step of large degree of freedom (order of $n$) \ 
(b) up step of small degree of freedom (order of $1$) at which 
\begin{align*}
&\text{(a) take new $j$ outside the letters already used and multiply $[j\ \ast]$}  \\ 
&\text{(b) multiply one of the other $[j\ \ast]$'s} 
\end{align*}
In this sense a down step is always of small degree of freedom. 
For each $\rho\in\overline{\mathcal{P}}$, though there may be several $2k$-walks terminating 
at $\rho$, the numbers of up steps and down steps along a walk are uniquely determined
\footnote{Setting them as $u$ and $d$ respectively, we have 
$u+d=2k$, $u-d= |\rho|-l(\rho)$.}.  
We can reach $\rho$ even if we use as up steps only those of large degree of freedom. 
Once a cycle ($\in\widetilde{\mathfrak{S}}_n$) not containing $\ast$ is produced 
on the way of a $2k$-walk of \eqref{eq:3-61}, we have only to count those walks 
which never touch the produced cycle (of odd length) afterwards. 
An odd length cycle commutes with cyles which do not share common letters by 
\eqref{eq:3-5}. 

Let $[j\ \ast]$ be multiplied from the left as in \eqref{eq:3-61}. 
First at an up step, since $j$ is not contained in the letters already used, a cycle 
(of odd length) not containing $\ast$ commutes with $[j\ \ast]$, and hence 
\begin{equation}\label{eq:3-64}
[j\ \ast] [i_1\ \cdots\ i_p\ \ast] = [i_1\ \cdots\ i_p\ j\ \ast] 
\qquad (\forall p\in\{0,1,\cdots\}). 
\end{equation}
In other words, after an up step, the cycle containing $\ast$ grows by $1$ and 
contains $\ast$ at the right end. 
(If $p=0$, simply $[j\ \ast]$ is produced.) 
Next at a down step, we have 
\begin{equation}\label{eq:3-65}
[j\ \ast] [i_1\ \cdots\ i_p\ j\ i_{p+1}\ \cdots\ i_{p+q}\ \ast] = 
z^{pq} [i_1\ \cdots\ i_p\ \ast] [j\ i_{p+1}\ \cdots\ i_{p+q}]
\end{equation}
for $p, q\in\{0,1,2,\cdots\}$. 
Since we have only to count $[j\ i_{p+1}\ \cdots\ i_{p+q}]$ in the right hand 
side of \eqref{eq:3-65} which survives as an odd length cycle, we let $q$ be even 
and hence $z^{pq}=e$. 
Hence a down step gives rise to 
\[ 
[i_1\ \cdots\ i_p\ \ast] [ \text{odd length cycle} ].
\] 
When a walk reaches $\sigma^\circ$ after taking \eqref{eq:3-64} and \eqref{eq:3-65} 
by several steps, the final step must be \eqref{eq:3-65} of $p=0$ and produces 
an element of $\widetilde{\mathfrak{S}}_n$
\footnote{Of course this step may happen on the way besides the final step.}. 
The resulting element $x$ of type $\sigma^\circ$ is product of odd length cycles: 
\[ 
x = [1_1\ \cdots\ i_{\sigma^\circ_1}] [j_1\ \cdots\ j_{\sigma^\circ_2}] 
[k_1\ \cdots\ k_{\sigma^\circ_3}] \cdots. 
\] 
Since this element $x$ is conjugate with $[\sigma^\circ, 1^{n-|\sigma^\circ|}]$ in 
$\widetilde{\mathfrak{S}}_n$, it belongs to 
$\widetilde{C}_{(\sigma^\circ, 1^{n-|\sigma^\circ|})}$. 
Consequently, if we have ever a $2k$-walk of \eqref{eq:3-61} terminating at an 
element of $z \widetilde{C}_{(\sigma^\circ, 1^{n-|\sigma^\circ|})}$, the walk 
necessarily contains an up step of small degree of freedom. 
This implies that the ratio satisfies \eqref{eq:3-59}. 
\hfill $\blacksquare$

\medskip

Let us express the result of Proposition~\ref{prop:3-5} in terms of probabilitiy 
$M^{(n)}$ on $(\widetilde{\mathfrak{S}}_n)^\wedge_{\mathrm{spin}}$ 
corresponding to $f^{(n)}$. 
We have 
\begin{equation}\label{eq:3-68}
f^{(n)}(x) = \sum_{(\lambda,\gamma)\in 
(\widetilde{\mathfrak{S}}_n)^\wedge_{\mathrm{spin}}} 
M^{(n)}\bigl( \{(\lambda,\gamma)\}\bigr) 
\frac{\chi^{(\lambda,\gamma)}(x)}{\dim (\lambda,\gamma)}
\end{equation}
by \eqref{eq:1-25}. 
Let $\mathbb{E}^{(n)}$ denote the expectation of a random variable with respect to $M^{(n)}$: 
\[ 
\mathbb{E}^{(n)} [X] = \sum_{(\lambda,\gamma) \in 
(\widetilde{\mathfrak{S}}_n)^\wedge_{\mathrm{spin}}} X(\lambda, \gamma) 
M^{(n)}(\{(\lambda,\gamma)\}).
\] 
For $\lambda\in \mathcal{SP}_n$, we pick up the nonzero coordinates of valleys of the 
profile of $D(\lambda)$ as 
\begin{equation}\label{eq:3-70}
-x_r-1 < -x_{r-1}-1< \cdots < -x_1-1 < x_1 < \cdots < x_{r-1} < x_r, 
\qquad x_1>0 
\end{equation}
from its interlacing coordinates. 
Though $\mathrm{supp}\,\mathfrak{m}_{D(\lambda)}$ may contain $0$ 
besides \eqref{eq:3-70}, it has no effects on integrals like moments. 
We have 
\begin{align}
\frac{\chi^{(\lambda,\gamma)}}{\dim (\lambda,\gamma)}
(\widetilde{E}_n \widetilde{J}_{n+1}^{\;2k}) 
&= \sum_{\mu\in \mathcal{SP}_{n+1}: \mu\searrow\lambda} 
\Bigl( \frac{c(\mu/\lambda)\bigl( c(\mu/\lambda)+1\bigr)}{2}
\Bigr)^k \mathfrak{m}_{D(\lambda)} \bigl( \bigl\{ 
c(\mu/\lambda), -c(\mu/\lambda)-1\bigr\} \bigr) 
\notag \\ 
&= \frac{1}{2^k} \sum_{i=1}^r \bigl( x_i(x_i+1)\bigr)^k 
\bigl\{ \mathfrak{m}_{D(\lambda)}(\{x_i\}) + 
\mathfrak{m}_{D(\lambda)}(\{ -x_i-1\}) \bigr\}
\label{eq:3-71}
\end{align}
by Theorem~\ref{th:JMtr}. 
By using \eqref{eq:3-68} and \eqref{eq:3-71}, we restate 
\eqref{eq:3-30} and \eqref{eq:3-32} of Proposition~\ref{prop:3-5} 
as follows. 
Note that a normalized irreducible character is multiplicative on the center. 

\begin{corollary}\label{cor:3-8}
Assume $M^{(n)}$ satisfies approximate factorization property of 
Definition~\ref{def:AFP} 
and \eqref{eq:1-30}. 
We have for $\forall k\in\mathbb{N}$ 
\begin{align}
&\lim_{n\to\infty} 
\mathbb{E}^{(n)}\Bigl[ \frac{1}{n^k 2^k} \sum_{i=1}^r 
\bigl( x_i(x_i+1)\bigr)^k 
\mathfrak{m}_{D(\lambda)}(\{x_i, -x_i-1\}) \Bigr] =m_{2k}, 
\label{eq:3-72} \\ 
&\lim_{n\to\infty} 
\mathbb{E}^{(n)}\Bigl[ \Bigl\{ \frac{1}{n^k 2^k} \sum_{i=1}^r 
\bigl( x_i(x_i+1)\bigr)^k 
\mathfrak{m}_{D(\lambda)}(\{x_i, -x_i-1\}) -m_{2k}\Bigr\}^2 
\Bigr] =0 
\label{eq:3-73}
\end{align}
under \eqref{eq:3-70}. 
\hfill $\square$
\end{corollary}

Let $M_j(\mu)$ denote the $j$th moment of probability $\mu$ on $\mathbb{R}$. 
Recalling \eqref{eq:1-1}, we have 
\begin{equation}\label{eq:3-74}
M_{2k}(\mathfrak{m}_{D(\lambda)^{\sqrt{2}}}) = 
\frac{1}{2^k} M_{2k}(\mathfrak{m}_{D(\lambda)}) = 
\frac{1}{2^k}
\sum_{i=1}^r \bigl\{ x_i^{2k} \mathfrak{m}_{D(\lambda)}(\{x_i\}) + 
 (x_i+1)^{2k} \mathfrak{m}_{D(\lambda)}(\{-x_i-1\})\bigr\}.
\end{equation}
In order to get law of large numbers for the moment sequence of 
$\mathfrak{m}_{D(\lambda)}$ from Corollary~\ref{cor:3-8}, we need 
two modifications: 

\noindent 
-- the $2k$th-moment of $\mathfrak{m}_{D(\lambda)}$ is almost given by 
\eqref{eq:3-71} but not exactly, 

\noindent 
-- the odd moments of $\mathfrak{m}_{D(\lambda)}$ should be treated also 
because it is not exactly symmetric. 

\noindent 
Let us settle these small problems. 

\begin{lemma}\label{lem:3-9}
In the situation of Corollary~\ref{cor:3-8}, we have 
\begin{align*}
&\lim_{n\to\infty} \mathbb{E}^{(n)}\bigl[ 
M_{2k}(\mathfrak{m}_{D(\lambda)^{\sqrt{2n}}}) \bigr] = 
m_{2k}, \\ 
&\lim_{n\to\infty} \mathbb{E}^{(n)}\bigl[ \bigl( 
M_{2k}(\mathfrak{m}_{D(\lambda)^{\sqrt{2n}}}) - m_{2k}\bigr)^2 
\bigr] = 0 
\end{align*}
for $\forall k\in\mathbb{N}$. 
\hfill $\square$
\end{lemma}

\textit{Proof} \ 
Comparing \eqref{eq:3-72} and \eqref{eq:3-73} with \eqref{eq:3-74}, 
we have ony to verify $L^2$-convergence of the following random variables 
to $0$ as $n\to\infty$: 
\begin{align*}
&\frac{1}{n^k 2^k} \sum_{i=1}^r \bigl( x_i(x_i+1)\bigr)^k 
\mathfrak{m}_{D(\lambda)} (\{ x_i, -x_i-1\}) - 
M_{2k}(\mathfrak{m}_{D(\lambda)^{\sqrt{2n}}})  \\ 
&\quad = \frac{1}{n^k 2^k} \sum_{i=1}^r \bigr\{\bigl( x_i(x_i+1)\bigr)^k 
-x_i^{2k}\bigr\} \mathfrak{m}_{D(\lambda)} (\{ x_i\})  \\ 
&\qquad + 
\frac{1}{n^k 2^k} \sum_{i=1}^r \bigr\{\bigl( x_i(x_i+1)\bigr)^k 
-(x_i+1)^{2k}\bigr\} \mathfrak{m}_{D(\lambda)} (\{ -x_i-1\}). 
\end{align*}
It suffices to show $L^2$-convergence of 
\begin{equation}\label{eq:3-78}
\frac{1}{n^k} \int_{\mathbb{R}} |x|^p 
\mathfrak{m}_{D(\lambda)} (dx) 
\end{equation}
to $0$ as $n\to\infty$ for $\forall p\in\{0,1,\cdots, 2k-1\}$. 
H\"{o}lder's inequality gives 
\begin{equation}\label{eq:3-79}
n^{-\frac{p}{2}} \int |x|^p \mathfrak{m}_{D(\lambda)}(dx) \leqq 
\Bigl( n^{-k} \int x^{2k} \mathfrak{m}_{D(\lambda)}(dx) \Bigr) 
\vee 1. 
\end{equation}
The right hand side of \eqref{eq:3-79} is $L^2$-bounded by \eqref{eq:3-73}. 
This impies \eqref{eq:3-78} converges to $0$ in $L^2$. 
\hfill $\blacksquare$

\medskip

To treat odd moments, we mention an observation. 

\begin{lemma}\label{lem:3-10}
We have 
\begin{equation}\label{eq:3-80}
x_k \mathfrak{m}_{D(\lambda)}(\{x_k\}) = 
(x_k+1) \mathfrak{m}_{D(\lambda)}(\{-x_k-1\}), \qquad 
k\in \{1,\cdots, r\}
\end{equation}
under \eqref{eq:3-70}. 
\hfill $\square$
\end{lemma}

\textit{Proof} \ 
Recall the cases of (1) and (iib) in the proof of Lemma~\ref{lem:2-7}. 
We can read out \eqref{eq:3-80} from \eqref{eq:2-48} and \eqref{eq:2-49} 
for case (i), and from \eqref{eq:2-60} and \eqref{eq:2-61} for case (iib) 
respectively. 
\hfill $\blacksquare$

\begin{lemma}\label{lem:3-11}
In the situation of Corollary~\ref{cor:3-8}, we have 
\[ 
\lim_{n\to\infty} \mathbb{E}^{(n)}\Bigl[ 
M_{2k+1}(\mathfrak{m}_{D(\lambda)^{\sqrt{2n}}})^2 
\Bigr] =0
\]
%
for $\forall k\in\mathbb{N}$. 
\hfill $\square$
\end{lemma}

\textit{Proof} \ 
Applying 
Lemma~\ref{lem:3-10}, we have 
\begin{align}
M_{2k+1}(\mathfrak{m}_{D(\lambda)^{\sqrt{2n}}}) &= 
\frac{1}{(2n)^{k+\frac{1}{2}}} \sum_{i=1}^r \bigl\{ 
x_i^{2k+1} \mathfrak{m}_{D(\lambda)}(\{x_i\}) + 
(-x_i-1)^{2k+1} \mathfrak{m}_{D(\lambda)}(\{-x_i-1\})\bigr\} 
\notag \\ 
&= \frac{1}{(2n)^{k+\frac{1}{2}}} \sum_{i=1}^r x_i \bigl\{ 
x_i^{2k}- (x_i+1)^{2k}\bigr\} \mathfrak{m}_{D(\lambda)}(\{x_i\}). 
\label{eq:3-82}
\end{align}
Absolute value of $\sum\cdots$ in the right hand side is bounded by a 
non-negative linear combination of absolute moments up to $2k$th degree 
of $\mathfrak{m}_{D(\lambda)}$. 
Those divided by $n^k$ are all $L^2$-bounded by \eqref{eq:3-79}. 
Hence \eqref{eq:3-82} converges to $0$ in $L^2$. 
\hfill $\blacksquare$

\medskip 

\textit{Proof of Theorem~\ref{th:static}} \ 
Lemma~\ref{lem:3-9} and Lemma~\ref{lem:3-11} yield 
\begin{equation}\label{eq:3-83}
M_k(\mathfrak{m}_{D(\lambda)^{\sqrt{2n}}}) \ 
\xrightarrow[\,n\to\infty\,] \ m_k \qquad \text{in} \quad 
L^2\bigl( (\widetilde{\mathfrak{S}}_n)^\wedge_{\mathrm{spin}}, M^{(n)}
\bigr) 
\end{equation}
for $\forall k\in\mathbb{N}$. 
Since \eqref{eq:3-83} implies in particular 
\begin{equation}\label{eq:3-84}
m_k = \lim_{n\to\infty} \mathbb{E}^{(n)} \bigl[ 
M_k(\mathfrak{m}_{D(\lambda)^{\sqrt{2n}}}) \bigr], \qquad 
k\in\mathbb{N},  
\end{equation}
we see $(m_k)$ gives a moment sequence, which means the associated Hankel matrices 
are non-negative definite. 
Since we can take $\omega^{(n)}\in \mathbb{D}$ such that 
\begin{equation}\label{eq:3-85}
\mathfrak{m}_{\omega^{(n)}} = \sum_{(\lambda,\gamma)} 
M^{(n)}(\{(\lambda,\gamma)\}) \mathfrak{m}_{D(\lambda)^{\sqrt{2n}}},  
\end{equation}
we may write \eqref{eq:3-84} as 
\[ 
m_k = \lim_{n\to\infty} M_k (\mathfrak{m}_{\omega^{(n)}}), 
\qquad k\in\mathbb{N}. 
\] 

Case (1). Since \eqref{eq:1-31} assures 
\[ 
\exists C >0, \qquad |m_k|\leqq C^k, 
\] 
there exists a unique probability $\mu$ on $\mathbb{R}$ supported by a compact 
subset which has $m_k$ as its $k$th moment: 
\begin{equation}\label{eq:3-88}
m_k = M_k(\mu), \qquad k\in\mathbb{N}.
\end{equation}
Then, the uniformly convergent limit is captured in $\mathbb{D}$ as 
\begin{equation}\label{eq:3-89}
\omega = \lim_{n\to\infty} \omega^{(n)} \in \mathbb{D}, \qquad 
\mathfrak{m}_\omega =\mu.
\end{equation}
With help of Chebyshev's inequality, \eqref{eq:3-83} yields 
\begin{equation}\label{eq:3-90}
\lim_{n\to\infty} M^{(n)}\Bigl( \Bigl\{ (\lambda, \gamma)\in 
(\widetilde{\mathfrak{S}}_n)^\wedge_{\mathrm{spin}} \,\Big|\, 
| M_k(\mathfrak{m}_{D(\lambda)^{\sqrt{2n}}}) - 
M_k(\mathfrak{m}_\omega) | > \epsilon \Bigr\}\Bigr) =0
\end{equation}
for $\forall\epsilon >0$ and $\forall k\in\mathbb{N}$. 
The topology of $\mathbb{D}$ assures that 
$\forall\epsilon >0$, \ $\exists\delta >0, \ k_1, \cdots, k_p$ such that 
\[ 
\max_{i\in\{1,\cdots,p\}} 
| M_{k_i}(\mathfrak{m}_{D(\lambda)^{\sqrt{2n}}}) - 
M_{k_i}(\mathfrak{m}_\omega) | <\delta \ \Longrightarrow \ 
\sup_{x\in\mathbb{R}} \bigl| D(\lambda)^{\sqrt{2n}}(x) - 
\omega (x)\bigr| < \epsilon.
\] 
Hence we get weak law of large numbers 
\begin{equation}\label{eq:3-92}
\lim_{n\to\infty} M^{(n)}\Bigl( \Bigl\{ (\lambda, \gamma)\in 
(\widetilde{\mathfrak{S}}_n)^\wedge_{\mathrm{spin}} \,\Big|\, 
\sup_{x\in\mathbb{R}} \bigl| D(\lambda)^{\sqrt{2n}}(x) - 
\omega (x)\bigr| > \epsilon \Bigr\}\Bigr) =0.
\end{equation}

Case (2). 
Since the moment problem is determinate, $\{ \mathfrak{m}_{\omega^{(n)}} \}$ 
of \eqref{eq:3-85} weakly converges to the uniqulely determined probability $\mu$ 
of \eqref{eq:3-88}. 
Then, since corresponding limit is captured in $\mathcal{D}$, we have 
\eqref{eq:3-89} such that $\omega\in\mathcal{D}$. 
Through the homemomorphism $\mathcal{M} \cong \mathcal{D}$
\footnote{See \cite{Ker98} and \S\S\ref{subsec:A2}.}, 
\eqref{eq:3-90} yields 
for $\forall \epsilon >0$ and $\forall r>0$ 
\begin{equation}\label{eq:3-93}
\lim_{n\to\infty} M^{(n)}\Bigl( \Bigl\{ (\lambda, \gamma)\in 
(\widetilde{\mathfrak{S}}_n)^\wedge_{\mathrm{spin}} \,\Big|\, 
\sup_{x\in [-r, r]} \bigl| D(\lambda)^{\sqrt{2n}}(x) - 
\omega (x)\bigr| \leqq \epsilon \Bigr\}\Bigr) =1.
\end{equation}
Since $\omega\in\mathcal{D}$ satisfies 
\[ 
\lim_{x\pm \infty} (\omega(x)-|x|) =0, \qquad 
\omega(x) \geqq |x|, 
\] 
the range of $\sup$ in \eqref{eq:3-93} can be modified to the whole $\mathbb{R}$. 
We thus get \eqref{eq:3-92} in Case (2), too. 
\hfill $\blacksquare$

\section{Proof for dynamic model}\label{sect:4}

Let us proceed to prove Theorem~\ref{th:dynamic}. 
In order to see initial approximate factorization proeprty is propagated at 
macroscopic time $t>0$, it is important to observe the relation between 
the transition matrix of Res-Ind chain and the irreducible character values. 

In general, for a finite group $G$ and its subgroup $H$, let us consider 
Res-Ind chain on $\widehat{G}$ whose transition matrix is given by 
$P$ of \eqref{eq:1-39}. 
The irreducible character value at an element of conjugacy class $C$ of $G$ 
is denoted by $\chi^\xi_C$ for $\xi\in \widehat{G}$. 
Take a column vector 
\begin{equation}\label{eq:4-1}
\chi_C = \bigl[ \frac{1}{\dim\xi} \chi^\xi_C \bigr]_{\xi\in\widehat{G}}. 
\end{equation}

\begin{lemma}\label{lem:4-1}
We have 
\[
P \chi_C = \frac{|C\cap H|}{|C|} \chi_C.
\] 
Namely, $\chi_C$ is an eigenvector of $P$. 
\hfill $\square$
\end{lemma}

\textit{Proof} is given by the  induced character formula. 
See \cite{Hor15}, \cite[Lemma~5.2]{Hor16}. 
\hfill $\blacksquare$

\begin{lemma}\label{lem:4-2}
For transition matrix $P^{(n)} = P_{\mathrm{spin}}$ of \eqref{eq:1-46} on 
$(\widetilde{\mathfrak{S}}_n)^\wedge_{\mathrm{spin}}$ and $\rho\in\overline{\mathcal{OP}}$, 
we have
\footnote{Here $(\lambda,\gamma)$ is a Nazarov parameter of \eqref{eq:1-16}, 
$[\rho,1^{n-|\rho|}]$ is an element of $\widetilde{\mathfrak{S}}_n$ 
defined by \eqref{eq:3-2}.}
\begin{equation}\label{eq:4-3}
P^{(n)} \Bigl[ \frac{\chi^{(\lambda,\gamma)}([\rho, 1^{n-|\rho|}])}
{\dim(\lambda,\gamma)} \Bigr]_{(\lambda,\gamma)\in 
(\widetilde{\mathfrak{S}}_n)^\wedge_{\mathrm{spin}}} = 
\bigl( 1-\frac{|\rho|}{n}\bigr) 
\Bigl[ \frac{\chi^{(\lambda,\gamma)}([\rho, 1^{n-|\rho|}])}
{\dim(\lambda,\gamma)} \Bigr]_{(\lambda,\gamma)\in 
(\widetilde{\mathfrak{S}}_n)^\wedge_{\mathrm{spin}}}, 
\end{equation}
which is valid also for $z [\rho, 1^{n-|\rho|}]$ instead of 
$[\rho, 1^{n-|\rho|}]$. 
\hfill $\square$ 
\end{lemma}

\textit{Proof} \ 
Apply Lemma~\ref{lem:4-1} for $G= \widetilde{\mathfrak{S}}_n$ and 
$H= \widetilde{\mathfrak{S}}_{n-1}$. 
Since $P$ is block diagonal as \eqref{eq:1-46}, the vector consisting of the 
$(\widetilde{\mathfrak{S}}_n)^\wedge_{\mathrm{spin}}$-entries of \eqref{eq:4-1} 
is an eigenvector of $P^{(n)}$. 
Here we take $\widetilde{C}_{(\rho, 1^{n-|\rho|})}$ of \eqref{eq:3-3} 
as conjugacy class $C$ and $[\rho, 1^{n-|\rho|}]$ of \eqref{eq:3-2} as 
a representative. 
Then, since $\rho\in\overline{\mathcal{OP}}$, we have 
\begin{align*}
&|\widetilde{C}_{(\rho, 1^{n-|\rho|})}| = |C_{(\rho, 1^{n-|\rho|})}| = 
\frac{1}{z_\rho} n^{\downarrow |\rho|}, \\ 
&|\widetilde{C}_{(\rho, 1^{n-|\rho|})}\cap \widetilde{\mathfrak{S}}_{n-1}| 
= | \widetilde{C}_{(\rho, 1^{n-1-|\rho|})}| = |C_{(\rho, 1^{n-1-|\rho|})}| 
= \frac{1}{z_\rho} (n-1)^{\downarrow |\rho|} 
\end{align*}
and, by taking the ratio, 
\[
\frac{(n-1)^{\downarrow |\rho|}}{n^{\downarrow |\rho|}} = 
\frac{n-1-|\rho|+1}{n} = 1- \frac{|\rho|}{n}. 
\] 
This implies \eqref{eq:4-3}. 
\hfill $\blacksquare$ 

\medskip

Let $f_s^{(n)}$ be the spin normalized central positive-definite function 
corresponding to distribution $M_s^{(n)}$ at time $s$ given by \eqref{eq:1-51}. 

\begin{lemma}\label{lem:4-3}
Let $x\in \widetilde{\mathfrak{S}}_\infty \setminus \{e\}$ have type 
$\rho\in\overline{\mathcal{OP}}$. 
In other words, let $\mathrm{type}\,x = (\rho, 1^{n-|\rho|})$ hold 
under \eqref{eq:1-28} for sufficiently large $n\in\mathbb{N}$. 
Then 
\begin{equation}\label{eq:4-7}
f_s^{(n)}(x) = f_0^{(n)}(x) \sum_{j=0}^\infty \bigl( 1-\frac{|\rho|}{n} 
\bigr)^j \int_{[0, s]} \psi \bigl( (s-u, \infty)\bigr) \psi^{\ast j}(du) 
\end{equation}
holds. 
\hfill $\square$
\end{lemma}

\textit{Proof} \ 
We get from \eqref{eq:1-25} and \eqref{eq:1-51} 
\begin{align}
f_s^{(n)}(x) &= 
\sum_{(\lambda,\gamma)\in (\widetilde{\mathfrak{S}}_n)^\wedge_{\mathrm{spin}}} 
M_s^{(n)}(\{(\lambda,\gamma)\}) 
\frac{\chi^{(\lambda,\gamma)}(x)}{\dim(\lambda,\gamma)} 
\notag \\ 
&= \sum_{j=0}^\infty 
\Bigl( \int_{[0, s]} \psi \bigl( (s-u, \infty)\bigr) \psi^{\ast j}(du)\Bigr) \! 
\sum_{(\lambda,\gamma)\in (\widetilde{\mathfrak{S}}_n)^\wedge_{\mathrm{spin}}} 
\!\! (M_0^{(n)} P^{(n) j})_{(\lambda,\gamma)} 
\frac{\chi^{(\lambda,\gamma)}(x)}{\dim(\lambda,\gamma)}. 
\label{eq:4-8} 
\end{align}
Since Lemma~\ref{lem:4-2} and \eqref{eq:1-25} yield 
\begin{align*}
M_0^{(n)} P^{(n) j} \bigl[ 
\frac{\chi^{(\lambda,\gamma)}(x)}{\dim(\lambda,\gamma)} 
\bigr]_{(\lambda,\gamma)\in (\widetilde{\mathfrak{S}}_n)^\wedge_{\mathrm{spin}}} 
&= \bigl( 1- \frac{|\rho|}{n}\bigr)^j M_0^{(n)} \bigl[ 
\frac{\chi^{(\lambda,\gamma)}(x)}{\dim(\lambda,\gamma)} 
\bigr]_{(\lambda,\gamma)\in (\widetilde{\mathfrak{S}}_n)^\wedge_{\mathrm{spin}}} 
\\ 
&= \bigl( 1- \frac{|\rho|}{n}\bigr)^j f_0^{(n)}(x),  
\end{align*}
the right hand side of \eqref{eq:4-7} is equal to \eqref{eq:4-8}. 
\hfill $\blacksquare$ 

\begin{lemma}\label{lem:4-4} 
Assume \eqref{eq:1-53} for the pausing time distribution $\psi$ of $(X_s^{(n)})$ 
where the mean of $\psi$ is $m$. 
Then 
\begin{equation}\label{eq:4-10}
\sum_{j=0}^\infty \bigl( 1-\frac{k}{n} 
\bigr)^j \int_{[0, tn]} \psi \bigl( (tn-u, \infty)\bigr) \psi^{\ast j}(du) \ 
\xrightarrow[\,n\to\infty\,] \ e^{-\frac{kt}{m}}
\end{equation}
holds for $t>0, \ k\in\{2,3,\cdots\}$. 
\hfill $\square$  
\end{lemma}

\textit{Proof} is the same as that of Proposition~2.1 of \cite{Hor20}
\footnote{The proof is carried out through some estimates for characteristic 
function $\varphi$. 
Usefulness of such Fourier-analytic methods for treating general pausing time 
distributions is illustrated in \cite{Wei94}.}. 
\hfill $\blacksquare$ 

\medskip

We are ready to prove Theorem\ref{th:dynamic}. 

\medskip

\textit{Proof of Theorem\ref{th:dynamic}}

\noindent\textit{Step~1} \ 
We show the sequence $\{f_{tn}^{(n)}\}_{n\in\mathbb{N}}$ of 
spin normalized central positive-definite functions satisfies 
approximate factorization property for $t>0$. 
Let $x, y\in \widetilde{\mathfrak{S}}_\infty \setminus\{e\}$ satisfy 
$\mathrm{supp}\,x\cap \mathrm{supp}\,y =\varnothing$ and have 
types as \eqref{eq:1-27} with $\rho, \sigma\in \overline{\mathcal{OP}}$. 
For simplicity, $a(k,t,n)$ denotes the left hand side of \eqref{eq:4-10}. 
Lemma~\ref{lem:4-3} yields 
\begin{align}
f_{tn}^{(n)}(xy) &= a( |\rho|+|\sigma|, t, n) f_0^{(n)}(xy), \notag \\ 
f_{tn}^{(n)}(x) f_{tn}^{(n)}(y) &= a(|\rho|, t, n) a(|\sigma|, t, n) 
f_0^{(n)}(x) f_0^{(n)}(y). \label{eq:4-11}
\end{align}
Since $\{f_0^{(n)}\}_{n\in\mathbb{N}}$ satisfies approximate factorization 
property, we have 
\begin{equation}\label{eq:4-12}
f_0^{(n)}(x) = O\bigl( n^{-\frac{|\rho|-l(\rho)}{2}} \bigr), \qquad 
f_0^{(n)}(y) = O\bigl( n^{-\frac{|\sigma|-l(\sigma)}{2}} \bigr)
\end{equation}
similarly to Lemma~\ref{lem:3-2}. 
Difference of the two of \eqref{eq:4-11} is written as 
\begin{multline}\label{eq:4-13}
a(|\rho|+|\sigma|, t, n) \bigl\{ f_0^{(n)}(xy) - f_0^{(n)}(x) f_0^{(n)}(y) 
\bigr\} \\ 
+ \bigl\{ a(|\rho|+|\sigma|, t, n) - a(|\rho|, t, n) a(|\sigma|, t, n) \bigr\} 
f_0^{(n)}(x) f_0^{(n)}(y). 
\end{multline}
Approximate factorization property with Lemma~\ref{lem:4-4} and \eqref{eq:4-12} 
yields 
\[ 
\text{\eqref{eq:4-13}} = 
o\bigl( n^{-\frac{|\rho|-l(\rho)+|\sigma|-l(\sigma)}{2}} \bigr) + 
o(1) O \bigl( n^{-\frac{|\rho|-l(\rho)+|\sigma|-l(\sigma)}{2}} \bigr). 
\] 
This implies $\{f_{tn}^{(n)}\}$ satisfies approximate factorization property, too. 

\noindent\textit{Step~2} \ 
We see 
\begin{equation}\label{eq:4-15}
\lim_{n\to\infty} n^{\frac{k-1}{2}} f_{tn}^{(n)} \bigl( 
[1\ 2\ \cdots\ k]\bigr) = r_{k+1} e^{-\frac{kt}{m}} 
\end{equation}
holds for $k\in\mathbb{N}$, odd, $\geqq 3$. 
In fact, \eqref{eq:1-52}, Lemma~\ref{lem:4-3} and Lemma~\ref{lem:4-4} yield 
%
\[ 
n^{\frac{k-1}{2}} f_{tn}^{(n)} \bigl( [1\ 2\ \cdots\ k]\bigr) = 
n^{\frac{k-1}{2}} f_0^{(n)} \bigl( [1\ 2\ \cdots\ k]\bigr) a(k,t,n) 
\ \xrightarrow[\,n\to\infty\,] \ r_{k+1} e^{-\frac{kt}{m}}. 
\]
%
\textit{Step~3} \ 
Assume the situation of (1), that is, \eqref{eq:1-54}. 
If $r_{k+1}$ is replaced by $r_{k+1} e^{-\frac{kt}{m}}$ of \eqref{eq:4-15}, 
the condition like \eqref{eq:1-31} still holds. 
Form sequence $(r_{t,j})_{j\in\mathbb{N}}$ from 
$r_{k+1} e^{-\frac{kt}{m}}$ of \eqref{eq:4-15} by inserting $1$ as 
the second and $0$ as those of odd degrees. 
Then the free cumulant-moment formula gives us sequence $(m_{t,j})$: 
\begin{equation}\label{eq:4-161}
m_{t,j} = \sum_{\rho\in \mathcal{P}_j} |NC(\rho)| r_{t, \rho_1}\cdots 
r_{t, \rho_{l(\rho)}}.
\end{equation}
We apply (1) of Theorem~\ref{th:static} to $\{f_{tn}^{(n)}\}$. 
Let $\mu_t$ denote the obtained probability on $\mathbb{R}$ with 
compact support. 
The moment sequence of $\mu_t$ is $(m_{t,j})$: 
\begin{equation}\label{eq:4-162}
M_j(\mu_t) = m_{t,j}. 
\end{equation}
The transition measure of the limit shape of $D(\lambda)^{\sqrt{2n}}$ 
under $\{M_{tn}^{(n)}\}$ is $\mu_t$. 
Its free cumulants are 
\begin{align}
&R_1(\mu_t) = R_3(\mu_t) = \cdots =0, \qquad R_2(\mu_t) =1, \notag \\ 
&R_j(\mu_t) = r_j e^{-\frac{(j-1)t}{m}} \qquad (j\geqq 4, \ \text{even}) 
\label{eq:4-19} 
\end{align}
from \eqref{eq:4-161} and \eqref{eq:4-162}. 
Note \eqref{eq:4-19} holds also for $t=0$. 
Take $\omega_t\in\mathbb{D}$ corresponding to $\mu_t$ by 
the Markov transform: 
\begin{equation}\label{eq:4-17}
\mathfrak{m}_{\omega_t} = \mu_t. 
\end{equation}
We get from \eqref{eq:1-32} 
\[ 
\lim_{n\to\infty} M_{tn}^{(n)}\Bigl( \Bigl\{ (\lambda,\gamma)\in 
(\widetilde{\mathfrak{S}}_n)^\wedge_{\mathrm{spin}} \,\Big|\, \sup_{x\in\mathbb{R}} 
\bigl| D(\lambda)^{\sqrt{2n}}(x)- \omega_t(x) \bigr| > \epsilon 
\Bigr\} \Bigr) =0.
\] 
By \eqref{eq:4-19} and \eqref{eq:4-17}, the free cumulant sequence of 
$\mathfrak{m}_{\omega_t}$ satisfies \eqref{eq:1-57}. 
The right hand side of \eqref{eq:1-56} is determined from free compression 
and free convolution for probabilities on $\mathbb{R}$ with compact supports, 
whose free cumulant sequence is given by \eqref{eq:1-57}. 

\noindent\textit{Step~4} 
Assume the situation of (2), that is, \eqref{eq:1-571}. 
For $(m_{t,j})$ determined by \eqref{eq:4-161}, noting 
\[ 
|NC(j)| \leqq \frac{2^{2j}}{j}, \qquad 
|\mathcal{P}_j| \leqq (\mathrm{const.})^{\sqrt{j}}
\] 
hold, we have 
\begin{align}
m_{t,2k} &= \sum_{\rho\in \mathcal{P}_{2k} \text{ : even length rows}} |NC(\rho)| 
r_4^{m_4(\rho)}r_6^{m_6(\rho)}\cdots e^{-\frac{3m_4(\rho)t}{m}} 
e^{-\frac{5m_6(\rho)t}{m}}\cdots \notag \\ 
&\leqq \sum_{\rho\in \mathcal{P}_{2k} \text{ : even length rows}} 
|NC(\rho)| C^{4m_4(\rho)} 4^{4m_4(\rho)} C^{6m_6(\rho)}
6^{m_6(\rho)}\cdots \notag \\ 
&\leqq \sum_{\rho\in \mathcal{P}_{2k} \text{ : even length rows}} 
|NC(\rho)| (C\cdot 2k)^{4m_4(\rho)+6m_6(\rho)+\cdots} 
\leqq (\mathrm{const.})^{2k} (2k)^{2k}. 
\label{eq:4-25}
\end{align}
Since \eqref{eq:4-25} implies Carleman's condition 
\[ 
\sum_k m_{t, 2k}^{-\frac{1}{2k}} = \infty, 
\] 
the moment problem for $(m_{t,j})$ is determinate. 
This enables us to apply (2) of Theorem~\ref{th:static}. 
Then, the argument proceeds similarly to Step~3 to reach \eqref{eq:1-57}. 
Note, however, the limit shape $\omega_t$ belongs to $\mathcal{D}$, 
not necessarily to $\mathbb{D}$. 
\hfill $\blacksquare$

\section{Examples}\label{sect:5}

\subsection{Vershik curve}

The limit shape of $Y(\lambda)$ for $\lambda\in \mathcal{P}_n$ with respect to the uniform 
probability $M_{\mathrm{U}}^{(n)}$ on $\mathcal{P}_n$ is given by the Vershik curve 
$\psi_{\mathrm{U}}$ of \eqref{eq:1-5}. 
Though the one with respect to the uniform probability $M_{\mathrm{US}}^{(n)}$ on 
$\mathcal{SP}_n$ is given by $\psi_{\mathrm{US}}$ of \eqref{eq:1-5}, 
it coincides with the right half of $\psi_{\mathrm{U}}$ if it is 
regarded as the limit shape of $S(\lambda)$ tranformed by \eqref{eq:1-6}. 
Let us reformulate the fact as the limit shape of $D(\lambda)$ in Figure~\ref{fig:1-3} as 
\begin{align}
&M_{\mathrm{US}}^{(n)} \Bigl( \Bigl\{ \lambda\in \mathcal{SP}_n \,\Big|\, \sup_{x\in\mathbb{R}}
\bigl| D(\lambda)^{\sqrt{2n}}(x) - \varOmega_{\mathrm{V}}(x)\bigr| >\epsilon\Bigr\}\Bigr) 
\ \xrightarrow[\, n\to\infty\,] \ 0, 
\label{eq:5-1} \\ 
&\varOmega_{\mathrm{V}}(x) = \frac{2\sqrt{6}}{\pi} \log 
\bigl( e^{\frac{\pi}{2\sqrt{6}}x} + e^{-\frac{\pi}{2\sqrt{6}}x}\bigr). 
\label{eq:5-2}
\end{align}
\begin{figure}[hpt] 
\centering 
\includegraphics[width=5cm]{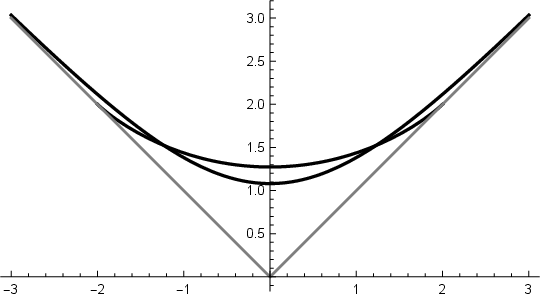}
\caption{Vershik curve \eqref{eq:5-2} and  VKLS curve \eqref{eq:1-10}}
\end{figure}
Since $\varOmega_{\mathrm{V}}$ belongs to $\mathcal{D}$ but not to $\mathbb{D}$ 
(that is, its transition measure has noncompact support), we should be careful. 
Rayleigh measure $\tau_{\mathrm{V}}$ corresponding to \eqref{eq:5-2} is the probability 
on $\mathbb{R}$ whose density is given by 
\begin{equation}\label{eq:5-3}
\frac{d\tau_{\mathrm{V}}}{dx} = \frac{\pi}{\sqrt{6}} 
\bigl( e^{\frac{\pi}{2\sqrt{6}}x} + e^{-\frac{\pi}{2\sqrt{6}}x}\bigr)^{-2}. 
\end{equation}
The moments of $\tau_{\mathrm{V}}$ are expressed in terms of the Bernoulli numbers $B_n$ 
(e.g. see \cite{GrRy07}). 
Determining $B_n$ by 
\begin{equation}\label{eq:5-4}
\frac{t}{e^t-1} = \sum_{n=0}^\infty B_n \frac{t^n}{n!}, 
\end{equation}
we have 
\[
|B_{2k}| = \frac{1}{\pi^{2k}} \int_0^\infty \frac{x^{2k}}{\sinh^2 x}\, dx 
= \frac{2^{2k-1}}{(2^{2k-1}-1)\pi^{2k}} \int_0^\infty \frac{x^{2k}}{\cosh^2 x}\, dx 
\]
as integral expressions. 
Hence, for $k\in\mathbb{N}$, \eqref{eq:5-3} satisfies  
\begin{equation}\label{eq:5-5}
M_{2k}(\tau_{\mathrm{V}}) = \int_{\mathbb{R}} \frac{\pi}{\sqrt{6}} 
\frac{x^{2k}}{( e^{\frac{\pi}{2\sqrt{6}}x} + e^{-\frac{\pi}{2\sqrt{6}}x})^2}\, dx 
= (2^{2k}-2) 6^k |B_{2k}| 
\end{equation}
and furthermore 
\begin{equation}\label{eq:5-6}
\bigl( 1- \frac{1}{2^{2k-1}}\bigr) \frac{2\cdot6^k}{\pi^{2k}}(2k)! 
\leqq M_{2k} (\tau_{\mathrm{V}}) \leqq \frac{2\cdot 6^k}{\pi^{2k}}(2k)!. 
\end{equation}
In particular, the moment problem for $\tau_{\mathrm{V}}$ is determinate. 
Let us compute the free cumulants of transition measure $\mathfrak{m}_{\mathrm{V}}$ 
corresponding to \eqref{eq:5-2} and verify the estimate of \eqref{eq:1-571} 
in Theorem~\ref{th:dynamic} (2) is valid. 

In general, Rayleigh measure $\tau$, transition measure $\mu$, $k$th moment 
$M_k(\tau)$, $k$th free cumulant $R_k(\mu)$ and generating function 
(Stieltjes transform) $G_\mu(z)$ are connected by 
\begin{align}
G_\mu(z) &= \frac{1}{z} \exp \sum_{k=1}^\infty \frac{M_k(\tau)}{k} 
\frac{1}{z^k}, \label{eq:5-7} \\ 
R_k(\mu) &= -\frac{1}{k-1} [z^{-1}] \bigl( \frac{1}{G_\mu(z)^{k-1}}\bigr) 
\qquad (k\geqq 2). \label{eq:5-8}
\end{align}
Here the supports of $\tau$ and $\mu$ are not necessarily compact. 
We recognize that \eqref{eq:5-7} and \eqref{eq:5-8} live in $\mathbb{C}((z^{-1}))$, 
the field of formal power series of $z^{-1}$. 
Putting $M_{2j-1}(\tau) =0$ ($j\in\mathbb{N}$) in \eqref{eq:5-7}, we get from 
\eqref{eq:5-8} 
\begin{equation}\label{eq:5-9}
R_{2k}(\mu) = \sum_{l=1}^k \frac{(1-2k)^{l-1}}{l!} 
\sum_{j_i\in\mathbb{N}: \, j_1+\cdots+j_l =k} 
\frac{M_{2j_1}(\tau)\cdots M_{2j_l}(\tau)}{2^l j_1\cdots j_l} 
\end{equation}
for $k\in\mathbb{N}$. 
Applying \eqref{eq:5-9} to $\tau =\tau_{\mathrm{V}}$ and 
$\mu=\mathfrak{m}_{\mathrm{V}}$ and using \eqref{eq:5-6}, we have 
\begin{align*}
|R_{2k}(\mathfrak{m}_{\mathrm{V}})| &\leqq \sum_{l=1}^k \frac{(2k-1)^{l-1}}{l!} 
\bigl( \frac{\sqrt{6}}{\pi}\bigr)^{2k} \sum_{j_i\in\mathbb{N}:\, j_1+\cdots+j_l =k} 
\frac{(2j_1)! \cdots (2j_l)!}{j_1\cdots j_l} \\ 
&\leqq \bigl( \frac{\sqrt{6}}{\pi}\bigr)^{2k} 2^{2k} \sum_{l=1}^k 
\frac{(2k-1)^{l-1}}{l!} (k-l+1)^{2k-l} \binom{k+l-1}{l-1}, 
\end{align*}
by Stirling's formula 
\begin{align}
&\leqq \bigl( \frac{\sqrt{6}}{\pi}\bigr)^{2k} 2^{2k} k 
\max_{1\leqq l\leqq k} \frac{(2k-1)^{l-1} (k-l+1)^{2k-l} (k+l-1)^{l-1}}
{(l+1)^{l+\frac{1}{2}} l^{l-\frac{1}{2}}} 
\frac{e^{2l+1}}{2\pi} 
\notag \\ 
&\leqq \bigl( \frac{\sqrt{6}}{\pi}\bigr)^{2k} 2^{2k} k 
\frac{e^{2k+1}}{2\pi} \max_{1\leqq l\leqq k} 
\frac{2^k k^l k^{2k-l} 2^k k^l}{l^{2l}} 
= \frac{k}{2\pi} 2^{4k} e^{2k+1} \bigl( \frac{\sqrt{6}}{\pi}\bigr)^{2k} 
k^{2k} \max_{1\leqq l\leqq k}  \frac{k^l}{l^{2l}}.
\label{eq:5-11}
\end{align}
Since $k\geqq 8$ yields 
\[ 
\max_{1\leqq x\leqq k} \frac{k^x}{x^{2x}} = e^{2\sqrt{k}/e}, 
\] 
\eqref{eq:5-11} tells that 
\begin{equation}\label{eq:5-16}
R_{2k}(\mathfrak{m}_{\mathrm{V}}) 
= (-1)^{k+1} 6^k \sum_{l=1}^k \frac{(2k-1)^{l-1}}{l!} \!\!
\sum_{j_i\in\mathbb{N}: j_1+\cdots+j_l =k} \!\!\!\!
\frac{(2^{2j_1-1}-1)B_{2j_1}}{j_1} \cdots \frac{(2^{2j_l-1}-1)B_{2j_l}}{j_l}
\end{equation}
obtained by \eqref{eq:5-5} and \eqref{eq:5-9} satisfies the estimate of \eqref{eq:1-571}
\footnote{Note $|B_{2n}| = (-1)^{n-1} B_{2n}$ ($n\in\mathbb{N}$).}. 
Then, moment $M_{2k}(\mathfrak{m}_{\mathrm{V}})$ satisfies the same kind of estimate 
by \eqref{eq:4-25}. 
In particular, the corresponding moment problem is determinate. 

Let us ask whether one can take the Vershik curve $\varOmega_{\mathrm{V}}$ 
as an initial condition $\omega_0$ in Theorem~\ref{th:dynamic}. 
That will be realized by 
considering delta measures as initial distributions (though less interesting). 
Approximate $\varOmega_{\mathrm{V}}$ by a sequence of rectangular diagrams 
$\{ \omega^{(n)}\}_n$ in $\mathcal{D}$. 
Slightly modifying $\omega^{(n)}$ if necessary, we may set 
\begin{equation}\label{eq:5-14}
\omega^{(n)} = D(\lambda^{(n)})^{\sqrt{2n}}, \qquad 
\lambda^{(n)}\in \mathcal{SP}_n. 
\end{equation}
Then, $\{\mathfrak{m}_{\omega^{(n)}}\}$ converges weakly to 
$\mathfrak{m}_{\mathrm{V}}$ as $n\to\infty$ (see \S\S\ref{subsec:A2}). 
This yields convergence of the moments and hence of the free cumulants: 
\ for $\forall k\in\mathbb{N}$, 
\begin{equation}\label{eq:5-15}
\lim_{n\to\infty} M_{2k}(\mathfrak{m}_{\omega^{(n)}}) = 
M_{2k}(\mathfrak{m}_{\mathrm{V}}), \qquad 
\lim_{n\to\infty} R_{2k}(\mathfrak{m}_{\omega^{(n)}}) = 
R_{2k}(\mathfrak{m}_{\mathrm{V}}).
\end{equation}
It is not obvious that convergence of the spin irreducible character values of 
$\widetilde{\mathfrak{S}}_n$ follows from \eqref{eq:5-15} and \eqref{eq:5-14}. 
However, if we take for granted the spin Kerov polynomials connecting the 
spin irreducible character values of $\widetilde{\mathfrak{S}}_n$ and 
free cumulants (\cite{Mat18}), the convegence \eqref{eq:5-15} of 
free cumulants yields 
\[ 
n^{\frac{2k-1-1}{2}} \frac{\chi^{(\lambda^{(n)},\gamma)}}
{\dim (\lambda^{(n)},\gamma)} \bigl( [ 1\ 2\ \cdots \ 2k\!-\!1]\bigr) 
\ \xrightarrow[\,n\to\infty\,] \ R_{2k}(\mathfrak{m}_{\mathrm{V}}), \qquad 
k\in\mathbb{N}. 
\] 
Then, by using multiplicativity of a normalized irreducible character on the center 
and an elementary argument on conjugacy classes (as in \eqref{eq:3-14} and 
\eqref{eq:3-15}), we see that 
$\{\chi^{(\lambda^{(n)},\gamma)}/ \dim (\lambda^{(n)},\gamma)\}$ 
satisfies \eqref{eq:3-8} and approximate factorization property. 
Hence \textit{Theorem~\ref{th:dynamic} is applicable with Vershik curve 
$\varOmega_{\mathrm{V}}$ as an initial shape}. 
In particular, $\omega_t$ interpolates between the Vershik curve ($t=0$) and 
the VKLS curve ($t=\infty$). 
While we know free cumulants $R_k(\mathfrak{m}_{\omega_t})$ for $\omega_t$, 
an explicit form of curve $\omega_t$ itself is unclear. 

Since the Vershik curve is the limit shape of Young diagrams in a uniform ensemble, 
it may be desirable to take uniform distributions as initial distributions under which 
the initial shape is given by the Vershik curve. 
There are two parameters $(\lambda, 1)$ and $(\lambda, -1)$ of spin 
irreducible representations for $\lambda\in \mathcal{SP}_n^-$. 
Hence, when the uniform distribution on $\mathcal{SP}_n$ is interpreted as a probability 
on $(\widetilde{\mathfrak{S}}_n)^\wedge_{\mathrm{spin}}$, we regard it as 
\begin{equation}\label{eq:5-18}
M^{(n)}(\{(\lambda, \gamma)\}) = \begin{cases} 1/|\mathcal{SP}_n|, & 
\lambda\in \mathcal{SP}_n^+ \\ 1/ 2|\mathcal{SP}_n|, & \lambda\in \mathcal{SP}_n^-. \end{cases} 
\end{equation}
Let $f^{(n)}$ be the function corresponding to $M^{(n)}$ of \eqref{eq:5-18} 
through \eqref{eq:1-25}. 
We have 
\begin{align}
&f^{(n)}([1\ 2\ \cdots \ k]) = 
\sum_{(\lambda,\gamma)\in (\widetilde{\mathfrak{S}}_n)^\wedge_{\mathrm{spin}}} 
M^{(n)}(\{(\lambda,\gamma)\}) 
\frac{\chi^{(\lambda,\gamma)}([1\ 2\ \cdots\ k])}{\dim(\lambda,\gamma)} 
= \sum_{(\lambda,\gamma): \, \lambda\in \mathcal{SP}_n^+} + 
\sum_{(\lambda,\gamma): \, \lambda\in \mathcal{SP}_n^-} 
\notag \\ 
&= \sum_{\lambda\in \mathcal{SP}_n^+} \frac{1}{|\mathcal{SP}_n|} 
\frac{\chi^{(\lambda,\gamma)}}{\dim(\lambda,\gamma)}([1\ 2\ \cdots\ k])
+ \sum_{\lambda\in \mathcal{SP}_n^-} \frac{1}{2 |\mathcal{SP}_n|} 
\Bigl( \frac{\chi^{(\lambda, 1)}}{\dim(\lambda, 1)} + 
\frac{\chi^{(\lambda, -1)}}{\dim(\lambda, -1)}\Bigr) 
([1\ 2\ \cdots\ k])
\notag \\ 
&= \sum_{\lambda\in \mathcal{SP}_n} \frac{1}{|\mathcal{SP}_n|} 
\frac{\chi^{(\lambda,\gamma)}}{\dim(\lambda,\gamma)}([1\ 2\ \cdots\ k])
\label{eq:5-19}
\end{align}
for $k\in\mathbb{N}$, odd, $\geqq 3$. 
Here in \eqref{eq:5-19}, we take just one of $\gamma =\pm 1$ if 
$\lambda\in \mathcal{SP}_n^-$. 
Since we know the Vershik curves appear as the LLN limits as shown in \cite{Ver96}, 
we expect $f^{(n)}$ to be adapted in the framework of Theorem~\ref{th:static}. 
Unfortunately, however, we do not have a rigorous proof. 

\medskip

\noindent\textbf{Problem~1a} \ 
Tell whether $f^{(n)}$ defined from spin irreducible characters 
$\chi^{(\lambda,\gamma)}$ of $\widetilde{\mathfrak{S}}_n$ by 
the rightmost side of \eqref{eq:5-19} 
satisfies approximate factorization property of Definition~\ref{def:AFP} or not.

\medskip

Also, \eqref{eq:1-30} seems to hold for \eqref{eq:5-19}, where $r_{k+1}$ is the 
$(k+1)$th free cumulant of $\mathfrak{m}_{\mathrm{V}}$ given by \eqref{eq:5-16}. 

\medskip

\noindent\textbf{Problem~1b} \ 
For spin irreducible characters $\chi^{(\lambda,\gamma)}$ of 
$\widetilde{\mathfrak{S}}_n$, verify 
\begin{align}
&\lim_{n\to\infty} \frac{n^{k-1}}{|\mathcal{SP}_n|} \sum_{\lambda\in \mathcal{SP}_n} 
\frac{\chi^{(\lambda,\gamma)}([1\ 2\ \cdots\ 2k\!-\!1])}{\dim(\lambda,\gamma)} 
\notag \\ 
&= (-1)^{k+1} 6^k \!\sum_{l=1}^k \frac{(2k-1)^{l-1}}{l!} \!
\sum_{j_i\in\mathbb{N}:\, j_1+\cdots +j_l=k} 
\frac{(2^{2j_1-1}-1) B_{2j_1}}{j_1} \cdots 
\frac{(2^{2j_l-1}-1) B_{2j_l}}{j_l} 
\label{eq:5-20}
\end{align}
for $k\in\mathbb{N}, \ k\geqq 2$, where 
$B_{2j}$ is the Bernoulli number of \eqref{eq:5-4}. 

\medskip

These problems are formulated also for irreducible representations of $\mathfrak{S}_n$ 
by recalling the concentration phenomenon of \eqref{eq:1-2} and \eqref{eq:1-3} 
(and \eqref{eq:5-1}). 
Let $\chi^\lambda$ be the irreducible character of $\mathfrak{S}_n$ 
corresponding to $\lambda\in \mathcal{P}_n$. 
Cycles of $\mathfrak{S}_n$ are denoted as $(1\ 2\ \cdots \ r)$ etc. 
Since $\chi^{^t\lambda} = \mathrm{sgn}\cdot \chi^\lambda$ holds, 
we have $\chi^\lambda ((1\ 2\ \cdots \ 2k)) =0$ for $^t\lambda=\lambda$. 
This yields 
\[ 
\sum_{\lambda\in \mathcal{P}_n} \frac{\chi^\lambda ((1\ 2\ \cdots \ 2k))}
{\dim\lambda} =0, 
\qquad k\in\mathbb{N}.
\] 
Set 
\begin{equation}\label{eq:5-23}
f^{(n)}(x) = \sum_{\lambda \in \mathcal{P}_n} \frac{1}{|\mathcal{P}_n|} 
\frac{\chi^\lambda (x)}{\dim\lambda}, \qquad x\in \mathfrak{S}_n. 
\end{equation}

\noindent\textbf{Problem~2a} \ 
Tell whether $f^{(n)}$ defined from irreducible characters $\chi^\lambda$ of 
$\mathfrak{S}_n$ by \eqref{eq:5-23} 
satisfies approximate factorization property or not. 
Here in \eqref{eq:1-29}, let $\rho, \sigma\in \mathcal{P}$ (not restricted on $\mathcal{OP}$). 

\smallskip

\noindent\textbf{Problem~2b} \ 
For irreducible characters $\chi^\lambda$ of $\mathfrak{S}_n$, verify 
\[ 
\lim_{n\to\infty} \frac{n^{k-1}}{|\mathcal{P}_n|} 
\sum_{\lambda\in \mathcal{P}_n} \frac{\chi^\lambda ((1\ 2\ \cdots \ 2k\!-\!1))}
{\dim\lambda} = \text{the right hand side of \eqref{eq:5-20}}
\] 
for $k\in\mathbb{N}, \ k\geqq 2$. 

\medskip

The formulas of Problem~1b and 2b, connecting irreducible character values with 
Bernoulli numbers, might be preferably shown more directly.

\subsection{Spin characters of $\widetilde{\mathfrak{S}}_\infty$}\label{subsec:5-2}

Let $\mathcal{K}(\widetilde{\mathfrak{S}}_\infty)$ denote the set of spin normalized 
central positive-definite functions on $\widetilde{\mathfrak{S}}_\infty$. 
If $f\in \mathcal{K}(\widetilde{\mathfrak{S}}_\infty)$ is factorizable, that is, 
\begin{equation}\label{eq:5-2-1}
x, y \in \widetilde{\mathfrak{S}}_\infty \setminus \{e\}, \quad 
\mathrm{supp}\,x \cap \mathrm{supp}\,y =\varnothing \ 
\Longrightarrow \ f(xy) = f(x)f(y),  
\end{equation}
it produces a sequence which satisfies approximate factorization property of 
\eqref{eq:1-29} without error terms. 
For $f\in \mathcal{K}(\widetilde{\mathfrak{S}}_\infty)$, it satisfies 
\eqref{eq:5-2-1} if and only if it is indecomposable, that is, it is an extremal 
point of $\mathcal{K}(\widetilde{\mathfrak{S}}_\infty)$
\footnote{See \cite{Naz92}, \cite[Theorem~10.5.1]{Hir18}.}. 
Parametrization of the extremal points of 
$\mathcal{K}(\widetilde{\mathfrak{S}}_\infty)$ 
($=$ the set of spin characters of $\widetilde{\mathfrak{S}}_\infty$) 
and a character formula were given by Nazarov in \cite{Naz92} as follows
\footnote{The famous result for the infinite symmetric group $\mathfrak{S}_\infty$ 
is due to Thoma (\cite{Tho64}). 
The set $\Delta$ of \eqref{eq:5-2-2} is an analogue of the Thoma simplex. 
In \cite{HiHoHi13} and \cite{HiHo22}, we studied the spin characters of 
infinite complex reflection groups. 
In general, factorizability of a positive-definite function on an inductive limit group is 
reflected upon ergodicity of the corresponding measure on the path space of its 
branching graph. 
Such aspects for wreath product groups were investigated in detail in a series 
of our works (\cite{HiHiHo09}, \cite{HoHiHi08}, \cite{HoHi14}).}.
See also \cite{Iva99}. 
Set 
\begin{equation}\label{eq:5-2-2}
\Delta = \bigl\{ \alpha = (\alpha_i)_{i=1}^\infty \,\big|\, 
\alpha_1\geqq\alpha_2\geqq\alpha_3\geqq\cdots\geqq 0, \ 
\sum_{i=1}^\infty \alpha_i \leqq 1\bigr\}, 
\end{equation}
which parametrizes the set of spin characters of $\widetilde{\mathfrak{S}}_\infty$. 
The spin character $f_\alpha$ assigned to $\alpha\in\Delta$ is 
\begin{equation}\label{eq:5-2-3}
f_\alpha ([1\ 2\ \cdots\ k]) = 2^{-\frac{k-1}{2}} \sum_{i=1}^\infty 
\alpha_i^k, \qquad k\in\mathbb{N}, \ \text{odd}, \ \geqq 3. 
\end{equation}
Moreover, letting 
\begin{equation}\label{eq:5-2-4}
\lambda^{(n)}\in \mathcal{SP}_n, \quad \lim_{n\to\infty} \frac{\lambda^{(n)}_i}
{n} = \alpha_i \qquad (i\in\mathbb{N}),  
\end{equation}
one has 
\[ 
f_\alpha (x) = \lim_{n\to\infty} \frac{\chi^{(\lambda^{(n)},\gamma)}(x)}
{\dim(\lambda^{(n)},\gamma)}, \qquad \mathrm{type}\,x \in \mathcal{OP} 
\] 
where the convergence is taken pointwise on $\widetilde{\mathfrak{S}}_\infty$. 
When limit shapes of $Y(\lambda)$ and $S(\lambda)$ are discussed for diagrams 
of size $n$, typical length of rows and columns is of $\sqrt{n}$ order. 
Taking \eqref{eq:5-2-4} into account, we set a situation in which $\alpha_i$ 
decreases with $n$ by $1/\sqrt{n}$ order. 
In the case of $\mathfrak{S}_\infty$, Biane treated limit shapes of Young diagrams 
in this framework in \cite{Bia01}. 
Some examples were mentioned in \cite{Hor13} and \cite{Hor16}. 

For $\alpha\in\Delta$, set 
\[ 
\nu_\alpha = \bigl( 1- \sum_{i=1}^\infty \alpha_i\bigr) \delta_0 + 
\sum_{i=1}^\infty \frac{\alpha_i}{2} (\delta_{\alpha_i}+\delta_{-\alpha_i}) 
\] 
as a symmetric probability on $[-1, 1]$ and call it a Thoma measure. 
Its moments are given by 
\begin{equation}\label{eq:5-2-7}
M_{2k}(\nu_\alpha) = \sum_{i=1}^\infty \alpha_i^{2k+1}, \qquad 
M_{2k-1}(\nu_\alpha) =0, \qquad k\in\mathbb{N}.
\end{equation}
Let us investigate the condition of \eqref{eq:1-52}. 
We get from \eqref{eq:5-2-3} and \eqref{eq:5-2-7} 
\begin{multline}\label{eq:5-2-8}
n^{\frac{k-1}{2}} f_\alpha\big|_{\widetilde{\mathfrak{S}}_n} 
([1\ 2\ \cdots\ k]) \\ 
= n^{\frac{k-1}{2}} 2^{-\frac{k-1}{2}} 
\sum_{i=1}^\infty \alpha_i^k = 
\bigl( \sqrt{\frac{n}{2}}\bigr)^{k-1} M_{k-1}(\nu_\alpha) 
= M_{k-1}\Bigl( \nu_\alpha \bigl( \sqrt{\frac{2}{n}}\,\cdot\,\bigr) 
\Bigr).
\end{multline}
for $k\in\mathbb{N}$, odd, $\geqq 3$. 
Here $\nu_\alpha \bigl( \sqrt{\frac{2}{n}} dx\bigr)$ is the probability on 
$\mathbb{R}$ obtained by extending $\nu_\alpha(dx)$ by $\sqrt{n/2}$. 
If 
\begin{equation}\label{eq:5-2-9}
\alpha^{(n)} = (\alpha^{(n)}_i)_{i=1}^\infty \in \Delta, \qquad 
\alpha^{(n)}_1 = O(1/\sqrt{n}) 
\end{equation}
is assumed, then 
\[ 
\mathrm{supp}\,\nu_{\alpha^{(n)}} \bigl( \sqrt{\frac{2}{n}}\,\cdot\,\bigr) 
\subset \bigl[ -\sqrt{\frac{n}{2}} \alpha^{(n)}_1, 
\sqrt{\frac{n}{2}} \alpha^{(n)}_1 \bigr] 
\] 
are uniformly bounded. 
Furthermore, if the weak convergence 
\begin{equation}\label{eq:5-2-11}
\nu_{\alpha^{(n)}} \bigl( \sqrt{\frac{2}{n}} \,\cdot\,\bigr) 
\ \xrightarrow[\,n\to\infty\,] \ \nu 
\end{equation}
is assumed, then $\mathrm{supp}\,\nu$ is compact, and \eqref{eq:5-2-11} 
implies convergence of the moments, too. 
Hence \eqref{eq:5-2-8} yields 
\begin{equation}\label{eq:5-2-12}
n^{\frac{k-1}{2}} f_{\alpha^{(n)}}\big|_{\widetilde{\mathfrak{S}}_n} 
\bigl( [1\ 2\ \cdots\ k]\bigr) \ \xrightarrow[\,n\to\infty\,] \ 
M_{k-1}(\nu)
\end{equation}
for $k\in\mathbb{N}$, odd, $\geqq 3$. 
It is obvious that the sequence 
$\{ f_{\alpha^{(n)}}\big|_{\widetilde{\mathfrak{S}}_n}\}$ 
satisfies approximate factorization property \eqref{eq:1-29}. 
Comparing \eqref{eq:5-2-12} with \eqref{eq:1-52}, we have 
\begin{equation}\label{eq:5-2-13}
r_{k+1} = M_{k-1}(\nu). 
\end{equation}
Since $\mathrm{supp}\,\nu$ is compact, \eqref{eq:1-54} is satisfied. 
The situation of Theorem~\ref{th:dynamic} (1) is thus realized under the 
condition \eqref{eq:1-53} for pausing time. 
Hence we capture limit shape $\omega_t$, taking as an initial condition the 
sequence of probabilities on $(\widetilde{\mathfrak{S}}_n)^\wedge_{\mathrm{spin}}$ 
corresponding to 
$\{ f_{\alpha^{(n)}}\big|_{\widetilde{\mathfrak{S}}_n}\}$. 
The $R$ transform of Voiculescu of $\mathfrak{m}_{\omega_t}$ is 
computed from \eqref{eq:1-57}. 
Noting $R_k(\mathfrak{m}_{\omega_0}) = r_k$ and 
term-wise integration is valid, we get from \eqref{eq:5-2-13} 
\begin{align}
R_{\mathfrak{m}_{\omega_t}}(\zeta) 
&= \sum_{k=1}^\infty R_k(\mathfrak{m}_{\omega_t}) \zeta^{k-1} 
= \zeta + \sum_{j=2}^\infty r_{2j} \, e^{-\frac{(2j-1)t}{m}} \zeta^{2j-1} 
= \zeta + \sum_{j=2}^\infty (e^{-\frac{t}{m}}\zeta)^{2j-1}M_{2j-2}(\nu) \notag \\ 
&= \zeta +e^{-\frac{t}{m}}\zeta \int_{\mathbb{R}} 
\sum_{j=2}^\infty (e^{-\frac{t}{m}}\zeta x)^{2j-2}\nu(dx) 
=\int_{\mathbb{R}} \frac{\zeta\{ 1- (1-e^{-\frac{t}{m}})
(e^{-\frac{t}{m}}\zeta x)^2\}}{1- ( e^{-\frac{t}{m}}\zeta x)^2} \,\nu(dx)
\label{eq:5-2-14}
\end{align}
if $|\zeta|$ is small enough to yield $|\zeta| < 1/a$ with $\mathrm{supp}\,\nu \subset [-a, a]$. 
Let us see some examples of the parameter sequences $\{\alpha^{(n)}\}$ 
which satisfy \eqref{eq:5-2-9} and \eqref{eq:5-2-11}. 

\medskip

\textbf{1.} If $\alpha^{(n)}\equiv 0$, then 
$\nu_{\alpha^{(n)}} = \delta_0$ and 
$R_{\mathfrak{m}_{\omega_0}} = \zeta$. 
Hence $\omega_0$ agrees with $\varOmega_{\mathrm{VKLS}}$. 
This is a trivial case of no time evolution. 

\textbf{2.} Let $\alpha^{(n)}$ obey a uniform distribution, that is, 
\[ 
\alpha^{(n)} = \bigl( \frac{1}{N}, \cdots, \frac{1}{N}, 0, 0, \cdots \bigr) 
\qquad (\text{the first $N$ entries are $1/N$}).
\] 
We assume 
\[ 
\lim_{n\to\infty} \frac{\sqrt{n}}{\sqrt{2}N} = c >0 
\] 
to be adapted to \eqref{eq:5-2-9}
\footnote{The case of \textbf{2.} is treated in \cite{MaSn20} under the 
name of Schur--Weyl measures.}. 
Since \eqref{eq:5-2-3} and \eqref{eq:5-2-12} yield 
\[ 
M_{k-1}(\nu) = \lim_{n\to\infty} n^{\frac{k-1}{2}} 
2^{-\frac{k-1}{2}} \frac{N}{N^k} = c^{k-1} 
\] 
for $k\in\mathbb{N}$, odd, $\geqq 3$, we have \ 
%
$\nu = (\delta_c+\delta_{-c})/2$.
%
The $R$ transform 
$R_{\mathfrak{m}_{\omega_t}}(\zeta)$ is obtained from \eqref{eq:5-2-14}. 
In particular, we get 
\begin{equation}\label{eq:5-2-19}
R_{\mathfrak{m}_{\omega_0}}(\zeta) = \frac{\zeta}{1- c^2\zeta^2} 
\end{equation}
for $\mathfrak{m}_{\omega_0}$. 
Since the Stieltjes transform of $\mathfrak{m}_{\omega_0}$, 
$\zeta = G_{\mathfrak{m}_{\omega_0}}(z)$, is the inverse of  
$\zeta^{-1}+ R_{\mathfrak{m}_{\omega_0}}(\zeta)$, \eqref{eq:5-2-19} 
yields 
\begin{equation}\label{eq:5-2-20}
\frac{1}{\zeta} + \frac{\zeta}{1- c^2\zeta^2} =z, \qquad 
\text{hence} \quad c^2z\zeta^3 + (1-c^2)\zeta^2 -z\zeta +1 =0.
\end{equation}
Setting the real and imaginary parts as $z=x+iy$, $\zeta =\xi+i\eta$ in 
\eqref{eq:5-2-20}, we get the density of the absolutely continuous part 
of $\mathfrak{m}_{\omega_0}$ by 
\begin{equation}\label{eq:5-2-21}
\lim_{y\searrow 0} (-\frac{\eta}{\pi}). 
\end{equation}
Let us give a heuristic computation in the case of $c=1$. 
Putting $c=1$ in \eqref{eq:5-2-20} with \eqref{eq:5-2-21}, we have 
\begin{equation}\label{eq:5-2-22} 
x\xi^3- 3x\xi\eta^2-x\xi+1=0, \qquad 
3x\xi^2\eta-x\eta- x\eta^3 =0. 
\end{equation}
For $\eta\neq 0$, \eqref{eq:5-2-22} yields 
\[ 
x = \pm \frac{3\sqrt{3}}{2} \frac{1}{(4\eta^2+1)\sqrt{\eta^2+1}}, 
\qquad x\neq 0.
\] 
Taking \eqref{eq:5-2-21} into account, let $-\eta/\pi =u$ be a new variable to have 
\begin{equation}\label{eq:5-2-24}
x = \pm \frac{3\sqrt{3}}{2} \frac{1}{(4\pi^2u^2+1)\sqrt{\pi^2u^2+1}}, 
\qquad u\geqq 0. 
\end{equation}
Then we obtain the density from \eqref{eq:5-2-24} as 
\begin{equation}\label{eq:5-2-25}
u = u(x), \qquad x\in [-3\sqrt{3}/2, 3\sqrt{3}/2] .
\end{equation}
Indeed, $\mathfrak{m}_{\omega_0}$ is absolutely continuous because of 
\[ 
\int_{-3\sqrt{3}/2}^{3\sqrt{3}/2} u(x) dx =1.
\] 

\begin{figure}[hbt]
\centering
\input{SW1-fig}
\caption{The curve of \eqref{eq:5-2-25}}
\end{figure}

\textbf{3.} Let $\alpha^{(n)}$ obey a geometric distribution, that is, 
\[ 
\alpha^{(n)}_i = (1-q) q^{i-1},\qquad i\in\mathbb{N}
\] 
for $0<q\leqq 1$. 
We assume 
\[ 
\lim_{n\to\infty} \sqrt{\frac{n}{2}} (1-q) = r >0 
\] 
to be adapted to \eqref{eq:5-2-9}. 
Since \eqref{eq:5-2-3} and \eqref{eq:5-2-12} yield 
\[ 
M_{k-1}(\nu) = \lim_{n\to\infty} n^{\frac{k-1}{2}} 
2^{-\frac{k-1}{2}} \sum_{i=1}^\infty (1-q)^k q^{(i-1)k} 
= \frac{r^{k-1}}{k} = \int_{-r}^r \frac{x^{k-1}}{2r} dx 
\] 
for $k\in\mathbb{N}$, odd, $\geqq 3$, we see $\nu$ is uniformly distributed 
on $[-r, r]$. 
Computing the $R$ transform from \eqref{eq:5-2-14}, we have 
\begin{equation}\label{eq:5-2-30}
R_{\mathfrak{m}_{\omega_t}}(\zeta) = (1-e^{-\frac{t}{m}})\zeta 
+\frac{1}{2r} \log \frac{1+ r e^{-\frac{t}{m}}\zeta}{1- r e^{-\frac{t}{m}}\zeta}. 
\end{equation}
The Stieltjes transform of $\mathfrak{m}_{\omega_t}$, 
$\zeta = G_{\mathfrak{m}_{\omega_t}}(z)$, satisfies 
$\zeta^{-1}+ R_{\mathfrak{m}_{\omega_t}}(\zeta) = z$. 
It is an interesting, but difficult, problem to write down the $t$-dependence 
of $\mathfrak{m}_{\omega_t}$ from \eqref{eq:5-2-30} explicitly. 

We note that the examples of \textbf{2.} and \textbf{3.} are obtained by 
substituting $\alpha\in\Delta$ of \eqref{eq:5-2-2} with probability distributions on $\mathbb{N}$. 
Although $\alpha_i$ is taken to be decreasing, monotonicity is not essential 
because the character value of \eqref{eq:5-2-3} is a symmetric function in $\alpha_i$.

\section{Appendix}

\subsection{Transform of Vershik curve}\label{subsec:A1}

First it is easy to see $\psi_{\mathrm{US}}$ of \eqref{eq:1-5} is transformed to the right 
half of $\psi_{\mathrm{U}}$ of \eqref{eq:1-5} by \eqref{eq:1-6}. 
In fact, putting \eqref{eq:1-6}, that is, 
\begin{equation}\label{eq:a1-1}
y = \sqrt{2} v, \qquad x = \sqrt{2}(u-v)
\end{equation}
into $y = (2\sqrt{3}/\pi) \log (1+e^{-\pi x/ 2\sqrt{3}})$, 
we have 
\[ 
e^{-\pi u/\sqrt{6}}+ e^{-\pi v/\sqrt{6}} =1, \qquad \text{hence} \quad  
v = -\frac{\sqrt{6}}{\pi} \log (1- e^{-\frac{\pi}{\sqrt{6}}u}).
\] 
Next let us see PDE of \eqref{eq:1-12} is transformed to PDE of \eqref{eq:1-11} 
by \eqref{eq:1-6}. 
To avoid confusion, we use notations like 
$(\partial y/ \partial x)_t$, $(\partial y/ \partial t)_x$ for $y=y(t,x)$. 
By \eqref{eq:a1-1} we get 
\begin{multline}
\frac{dy}{\sqrt{2}} = dv = \bigl( \frac{\partial v}{\partial t}\bigr)_u dt + 
\bigl( \frac{\partial v}{\partial u}\bigr)_t du = 
\bigl( \frac{\partial v}{\partial t}\bigr)_u dt + \bigl( \frac{\partial v}{\partial u}\bigr)_t 
\bigl( \frac{dx}{\sqrt{2}}+\frac{dy}{\sqrt{2}}\bigr), 
\\ 
\text{hence} \quad \frac{1}{\sqrt{2}}\Bigl\{ 1- \bigl( \frac{\partial v}{\partial u}\bigr)_t 
\Bigr\} dy = \bigl( \frac{\partial v}{\partial t}\bigr)_u dt + 
\frac{1}{\sqrt{2}} \bigl( \frac{\partial v}{\partial u}\bigr)_t dx.  
\label{eq:a1-3}
\end{multline}
Comparing \eqref{eq:a1-3} with 
$\displaystyle dy = \bigl( \frac{\partial y}{\partial t}\bigr)_x dt + 
\bigl( \frac{\partial y}{\partial x}\bigr)_t dx$, 
we have 
\begin{equation}\label{eq:a1-4}
\bigl( \frac{\partial y}{\partial t}\bigr)_x = \sqrt{2} 
\frac{\bigl( \frac{\partial v}{\partial t}\bigr)_u}{1- \bigl( \frac{\partial v}{\partial u}\bigr)_t}, 
\qquad 
\bigl( \frac{\partial y}{\partial x}\bigr)_t = 
\frac{\bigl( \frac{\partial v}{\partial u}\bigr)_t}{1- \bigl( \frac{\partial v}{\partial u}\bigr)_t}.
\end{equation}
Differentiate once more the second of \eqref{eq:a1-4} to have 
\begin{equation}\label{eq:a1-5}
d \bigl( \frac{\partial y}{\partial x}\bigr) = 
d \Bigl( \frac{\frac{\partial v}{\partial u}}{1- \frac{\partial v}{\partial u}}\Bigr) 
= \Bigl( \frac{\partial}{\partial t} \Bigl( \frac{\frac{\partial v}{\partial u}}
{1- \frac{\partial v}{\partial u}}\Bigr) \Bigr)_u dt + 
\Bigl( \frac{\partial}{\partial u} \Bigl( \frac{\frac{\partial v}{\partial u}}
{1- \frac{\partial v}{\partial u}}\Bigr) \Bigr)_t du.
\end{equation}
Putting 
\[ 
du = \frac{dx}{\sqrt{2}}+ \frac{dy}{\sqrt{2}} = \frac{dx}{\sqrt{2}} + 
\bigl( \frac{\partial y}{\partial t}\bigr)_x \frac{dt}{\sqrt{2}} + 
\bigl( \frac{\partial y}{\partial x}\bigr)_t \frac{dx}{\sqrt{2}} 
\] 
into \eqref{eq:a1-5} and using 
\[ 
d \bigl( \frac{\partial y}{\partial x}\bigr) = 
\bigl( \frac{\partial^2 y}{\partial t \partial x}\bigr)_x dt + 
\bigl( \frac{\partial^2 y}{\partial x^2}\bigr)_t dx 
\] 
with \eqref{eq:a1-4}, we get 
\begin{equation}\label{eq:a1-8}
\frac{\partial^2 y}{\partial x^2} = \frac{1}{\sqrt{2}} \frac{\partial}{\partial u} 
\Bigl( \frac{\frac{\partial v}{\partial u}}{1- \frac{\partial v}{\partial u}} \Bigr) + 
\frac{1}{\sqrt{2}} \frac{\frac{\partial v}{\partial u}}{1- \frac{\partial v}{\partial u}} 
\cdot \frac{\partial}{\partial u} 
\Bigl( \frac{\frac{\partial v}{\partial u}}{1- \frac{\partial v}{\partial u}} \Bigr) = 
\frac{1}{\sqrt{2}} \frac{1}{1- \frac{\partial v}{\partial u}} 
\frac{\partial}{\partial u} \Bigl( \frac{\frac{\partial v}{\partial u}}
{1- \frac{\partial v}{\partial u}} \Bigr).
\end{equation}
Putting \eqref{eq:a1-8} and \eqref{eq:a1-4} into $y= y(t, x)$ of \eqref{eq:1-12}, 
we have 
\begin{align}
&\sqrt{2} \frac{\frac{\partial v}{\partial t}}{1- \frac{\partial v}{\partial u}} = 
\frac{1}{\sqrt{2}} \frac{1}{1- \frac{\partial v}{\partial u}} 
\frac{\partial}{\partial u} \Bigl( \frac{\frac{\partial v}{\partial u}}
{1- \frac{\partial v}{\partial u}} \Bigr) + \frac{\pi}{2\sqrt{3}} 
\frac{\frac{\partial v}{\partial u}}{1- \frac{\partial v}{\partial u}} \Bigl( 1+ 
\frac{\frac{\partial v}{\partial u}}{1- \frac{\partial v}{\partial u}}\Bigr), \notag 
\\ 
&\text{hence} \qquad \frac{\partial v}{\partial t} = \frac{1}{2} 
\frac{\partial}{\partial u} \Bigl( \frac{\frac{\partial v}{\partial u}}
{1- \frac{\partial v}{\partial u}} \Bigr) + \frac{1}{2} \frac{\pi}{\sqrt{6}} 
\frac{\frac{\partial v}{\partial u}}{1- \frac{\partial v}{\partial u}}. 
\label{eq:a1-9}
\end{align}
The equation \eqref{eq:a1-9} agrees with \eqref{eq:1-11} for $v =v(t,u)$, 
where $t$ is replaced by $t/2$. 
The $1/2$ multiple of time is not contradictory since it just reflects setting a constant 
in diffusive rescale of time and space.

\subsection{Kerov transition measure, continuous diagram}\label{subsec:A2}

Let us display Young diagram $Y(\lambda)$ for $\lambda\in \mathcal{P}$ in the $xy$ plane 
as in Figure~\ref{fig:1-1} right. 
Here each box has edges of length $\sqrt{2}$. 
Then, transition measure $\mathfrak{m}_{Y(\lambda)}$ is connected 
with the interlacing coordinates $x_1<y_1<\cdots <x_{r-1}<y_{r-1}<x_r$ 
consisting of valleys $x_i$ and peaks $y_i$ as 
\begin{equation}\label{eq:a2-1}
\frac{(z-y_1)\cdots (z-y_{r-1})}{(z-x_1)\cdots (z-x_r)} = 
\sum_{i=1}^r \frac{\mathfrak{m}_{Y(\lambda)}(\{x_i\})}{z-x_i} = 
\int_{\mathbb{R}} \frac{1}{z-x}\,\mathfrak{m}_{Y(\lambda)}(dx) 
\end{equation}
through partial fraction expansion. 
Transition measures are assigend to continuous diagrams also. 
Set 
\begin{align}
&\mathcal{D} = \Bigl\{ \omega : \mathbb{R} \longrightarrow \mathbb{R} \,\Big|\, 
|\omega(x)-\omega(y)| \leqq |x-y| \  (x, y\in\mathbb{R}), 
\notag \\ 
&\qquad\qquad\qquad\qquad \int_{-\infty}^{-1}(1+\omega^\prime(x))\frac{dx}{|x|} <\infty, \ 
\int_1^\infty (1-\omega^\prime(x))\frac{dx}{x} <\infty \Bigr\}, 
\notag\\ 
&\mathcal{F} = \Bigl\{ F : \mathbb{R} \longrightarrow \mathbb{R} \,\Big|\, 
\text{measurable}, \  0\leqq F(x)\leqq 1 \   (x\in\mathbb{R}), 
\notag \\ 
&\qquad\qquad\qquad\qquad \int_{-\infty}^0 \frac{F(x)}{1+|x|}\,dx <\infty, \ 
\int_0^\infty \frac{1-F(x)}{1+x}\,dx <\infty \Bigr\}, 
\notag \\ 
&\mathcal{M} = \{ \mu \text{ : probability on } \mathbb{R}\}.
\label{eq:a2-2}
\end{align}
Then, $\omega\in\mathcal{D}$, $F\in\mathcal{F}$ and $\mu\in\mathcal{M}$ 
have bijective correspondences by 
\begin{align}
&\int_{\mathbb{R}} \frac{\mu(dx)}{z-x} = \frac{1}{z} \exp \Bigl( 
- \int_{-\infty}^0 \frac{F(x)}{z-x}\,dx + \int_0^\infty \frac{1-F(x)}{z-x}\,dx \Bigr), 
\notag \\ 
&F(x) = \frac{1+\omega^\prime(x)}{2}, \quad 1-F(x) = \frac{1-\omega^\prime(x)}{2}, 
\quad \omega (x) = \int_0^x \bigl( 2 F(u)-1\bigr)\,du.
\label{eq:a2-3}
\end{align}
An element of $\mathcal{D}$ is called a continuous diagram, and an element of 
$\mathcal{F}$ a Rayleigh function. 
The correspondence 
$\mathcal{D} \leftrightarrow \mathcal{F} \leftrightarrow \mathcal{M}$ 
by \eqref{eq:a2-3} is often called the Markov transform. 
If sequence $\{\omega_n\}$ in $\mathcal{D}$ is related to sequence 
$\{\mu_n\}$ in $\mathcal{M}$ through the Markov transform, 
$\{\omega_n\}$ converges to $\omega\in\mathcal{D}$ compact-uniformly 
if and only if $\{\mu_n\}$ converges to $\mu\in\mathcal{M}$ weakly 
as $n\to\infty$. 
For details of these items about continuous diagrams, see \cite{Ker98} and \cite{Ker03}. 
The subset of $\mathcal{D}$ consisting of those $\omega$'s such that 
\[ 
\omega (x) = |x| \qquad \text{for sufficiently large } |x| 
\] 
is denoted by $\mathbb{D}$. 
Probability $\mu$ corresponding to $\omega\in\mathbb{D}$ has compact support. 
Young diagram $Y(\lambda)$ of $\lambda\in \mathcal{P}$ belongs to $\mathbb{D}$. 
Since the correspondence $Y(\lambda) \leftrightarrow \mathfrak{m}_\lambda$ 
is a special case of Markov transforms, we use $\mathfrak{m}_\omega$ to 
refer to the Markov transform of $\omega\in\mathcal{D}$. 
The VKLS curve \eqref{eq:1-10} belongs to $\mathbb{D}$, while the 
Vershik curve \eqref{eq:5-2} belongs to $\mathcal{D}\setminus\mathbb{D}$.

\subsection{Free cumulant}\label{subsec:A3}

Let $\mu$ be a probability on $\mathbb{R}$ having all moments. 
The $n$th moment $M_n = M_n(\mu)$ and the $k$th free cumulant $R_k =R_k(\mu)$ 
are connected by the free cumulant-moment formula. 
Let $NC(n)$ denote the set of noncrossing partitions of set $\{1,2,\cdots, n\}$. 
An element $\pi$ of $NC(n)$ is expressed as $\pi = \{ v_1, \cdots, v_{b(\pi)}\}$ 
by arranging the blocks of $\pi$. 
Here $b(\pi)$ is the number of blocks of $\pi$. 
For example
\footnote{The second one in Figure~\ref{fig:3-2}.}, 
$\pi = \{v_1, v_2\} = \bigl\{ \{1,2,5\}, \{3,4\}\bigr\} \in NC(5)$, $b(\pi)=2$, 
$|v_1|=3$, $|v_2|=2$. 
By using M\"{o}bius function $m(\pi)$ for $NC(n)$, set 
\begin{equation}\label{eq:a3-1}
R_n = \sum_{\pi=\{v_1, \cdots, v_{b(\pi)}\}\in NC(n)} m(\pi) 
\prod_{i=1}^{b(\pi)} M_{|v_i|}, \qquad n\in\mathbb{N}
\end{equation}
in a form of an inversion formula, which is equivalent to 
\begin{equation}\label{eq:a3-2}
M_n = \sum_{\pi=\{v_1, \cdots, v_{b(\pi)}\}\in NC(n)} 
\prod_{i=1}^{b(\pi)} R_{|v_i|}, \qquad n\in\mathbb{N}.
\end{equation}
We refer to \eqref{eq:a3-2} and \eqref{eq:a3-1} connecting the sequences 
$(M_n)$ and $(R_n)$ as the free cumulant-moment formula. 

Let us consider an algebraic probability space $(A, \phi)$ consisting of 
unital $\ast$-algebra $A$ over $\mathbb{C}$ and its state $\phi$. 
If self-adjoint element $a$ of $A$ and probability $\mu$ on $\mathbb{R}$ 
satisfy 
\[ 
\phi (a^n) = M_n(\mu), \qquad n\in\mathbb{N}, 
\] 
we write as $a\sim\mu$ and say $\mu$ is a distribution of $a$ or $a$ obeys $\mu$. 
For compactly supported probabilities $\mu$ and $\nu$ on $\mathbb{R}$, 
their free convolution $\mu\boxplus\nu$ is uniquely determined by taking 
free $a, b\in A$ and considering the distribution of $a+b$. 
Free convolution is characterized by the additive property of free cumulants: 
\begin{equation}\label{eq:a3-4}
R_k(\mu\boxplus\nu) = R_k(\mu) + R_k(\nu), \qquad k\in\mathbb{N}.
\end{equation}
Take a self-adjoint projection $q$ of $A$ such that $\phi(q)=c\neq 0$ and 
that $q$ and $a$ ($\sim \mu$) are free. 
The distribution of $qaq$ in probability space $(qAq, c^{-1}\phi|_{qAq})$ is 
uniquely determined by $\mu$ and $c$, called a free compression of $\mu$, 
and denoted by $\mu_c$. 
Free compression is characterized also by free cumulants: 
\begin{equation}\label{eq:a3-5}
R_k(\mu_c) = c^{k-1} R_k(\mu), \qquad k\in\mathbb{N}. 
\end{equation}
Even if probabilies on $\mathbb{R}$ do not necessarily have compact supports, 
if they have all moments, then free convolution and free compression can 
be defined through \eqref{eq:a3-4}, \eqref{eq:a3-5} and free cumulant-moment 
formula \eqref{eq:a3-1}. 
When moment problems are indeterminate, however, one should note that 
assigned probabilities are not uniquely determined. 

See \cite{NiSp06} for the above mentioned notions in free probability theory. 

\bigskip

\textbf{Acknowledgments} \ 
The author expresses deep appreciation to Professor Takeshi Hirai 
for instructing him in mathematics, especially representation theory, 
constantly for a long period.

\begin{figure}[hbt]
\centering
\input{fig-spin-branch}
\caption{Branching graph for symmetric groups: $(\widetilde{\mathfrak{S}}_n)^\wedge 
= (\widetilde{\mathfrak{S}}_n)^\wedge_{\mathrm{spin}}$ (left half) $\sqcup \, 
(\widetilde{\mathfrak{S}}_n)^\wedge_{\mathrm{ord}}$ (right half) with dimensions}
\label{fig:branch}
\end{figure}

\end{document}

%% file: Ylambda-fig.tex
{\unitlength 0.1in%
\begin{picture}(20.0000,4.0000)(4.0000,-8.0000)%
%
\special{pn 8}%
\special{pa 400 400}%
\special{pa 400 800}%
\special{fp}%
\special{pa 400 400}%
\special{pa 800 400}%
\special{fp}%
\special{pa 600 400}%
\special{pa 600 500}%
\special{fp}%
\special{pa 500 500}%
\special{pa 500 400}%
\special{fp}%
\special{pa 400 500}%
\special{pa 500 500}%
\special{fp}%
%
\special{pn 20}%
\special{pa 400 800}%
\special{pa 400 600}%
\special{fp}%
\special{pa 400 600}%
\special{pa 500 600}%
\special{fp}%
\special{pa 500 600}%
\special{pa 500 500}%
\special{fp}%
\special{pa 500 500}%
\special{pa 700 500}%
\special{fp}%
\special{pa 700 500}%
\special{pa 700 400}%
\special{fp}%
\special{pa 700 400}%
\special{pa 800 400}%
\special{fp}%
%
\special{pn 8}%
\special{pa 1000 400}%
\special{pa 1000 800}%
\special{fp}%
\special{pa 1000 800}%
\special{pa 1400 800}%
\special{fp}%
\special{pa 1200 800}%
\special{pa 1200 700}%
\special{fp}%
\special{pa 1100 700}%
\special{pa 1100 800}%
\special{fp}%
\special{pa 1000 700}%
\special{pa 1100 700}%
\special{fp}%
%
\special{pn 20}%
\special{pa 1000 400}%
\special{pa 1000 600}%
\special{fp}%
\special{pa 1000 600}%
\special{pa 1100 600}%
\special{fp}%
\special{pa 1100 600}%
\special{pa 1100 700}%
\special{fp}%
\special{pa 1100 700}%
\special{pa 1300 700}%
\special{fp}%
\special{pa 1300 700}%
\special{pa 1300 800}%
\special{fp}%
\special{pa 1300 800}%
\special{pa 1400 800}%
\special{fp}%
%
\special{pn 8}%
\special{pa 1600 400}%
\special{pa 2000 800}%
\special{fp}%
\special{pa 2000 800}%
\special{pa 2400 400}%
\special{fp}%
\special{pa 2100 700}%
\special{pa 2000 600}%
\special{fp}%
\special{pa 2200 600}%
\special{pa 2100 500}%
\special{fp}%
\special{pa 2000 600}%
\special{pa 1900 700}%
\special{fp}%
%
\special{pn 20}%
\special{pa 1600 400}%
\special{pa 1800 600}%
\special{fp}%
\special{pa 1800 600}%
\special{pa 1900 500}%
\special{fp}%
\special{pa 1900 500}%
\special{pa 2000 600}%
\special{fp}%
\special{pa 2000 600}%
\special{pa 2200 400}%
\special{fp}%
\special{pa 2200 400}%
\special{pa 2300 500}%
\special{fp}%
\special{pa 2300 500}%
\special{pa 2400 400}%
\special{fp}%
%
\special{pn 8}%
\special{pa 2400 800}%
\special{pa 1600 800}%
\special{fp}%
\end{picture}}%

%% file: Slambda-fig.tex
{\unitlength 0.1in%
\begin{picture}(20.0000,4.0000)(4.0000,-8.0000)%
%
\special{pn 8}%
\special{pa 400 400}%
\special{pa 400 800}%
\special{fp}%
\special{pa 400 400}%
\special{pa 800 400}%
\special{fp}%
\special{pa 700 400}%
\special{pa 700 500}%
\special{fp}%
\special{pa 700 500}%
\special{pa 400 500}%
\special{fp}%
\special{pa 500 600}%
\special{pa 500 400}%
\special{fp}%
\special{pa 600 400}%
\special{pa 600 600}%
\special{fp}%
\special{pa 600 600}%
\special{pa 500 600}%
\special{fp}%
%
\special{pn 8}%
\special{pa 1000 400}%
\special{pa 1000 800}%
\special{fp}%
\special{pa 1000 800}%
\special{pa 1400 800}%
\special{fp}%
\special{pa 1300 800}%
\special{pa 1300 700}%
\special{fp}%
\special{pa 1300 700}%
\special{pa 1000 700}%
\special{fp}%
\special{pa 1100 800}%
\special{pa 1100 600}%
\special{fp}%
\special{pa 1100 600}%
\special{pa 1200 600}%
\special{fp}%
\special{pa 1200 600}%
\special{pa 1200 800}%
\special{fp}%
%
\special{pn 8}%
\special{pa 1600 400}%
\special{pa 2000 800}%
\special{fp}%
\special{pa 2000 800}%
\special{pa 2400 400}%
\special{fp}%
\special{pa 2300 500}%
\special{pa 2200 400}%
\special{fp}%
\special{pa 2200 400}%
\special{pa 1900 700}%
\special{fp}%
\special{pa 2100 700}%
\special{pa 1900 500}%
\special{fp}%
\special{pa 1900 500}%
\special{pa 2000 400}%
\special{fp}%
\special{pa 2000 400}%
\special{pa 2200 600}%
\special{fp}%
\end{picture}}%

%% file: Dlambda-fig1.tex
{\unitlength 0.1in%
\begin{picture}(42.4000,11.1000)(4.0000,-15.1000)%
%
\special{pn 8}%
\special{pa 500 500}%
\special{pa 900 900}%
\special{fp}%
\special{pa 900 900}%
\special{pa 1300 500}%
\special{fp}%
\special{pa 1200 600}%
\special{pa 1100 500}%
\special{fp}%
\special{pa 1100 500}%
\special{pa 800 800}%
\special{fp}%
\special{pa 1000 800}%
\special{pa 800 600}%
\special{fp}%
\special{pa 800 600}%
\special{pa 900 500}%
\special{fp}%
\special{pa 900 500}%
\special{pa 1100 700}%
\special{fp}%
%
\special{pn 8}%
\special{pa 500 500}%
\special{pa 400 400}%
\special{fp}%
\special{pa 1300 500}%
\special{pa 1400 400}%
\special{fp}%
%
\special{pn 8}%
\special{pa 900 500}%
\special{pa 800 400}%
\special{fp}%
\special{pa 800 400}%
\special{pa 700 500}%
\special{fp}%
\special{pa 700 500}%
\special{pa 600 400}%
\special{fp}%
\special{pa 600 400}%
\special{pa 500 500}%
\special{fp}%
%
\special{pn 20}%
\special{pa 400 400}%
\special{pa 500 500}%
\special{fp}%
\special{pa 500 500}%
\special{pa 600 400}%
\special{fp}%
%
\special{pn 20}%
\special{pa 600 400}%
\special{pa 700 500}%
\special{fp}%
\special{pa 700 500}%
\special{pa 800 400}%
\special{fp}%
\special{pa 800 400}%
\special{pa 1000 600}%
\special{fp}%
\special{pa 1000 600}%
\special{pa 1100 500}%
\special{fp}%
\special{pa 1100 500}%
\special{pa 1200 600}%
\special{fp}%
\special{pa 1200 600}%
\special{pa 1400 400}%
\special{fp}%
%
\special{pn 8}%
\special{pa 1800 600}%
\special{pa 2200 1000}%
\special{fp}%
\special{pa 2200 1000}%
\special{pa 2600 600}%
\special{fp}%
\special{pa 2500 700}%
\special{pa 2400 600}%
\special{fp}%
\special{pa 2400 600}%
\special{pa 2100 900}%
\special{fp}%
\special{pa 2300 900}%
\special{pa 2100 700}%
\special{fp}%
\special{pa 2100 700}%
\special{pa 2200 600}%
\special{fp}%
\special{pa 2200 600}%
\special{pa 2400 800}%
\special{fp}%
%
\special{pn 8}%
\special{pa 1800 595}%
\special{pa 1700 495}%
\special{fp}%
\special{pa 2600 595}%
\special{pa 2700 495}%
\special{fp}%
%
\special{pn 8}%
\special{pa 2200 595}%
\special{pa 2100 495}%
\special{fp}%
\special{pa 2100 495}%
\special{pa 2000 595}%
\special{fp}%
\special{pa 2000 595}%
\special{pa 1900 495}%
\special{fp}%
\special{pa 1900 495}%
\special{pa 1800 595}%
\special{fp}%
%
\special{pn 20}%
\special{pa 1900 500}%
\special{pa 2000 600}%
\special{fp}%
\special{pa 2000 600}%
\special{pa 2100 500}%
\special{fp}%
\special{pa 2100 500}%
\special{pa 2300 700}%
\special{fp}%
\special{pa 2300 700}%
\special{pa 2400 600}%
\special{fp}%
\special{pa 2400 600}%
\special{pa 2500 700}%
\special{fp}%
\special{pa 2500 700}%
\special{pa 2700 500}%
\special{fp}%
%
\special{pn 8}%
\special{pa 4640 1510}%
\special{pa 4610 1490}%
\special{fp}%
%
\special{pn 20}%
\special{pa 1600 400}%
\special{pa 1800 600}%
\special{fp}%
\special{pa 1800 600}%
\special{pa 2000 400}%
\special{fp}%
\special{pa 2000 400}%
\special{pa 2100 500}%
\special{fp}%
\special{pa 2200 600}%
\special{pa 2300 500}%
\special{fp}%
\special{pa 2300 500}%
\special{pa 2400 600}%
\special{fp}%
%
\special{pn 8}%
\special{pa 1600 1000}%
\special{pa 2800 1000}%
\special{fp}%
%
\special{pn 20}%
\special{pa 2700 500}%
\special{pa 2800 400}%
\special{fp}%
%
\special{pn 8}%
\special{pa 2300 700}%
\special{pa 2300 1000}%
\special{dt 0.045}%
\special{pa 2000 1000}%
\special{pa 2000 600}%
\special{dt 0.045}%
\put(18.5000,-10.4000){\makebox(0,0)[lt]{$-2$}}%
\put(22.4000,-10.4000){\makebox(0,0)[lt]{$1$}}%
\end{picture}}%

%% file: content-fig.tex
{\unitlength 0.1in%
\begin{picture}(8.0000,6.0000)(4.0000,-10.0000)%
%
\special{pn 8}%
\special{pa 400 400}%
\special{pa 1200 400}%
\special{fp}%
\special{pa 400 400}%
\special{pa 400 1000}%
\special{fp}%
\special{pa 400 600}%
\special{pa 1000 600}%
\special{fp}%
\special{pa 1000 600}%
\special{pa 1000 400}%
\special{fp}%
\special{pa 800 400}%
\special{pa 800 800}%
\special{fp}%
\special{pa 800 800}%
\special{pa 600 800}%
\special{fp}%
\special{pa 600 800}%
\special{pa 600 400}%
\special{fp}%
\put(5.0000,-5.0000){\makebox(0,0){0}}%
\put(7.0000,-7.1000){\makebox(0,0){0}}%
\put(7.0000,-5.1000){\makebox(0,0){1}}%
\put(9.0000,-5.1000){\makebox(0,0){2}}%
%
\special{pn 8}%
\special{pa 1000 600}%
\special{pa 1000 800}%
\special{fp}%
\special{pa 1000 800}%
\special{pa 800 800}%
\special{fp}%
\put(8.9000,-7.0000){\makebox(0,0){1}}%
\end{picture}}%

%% file: hook-length-fig.tex
{\unitlength 0.1in%
\begin{picture}(12.0000,6.6000)(4.0000,-10.6000)%
%
\special{pn 8}%
\special{pa 500 600}%
\special{pa 600 500}%
\special{fp}%
\special{pa 500 600}%
\special{pa 800 900}%
\special{fp}%
\special{pa 800 900}%
\special{pa 1200 500}%
\special{fp}%
\special{pa 1200 500}%
\special{pa 1100 400}%
\special{fp}%
\special{pa 1100 400}%
\special{pa 700 800}%
\special{fp}%
\special{pa 600 500}%
\special{pa 900 800}%
\special{fp}%
\special{pa 400 1000}%
\special{pa 1300 1000}%
\special{fp}%
%
\special{pn 8}%
\special{pa 500 1000}%
\special{pa 500 600}%
\special{dt 0.045}%
\special{pa 1200 500}%
\special{pa 1200 1000}%
\special{dt 0.045}%
\put(4.8000,-10.5000){\makebox(0,0)[lt]{$\alpha$}}%
\put(11.8000,-10.6000){\makebox(0,0)[lt]{$\beta$}}%
\put(16.0000,-6.0000){\makebox(0,0)[lt]{$h(\ast) = \beta - \alpha -1$}}%
\put(8.0000,-8.0000){\makebox(0,0){$\ast$}}%
\end{picture}}%

%% file: zoneI-III-fig.tex
{\unitlength 0.1in%
\begin{picture}(32.0000,24.5000)(4.0000,-28.5000)%
%
\special{pn 8}%
\special{pa 400 2800}%
\special{pa 3600 2800}%
\special{fp}%
\special{pa 2000 2800}%
\special{pa 3600 1200}%
\special{fp}%
\special{pa 400 1200}%
\special{pa 2000 2800}%
\special{fp}%
\special{pa 2000 2800}%
\special{pa 2000 400}%
\special{fp}%
\special{pa 3400 1400}%
\special{pa 2800 800}%
\special{fp}%
\special{pa 1800 2600}%
\special{pa 3200 1200}%
\special{fp}%
\special{pa 3000 1000}%
\special{pa 1800 2200}%
\special{fp}%
\special{pa 1800 2200}%
\special{pa 2200 2600}%
\special{fp}%
\special{pa 1800 1800}%
\special{pa 2400 2400}%
\special{fp}%
\special{pa 1800 1800}%
\special{pa 2800 800}%
\special{fp}%
\special{pa 2000 1600}%
\special{pa 1800 1400}%
\special{fp}%
\special{pa 2000 1600}%
\special{pa 2600 2200}%
\special{fp}%
\special{pa 1800 1400}%
\special{pa 2400 800}%
\special{fp}%
\special{pa 2000 1200}%
\special{pa 1800 1000}%
\special{fp}%
\special{pa 1800 1000}%
\special{pa 2000 800}%
\special{fp}%
\special{pa 2000 800}%
\special{pa 2200 1000}%
\special{fp}%
\special{pa 2200 1000}%
\special{pa 3000 1800}%
\special{fp}%
\special{pa 2800 2000}%
\special{pa 2000 1200}%
\special{fp}%
\special{pa 2400 800}%
\special{pa 3200 1600}%
\special{fp}%
\put(26.8000,-9.3000){\makebox(0,0)[lt]{$\mu/\lambda$}}%
%
\special{pn 20}%
\special{pa 400 1200}%
\special{pa 800 800}%
\special{fp}%
\special{pa 800 800}%
\special{pa 1000 1000}%
\special{fp}%
\special{pa 1000 1000}%
\special{pa 600 1400}%
\special{fp}%
\special{pa 600 1400}%
\special{pa 400 1200}%
\special{fp}%
\special{pa 1000 1000}%
\special{pa 1200 800}%
\special{fp}%
\special{pa 1200 800}%
\special{pa 2000 1600}%
\special{fp}%
\special{pa 2000 1600}%
\special{pa 1800 1800}%
\special{fp}%
\special{pa 1800 1800}%
\special{pa 1000 1000}%
\special{fp}%
\special{pa 2000 2000}%
\special{pa 1800 1800}%
\special{fp}%
\special{pa 1800 1800}%
\special{pa 1400 2200}%
\special{fp}%
\special{pa 1400 2200}%
\special{pa 1600 2400}%
\special{fp}%
\special{pa 1600 2400}%
\special{pa 2000 2000}%
\special{fp}%
%
\special{pn 8}%
\special{pa 2000 800}%
\special{pa 1800 600}%
\special{fp}%
\special{pa 1800 600}%
\special{pa 1600 800}%
\special{fp}%
\special{pa 1600 800}%
\special{pa 1400 600}%
\special{fp}%
\special{pa 1400 600}%
\special{pa 1200 800}%
\special{fp}%
\special{pa 600 1000}%
\special{pa 600 1000}%
\special{fp}%
\special{pa 800 800}%
\special{pa 1000 600}%
\special{fp}%
\special{pa 1000 600}%
\special{pa 1200 800}%
\special{fp}%
%
\special{pn 4}%
\special{pa 1100 700}%
\special{pa 910 890}%
\special{fp}%
\special{pa 1130 730}%
\special{pa 940 920}%
\special{fp}%
\special{pa 1160 760}%
\special{pa 970 950}%
\special{fp}%
\special{pa 1070 670}%
\special{pa 880 860}%
\special{fp}%
\special{pa 1040 640}%
\special{pa 850 830}%
\special{fp}%
\special{pa 1010 610}%
\special{pa 820 800}%
\special{fp}%
%
\special{pn 8}%
\special{pa 2800 1200}%
\special{pa 2800 2800}%
\special{dt 0.045}%
\special{pa 1000 2800}%
\special{pa 1000 1000}%
\special{dt 0.045}%
\put(27.7000,-28.4000){\makebox(0,0)[lt]{$x_k$}}%
\put(10.5000,-28.4000){\makebox(0,0)[rt]{$-x_k-1$}}%
\put(19.7000,-28.5000){\makebox(0,0)[lt]{$0$}}%
\put(6.9000,-10.7000){\makebox(0,0){I}}%
\put(14.3000,-12.2000){\makebox(0,0){II}}%
\put(17.0000,-20.9000){\makebox(0,0){III}}%
\end{picture}}%

%% file: zoneVI-VIII-fig.tex
{\unitlength 0.1in%
\begin{picture}(60.0000,24.5000)(4.0000,-28.5000)%
%
\special{pn 8}%
\special{pa 400 2800}%
\special{pa 3600 2800}%
\special{fp}%
\special{pa 2000 2800}%
\special{pa 3600 1200}%
\special{fp}%
\special{pa 400 1200}%
\special{pa 2000 2800}%
\special{fp}%
\special{pa 2000 2800}%
\special{pa 2000 400}%
\special{fp}%
\special{pa 3400 1400}%
\special{pa 2800 800}%
\special{fp}%
\special{pa 1800 2600}%
\special{pa 3200 1200}%
\special{fp}%
\special{pa 3000 1000}%
\special{pa 1800 2200}%
\special{fp}%
\special{pa 1800 2200}%
\special{pa 2200 2600}%
\special{fp}%
\special{pa 1800 1800}%
\special{pa 2400 2400}%
\special{fp}%
\special{pa 1800 1800}%
\special{pa 2800 800}%
\special{fp}%
\special{pa 2000 1600}%
\special{pa 1800 1400}%
\special{fp}%
\special{pa 2000 1600}%
\special{pa 2600 2200}%
\special{fp}%
\special{pa 1800 1400}%
\special{pa 2400 800}%
\special{fp}%
\special{pa 2000 1200}%
\special{pa 1800 1000}%
\special{fp}%
\special{pa 1800 1000}%
\special{pa 2000 800}%
\special{fp}%
\special{pa 2000 800}%
\special{pa 2200 1000}%
\special{fp}%
\special{pa 2200 1000}%
\special{pa 3000 1800}%
\special{fp}%
\special{pa 2800 2000}%
\special{pa 2000 1200}%
\special{fp}%
\special{pa 2400 800}%
\special{pa 3200 1600}%
\special{fp}%
\put(26.8000,-9.3000){\makebox(0,0)[lt]{$\mu/\lambda$}}%
%
\special{pn 20}%
\special{pa 400 1200}%
\special{pa 800 800}%
\special{fp}%
\special{pa 800 800}%
\special{pa 1000 1000}%
\special{fp}%
\special{pa 1000 1000}%
\special{pa 600 1400}%
\special{fp}%
\special{pa 600 1400}%
\special{pa 400 1200}%
\special{fp}%
\special{pa 1000 1000}%
\special{pa 1200 800}%
\special{fp}%
\special{pa 1200 800}%
\special{pa 2000 1600}%
\special{fp}%
\special{pa 2000 1600}%
\special{pa 1800 1800}%
\special{fp}%
\special{pa 1800 1800}%
\special{pa 1000 1000}%
\special{fp}%
\special{pa 2000 2000}%
\special{pa 1800 1800}%
\special{fp}%
\special{pa 1800 1800}%
\special{pa 1400 2200}%
\special{fp}%
\special{pa 1400 2200}%
\special{pa 1600 2400}%
\special{fp}%
\special{pa 1600 2400}%
\special{pa 2000 2000}%
\special{fp}%
%
\special{pn 4}%
\special{pa 1100 700}%
\special{pa 910 890}%
\special{fp}%
\special{pa 1130 730}%
\special{pa 940 920}%
\special{fp}%
\special{pa 1160 760}%
\special{pa 970 950}%
\special{fp}%
\special{pa 1070 670}%
\special{pa 880 860}%
\special{fp}%
\special{pa 1040 640}%
\special{pa 850 830}%
\special{fp}%
\special{pa 1010 610}%
\special{pa 820 800}%
\special{fp}%
%
\special{pn 8}%
\special{pa 2800 1200}%
\special{pa 2800 2800}%
\special{dt 0.045}%
\special{pa 1000 2800}%
\special{pa 1000 1000}%
\special{dt 0.045}%
\put(27.7000,-28.4000){\makebox(0,0)[lt]{$x_k$}}%
\put(10.5000,-28.4000){\makebox(0,0)[rt]{$-x_k-1$}}%
\put(19.7000,-28.5000){\makebox(0,0)[lt]{$0$}}%
\put(6.9000,-10.7000){\makebox(0,0){VI}}%
\put(14.3000,-12.2000){\makebox(0,0){VII}}%
\put(17.0000,-20.9000){\makebox(0,0){VIII}}%
%
\special{pn 8}%
\special{pa 2400 800}%
\special{pa 2200 600}%
\special{fp}%
\special{pa 2200 600}%
\special{pa 2000 800}%
\special{fp}%
\special{pa 2000 800}%
\special{pa 1600 400}%
\special{fp}%
\special{pa 1600 400}%
\special{pa 1200 800}%
\special{fp}%
\special{pa 1200 800}%
\special{pa 1000 600}%
\special{fp}%
\special{pa 1000 600}%
\special{pa 800 800}%
\special{fp}%
%
\special{pn 8}%
\special{pa 4000 2800}%
\special{pa 6400 2800}%
\special{fp}%
\special{pa 5200 2800}%
\special{pa 5200 1000}%
\special{fp}%
%
\special{pn 8}%
\special{pa 5200 2800}%
\special{pa 6400 1600}%
\special{fp}%
\special{pa 4000 1600}%
\special{pa 5200 2800}%
\special{fp}%
%
\special{pn 8}%
\special{pa 5200 1600}%
\special{pa 5400 1400}%
\special{fp}%
\special{pa 5400 1400}%
\special{pa 5000 1000}%
\special{fp}%
\special{pa 5000 1000}%
\special{pa 4800 1200}%
\special{fp}%
\special{pa 4800 1200}%
\special{pa 5200 1600}%
\special{fp}%
\special{pa 5000 1400}%
\special{pa 5200 1200}%
\special{fp}%
%
\special{pn 20}%
\special{pa 4800 1200}%
\special{pa 4200 1800}%
\special{fp}%
\special{pa 4200 1800}%
\special{pa 4600 2200}%
\special{fp}%
\special{pa 4600 2200}%
\special{pa 5200 1600}%
\special{fp}%
\special{pa 5200 1600}%
\special{pa 4800 1200}%
\special{fp}%
\special{pa 5000 1400}%
\special{pa 4400 2000}%
\special{fp}%
%
\special{pn 4}%
\special{pa 5090 1090}%
\special{pa 4900 1280}%
\special{fp}%
\special{pa 5120 1120}%
\special{pa 4930 1310}%
\special{fp}%
\special{pa 5150 1150}%
\special{pa 4960 1340}%
\special{fp}%
\special{pa 5180 1180}%
\special{pa 4990 1370}%
\special{fp}%
\special{pa 5060 1060}%
\special{pa 4870 1250}%
\special{fp}%
\special{pa 5030 1030}%
\special{pa 4840 1220}%
\special{fp}%
\put(52.0000,-14.0000){\makebox(0,0){$\mu/\lambda$}}%
\put(46.2000,-15.7000){\makebox(0,0){IV}}%
\put(48.1000,-17.6000){\makebox(0,0){V}}%
\put(51.8000,-28.3000){\makebox(0,0)[lt]{$0$}}%
\end{picture}}%

%% file: sigma-maru-fig.tex
\unitlength 0.1in
\begin{picture}( 20.0000,  8.0000)(  4.0000,-12.0000)
%
\special{pn 8}%
\special{pa 400 400}%
\special{pa 400 1200}%
\special{fp}%
\special{pa 400 1200}%
\special{pa 800 1200}%
\special{fp}%
\special{pa 800 1200}%
\special{pa 800 400}%
\special{fp}%
\special{pa 400 400}%
\special{pa 1200 400}%
\special{fp}%
\special{pa 1200 400}%
\special{pa 1200 600}%
\special{fp}%
\special{pa 1200 600}%
\special{pa 600 600}%
\special{fp}%
\special{pa 600 400}%
\special{pa 600 800}%
\special{fp}%
\special{pa 600 800}%
\special{pa 1000 800}%
\special{fp}%
\special{pa 1000 800}%
\special{pa 1000 400}%
\special{fp}%
\special{pa 1800 400}%
\special{pa 2400 400}%
\special{fp}%
\special{pa 2400 400}%
\special{pa 2400 600}%
\special{fp}%
\special{pa 2400 600}%
\special{pa 1800 600}%
\special{fp}%
\special{pa 1800 400}%
\special{pa 1800 800}%
\special{fp}%
\special{pa 1800 800}%
\special{pa 2200 800}%
\special{fp}%
\special{pa 2200 800}%
\special{pa 2200 400}%
\special{fp}%
\special{pa 2000 400}%
\special{pa 2000 800}%
\special{fp}%
\put(10.0000,-10.0000){\makebox(0,0)[lt]{$\sigma$}}%
\put(22.0000,-10.0000){\makebox(0,0)[lt]{$\sigma^\circ$}}%
%
\special{pn 4}%
\special{pa 760 800}%
\special{pa 400 1160}%
\special{fp}%
\special{pa 800 820}%
\special{pa 420 1200}%
\special{fp}%
\special{pa 800 880}%
\special{pa 480 1200}%
\special{fp}%
\special{pa 800 940}%
\special{pa 540 1200}%
\special{fp}%
\special{pa 800 1000}%
\special{pa 600 1200}%
\special{fp}%
\special{pa 800 1060}%
\special{pa 660 1200}%
\special{fp}%
\special{pa 800 1120}%
\special{pa 720 1200}%
\special{fp}%
\special{pa 700 800}%
\special{pa 400 1100}%
\special{fp}%
\special{pa 640 800}%
\special{pa 400 1040}%
\special{fp}%
\special{pa 600 780}%
\special{pa 400 980}%
\special{fp}%
\special{pa 600 720}%
\special{pa 400 920}%
\special{fp}%
\special{pa 600 660}%
\special{pa 400 860}%
\special{fp}%
\special{pa 600 600}%
\special{pa 400 800}%
\special{fp}%
\special{pa 600 540}%
\special{pa 400 740}%
\special{fp}%
\special{pa 600 480}%
\special{pa 400 680}%
\special{fp}%
\special{pa 600 420}%
\special{pa 400 620}%
\special{fp}%
\special{pa 560 400}%
\special{pa 400 560}%
\special{fp}%
\special{pa 500 400}%
\special{pa 400 500}%
\special{fp}%
\special{pa 440 400}%
\special{pa 400 440}%
\special{fp}%
\end{picture}%

%% file: NC32-fig.tex
\unitlength 0.1in
\begin{picture}( 42.4500,  3.6500)(  3.5500, -7.6500)
%
\special{pn 8}%
\special{pa 400 400}%
\special{pa 400 700}%
\special{fp}%
\special{pa 400 400}%
\special{pa 700 400}%
\special{fp}%
\special{pa 700 400}%
\special{pa 700 700}%
\special{fp}%
\special{pa 550 700}%
\special{pa 550 400}%
\special{fp}%
\special{pa 850 400}%
\special{pa 850 700}%
\special{fp}%
\special{pa 1000 700}%
\special{pa 1000 400}%
\special{fp}%
\special{pa 1000 400}%
\special{pa 850 400}%
\special{fp}%
\special{pa 1300 700}%
\special{pa 1300 400}%
\special{fp}%
\special{pa 1300 400}%
\special{pa 1900 400}%
\special{fp}%
\special{pa 1900 400}%
\special{pa 1900 700}%
\special{fp}%
\special{pa 1750 700}%
\special{pa 1750 550}%
\special{fp}%
\special{pa 1750 550}%
\special{pa 1600 550}%
\special{fp}%
\special{pa 1600 550}%
\special{pa 1600 700}%
\special{fp}%
\special{pa 1450 700}%
\special{pa 1450 400}%
\special{fp}%
\special{pa 2200 700}%
\special{pa 2200 400}%
\special{fp}%
\special{pa 2200 400}%
\special{pa 2800 400}%
\special{fp}%
\special{pa 2800 400}%
\special{pa 2800 700}%
\special{fp}%
\special{pa 2650 700}%
\special{pa 2650 400}%
\special{fp}%
\special{pa 2500 700}%
\special{pa 2500 550}%
\special{fp}%
\special{pa 2500 550}%
\special{pa 2350 550}%
\special{fp}%
\special{pa 2350 550}%
\special{pa 2350 700}%
\special{fp}%
\special{pa 3100 700}%
\special{pa 3100 400}%
\special{fp}%
\special{pa 3100 400}%
\special{pa 3250 400}%
\special{fp}%
\special{pa 3250 400}%
\special{pa 3250 700}%
\special{fp}%
\special{pa 3400 700}%
\special{pa 3400 400}%
\special{fp}%
\special{pa 3400 400}%
\special{pa 3700 400}%
\special{fp}%
\special{pa 3700 400}%
\special{pa 3700 700}%
\special{fp}%
\special{pa 3550 700}%
\special{pa 3550 400}%
\special{fp}%
\special{pa 4000 700}%
\special{pa 4000 400}%
\special{fp}%
\special{pa 4600 400}%
\special{pa 4000 400}%
\special{fp}%
%
\special{pn 8}%
\special{pa 4600 400}%
\special{pa 4600 700}%
\special{fp}%
\special{pa 4450 700}%
\special{pa 4450 550}%
\special{fp}%
\special{pa 4450 550}%
\special{pa 4150 550}%
\special{fp}%
\special{pa 4150 550}%
\special{pa 4150 700}%
\special{fp}%
\special{pa 4300 700}%
\special{pa 4300 550}%
\special{fp}%
\put(4.0000,-8.5000){\makebox(0,0){1}}%
\put(5.5000,-8.5000){\makebox(0,0){2}}%
\put(7.0000,-8.5000){\makebox(0,0){3}}%
\put(8.5000,-8.5000){\makebox(0,0){4}}%
\put(10.0000,-8.5000){\makebox(0,0){5}}%
\end{picture}%

%% file: SW1-fig.tex
{\unitlength 0.1in%
\begin{picture}(26.0000,15.1000)(4.0000,-18.0000)%
\put(15.9000,-16.1000){\makebox(0,0)[rt]{O}}%
\put(15.6000,-4.0000){\makebox(0,0)[rt]{$u$}}%
\put(30.0000,-16.4000){\makebox(0,0)[rt]{$x$}}%
%
\special{pn 8}%
\special{pa 1600 1800}%
\special{pa 1600 290}%
\special{fp}%
\special{sh 1}%
\special{pa 1600 290}%
\special{pa 1580 357}%
\special{pa 1600 343}%
\special{pa 1620 357}%
\special{pa 1600 290}%
\special{fp}%
%
\special{pn 8}%
\special{pa 400 1600}%
\special{pa 3000 1600}%
\special{fp}%
\special{sh 1}%
\special{pa 3000 1600}%
\special{pa 2933 1580}%
\special{pa 2947 1600}%
\special{pa 2933 1620}%
\special{pa 3000 1600}%
\special{fp}%
\special{pn 8}%
\special{pa 2600 1600}%
\special{pa 2586 1572}%
\special{pa 2561 1552}%
\special{pa 2533 1537}%
\special{pa 2504 1524}%
\special{pa 2474 1511}%
\special{pa 2354 1467}%
\special{pa 2324 1457}%
\special{pa 2293 1446}%
\special{pa 2263 1436}%
\special{pa 2173 1403}%
\special{pa 2143 1391}%
\special{pa 2114 1379}%
\special{pa 2084 1367}%
\special{pa 2055 1354}%
\special{pa 1997 1326}%
\special{pa 1969 1310}%
\special{pa 1942 1294}%
\special{pa 1915 1277}%
\special{pa 1888 1259}%
\special{pa 1863 1239}%
\special{pa 1839 1219}%
\special{pa 1815 1197}%
\special{pa 1794 1173}%
\special{pa 1774 1148}%
\special{pa 1755 1122}%
\special{pa 1738 1095}%
\special{pa 1723 1067}%
\special{pa 1710 1038}%
\special{pa 1697 1008}%
\special{pa 1687 978}%
\special{pa 1678 947}%
\special{pa 1662 885}%
\special{pa 1655 854}%
\special{pa 1650 823}%
\special{pa 1642 759}%
\special{pa 1636 728}%
\special{pa 1631 696}%
\special{pa 1629 664}%
\special{pa 1629 632}%
\special{pa 1626 600}%
\special{pa 1620 504}%
\special{pa 1618 440}%
\special{pa 1616 409}%
\special{pa 1616 400}%
\special{fp}%
\special{pn 8}%
\special{pa 600 1600}%
\special{pa 614 1572}%
\special{pa 639 1552}%
\special{pa 667 1537}%
\special{pa 696 1524}%
\special{pa 726 1511}%
\special{pa 846 1467}%
\special{pa 876 1457}%
\special{pa 907 1446}%
\special{pa 937 1436}%
\special{pa 1027 1403}%
\special{pa 1057 1391}%
\special{pa 1086 1379}%
\special{pa 1116 1367}%
\special{pa 1145 1354}%
\special{pa 1203 1326}%
\special{pa 1231 1310}%
\special{pa 1258 1294}%
\special{pa 1285 1277}%
\special{pa 1312 1259}%
\special{pa 1337 1239}%
\special{pa 1361 1219}%
\special{pa 1385 1197}%
\special{pa 1406 1173}%
\special{pa 1426 1148}%
\special{pa 1445 1122}%
\special{pa 1462 1095}%
\special{pa 1477 1067}%
\special{pa 1490 1038}%
\special{pa 1503 1008}%
\special{pa 1513 978}%
\special{pa 1522 947}%
\special{pa 1538 885}%
\special{pa 1545 854}%
\special{pa 1550 823}%
\special{pa 1558 759}%
\special{pa 1564 728}%
\special{pa 1569 696}%
\special{pa 1571 664}%
\special{pa 1571 632}%
\special{pa 1574 600}%
\special{pa 1580 504}%
\special{pa 1582 440}%
\special{pa 1584 409}%
\special{pa 1584 400}%
\special{fp}%
\end{picture}}%

%% file: ms-2024nov13.bbl
\begin{thebibliography}{99}

\bibitem{Bia95} P. Biane, 
Permutation model for semi-circular systems and quantum random walks, 
Pacific J. Math. 171 (1995), 373--387.

\bibitem{Bia97} P. Biane, 
On the free convolution with a semi-circular distribution, 
Indiana Univ. Math. J. 46 (1997), 705--718.

\bibitem{Bia98} P. Biane, 
Representations of symmetric groups and free probability, 
Adv. Math. 138 (1998), 126--181.

\bibitem{Bia01} P. Biane, 
Approximate factorization and concentration for characters of 
symmetric groups, 
Internat. Math. Res. Notices (2001), 179--192.

\bibitem{Bia03} P. Biane, 
Characters of symmetric groups and free cumulants, 
In A.M. Vershik (ed.): Asymptotic Combinatorics with Applications to 
Mathematical Physics, Lect. Notes Math. Vol. 1815, pp.185--200, 
Springer, 2003.

\bibitem{Bor99} A.M. Borodin, 
Multiplicative central measures on the Schur graph, 
J. Math. Sci. Vol. 96, No.5 (1999), 3472--3477.

\bibitem{BoBuOl15} A. Borodin, A. Bufetov, G. Olshanski, 
Limit shapes for growing extreme characters of $U(\infty)$, 
Ann. Appl. Probab. 25 (2015), 2339--2381.

\bibitem{BoOl09} A. Borodin, G. Olshanski, 
Infinite-dimensional diffusions as limits of random walks on partitions, 
Probab. Theory Relat. Fields 144 (2009), 281--318.

\bibitem{DeS17} D. De Stavola, 
Asymptotic results for representations of finite groups, 
Ph.D. Thesis, Universitat Zurich, 2017.

\bibitem{DoFe16} M. Do\l\k{e}ga, V. F\'eray, 
Gaussian fluctuations of Young diagrams and structure constants of Jack characters, 
Duke Math. J. 165 (2016), 1193--1282

\bibitem{Ful04} J. Fulman, 
Stein's method, Jack measure, and the Metropolis algorithm, 
J. Combin. Theory Ser.A 108 (2004), 275--296.

\bibitem{Ful05} J. Fulman, 
Stein's method and Plancherel measure of the symmetric group, 
Trans. Amer. Math. Soc. 357 (2005), 555--570.

\bibitem{Fun16} T. Funaki, 
Lectures on Random Interfaces, 
Springer Briefs in Probability and Mathematical Statistics, Springer, 2016.

\bibitem{FuSa10} T. Funaki, M. Sasada,
Hydrodynamic limit for an evolutional  model of two-dimensional 
Young diagrams, 
Commun. Math. Phys. 299 (2010), 335--363.

\bibitem{GrRy07} I.S. Gradshteyn, I.M. Ryzhik, 
Tables of Integrals, Series, and Products, Seventh Edition, 
Academic Press, Elsevier, 2007.

\bibitem{Hir18} T. Hirai, 
Introduction to Spin Representations of Groups (in Japanese), Sugaku Shobo, 2018.

\bibitem{HiHo22} T. Hirai, A. Hora, 
Projective representations and spin characters of complex reflection groups 
$G(m,p,n)$ and $G(m,p,\infty)$, III, 
Kyoto J. Math. 62 (2022), no.1, 1--94.

\bibitem{HiHiHo09} T. Hirai, E. Hirai, A. Hora, 
Limits of characters of wreath products $\mathfrak{S}_n(T)$ of a compact group 
$T$ with the symmetric groups and characters of $\mathfrak{S}_\infty(T)$, I, 
Nagoya Math. J. 193 (2009), 1--93.

\bibitem{HiHoHi13} T. Hirai, A. Hora, E. Hirai, 
Projective representations and spin characters of complex reflection groups 
$G(m,p,n)$ and $G(m,p,\infty)$, 
MSJ Memoirs, Vol. 29, Math. Soc. Japan, 2013.

\bibitem{HoHu92} P.N. Hoffman, J.E. Humphreys, 
Projective Representations of the Symmetric Groups: 
Q-Functions and Shifted Tableaux, 
Oxford University Press, 1992.

\bibitem{Hor05} A. Hora, 
Remark on Biane's character formula and concentration phenomenon 
in asymptotic representation theory, 
In H.Heyer et al (eds.): Infinite Dimensional Harmonic Analysis III, 
pp. 141--159, World Scientific Publishing Co., 2005.

\bibitem{Hor06} A. Hora, 
Jucys-Murphy element and walks on modified Young graph, 
Banach Center Publications, Vol. 73 (2006), 223--235.

\bibitem{Hor13} A. Hora, 
Concentration phenomenon in ergodic ensembles of Young diagrams 
(in Japanese), 
Surikaisekikenkyusho (RIMS) Kokyuroku 1825 (2013), 75--90. 
http://hdl.handle.net/2433/194743 

\bibitem{Hor15} A. Hora, 
A diffusive limit for the profiles of random Young diagrams by way of 
free probability, 
Publ. RIMS Kyoto Univ. 51 (2015), 691--708.

\bibitem{Hor16} A. Hora, 
The Limit Shape Problem for Ensembles of Young Diagrams, 
Springer Briefs in Mathematical Physics, vol. 17, Springer, 2016.

\bibitem{Hor20} A. Hora, 
Effect of microscopic pausing time distributions on the dynamical 
limit shapes for random Young diagrams, 
Electron. J. Probab. 25 (2020), article no. 70, 1--21.

\bibitem{HoHi14} A. Hora, T. Hirai, 
Harmonic functions on the branching graph associated with the infinite wreath 
product of a compact group, 
Kyoto J. Math. 54 (2014), no.4, 775--817.

\bibitem{HoHiHi08} A. Hora, T. Hirai, E. Hirai, 
Limits of characters of wreath products $\mathfrak{S}_n(T)$ of a compact group 
$T$ with the symmetric groups and characters of $\mathfrak{S}_\infty(T)$, 
II: From a viewpoint of probability theory, 
J. Math. Soc. Japan 60 (2008), 1187--1217. 

\bibitem{HoOb07} A. Hora, N. Obata, 
Quantum Probability and Spectral Analysis of Graphs, 
Theoretical and Mathematical Physics, Springer, 2007.

\bibitem{Iva99} V.N. Ivanov, 
Dimensions of skew-shifted Young diagrams and projective characters 
of the infinite symmetric group,
J. Math. Sci. Vol. 96, No.5 (1999), 3517--3530.

\bibitem{Iva04} V.N. Ivanov, 
Gaussian limit for projective characters of large symmetric groups, 
J. Math. Sci. Vol. 121, No.3 (2004), 2330--2344.

\bibitem{Iva06} V. Ivanov, 
Plancherel measure on shifted Young diagrams, 
Amer. Math. Soc. Transl. (2) Vol. 217 (2006), 73--86.

\bibitem{Ker98} S. Kerov, 
Interlacing measures, 
Amer. Math. Soc. Transl. (2) Vol. 181 (1998), 35--83.

\bibitem{Ker99} S. Kerov, 
A differential model for the growth of Young diagrams, 
Amer. Math. Soc. Transl. (2) Vol. 188 (1999), 111--130.

\bibitem{Ker03} S.V. Kerov, 
Asymptotic Representation Theory of the Symmetric Group and Its 
Applications in Analysis, 
Transl. Math. Monographs Vol. 219, Amer. Math. Soc., 2003.

\bibitem{Kle05} A. Kleshchev, 
Linear and Projective Representations of Symmetric Groups, 
Cambridge University Press, 2005.

\bibitem{LoSh77} B.F. Logan, L.A. Shepp, 
A variational problem for random Young tableaux, 
Advances Math. 26 (1977), 206--222.

\bibitem{Mat18} S. Matsumoto, 
A spin analogue of Kerov polynomials, 
SIGMA, Vol. 14 (2018), 053, 1--13.

\bibitem{MaSn20} S. Matsumoto, P. \'Sniady, 
Random strict partitions and random shifted tableaux, 
Selecta Math. (2020) 26; 10, 1--59.

\bibitem{Naz90} M.L. Nazarov, 
Young's orthogonal form of irreducible projective representations of 
the symmetric group, 
J. London Math. Soc. (2) 42 (1990), 437--451.

\bibitem{Naz92} M.L. Nazarov, 
Projective representations of the infinite symmetric group, 
Advances Soviet Math. Vol. 9 (1992), 115--130.

\bibitem{NiSp06} A. Nica, R. Speicher, 
Lectures on the Combinatorics of Free Probability, 
London Math. Soc. Lect. Notes Series, vol. 335, Cambridge University Press, 2006.

\bibitem{Pet10} L. Petrov, 
Random walks on strict partitions, 
J. Math. Sci. Vol. 168, No.3 (2010), 437--463.

\bibitem{Sni06} P. \'Sniady, 
Gaussian fluctuations of characters of symmetric groups and of Young diagrams, 
Probab. Theory Relat. Fields 136 (2006), 263--297.

\bibitem{Str08} E. Strahov, 
A diffrential model for the deformation of the Plancherel growth process, 
Advances Math. 217 (2008), 2625--2663.

\bibitem{Tho64} E. Thoma, 
Die unzerlegbaren positiv-definiten Klassenfunktionen der abz\"{a}hlbar 
unendlichen, symmetrischen Gruppe, 
Math. Z. 85 (1964), 40--61.

\bibitem{Ver96} A.M. Vershik, 
Statistical mechanics of combinatorial partitions, and their limit shapes, 
Funct. Anal. Appl. 30 (1996), 90--105.

\bibitem{VeKe77} A.M. Vershik, S.V. Kerov, 
Asymptotics of the Plancherel measure of the symmetric group and the 
limiting form of Young tables, 
Soviet Math. Dokl. 18 (1977), 527--531.

\bibitem{Wei94} G.H. Weiss, 
Aspects and Applications of the Random Walk, North-Holland, 1994.

\end{thebibliography}
